# Reiteration Formulae for the Real Interpolation Method Including limiting 𝓛 or 𝓡 Spaces

L. R. Ya. Doktorski, P. Fernández-Martínez, T. Signes

ABSTRACT. We consider $K$-interpolation methods involving slowly varying functions. Let $\bar{A}^{\mathcal{L}}_{\theta,*}$ and $\bar{A}^{\mathcal{R}}_{\theta,*}$ ($0 \leq \theta \leq 1$) be the so called $\mathcal{L}$ or $\mathcal{R}$ limiting interpolation spaces which arise naturally in reiteration formulae for the limiting cases. We characterize the interpolation spaces $\left(\bar{A}^{\mathcal{L}}_{\theta_0,*}, *\right)_{\eta,r,a}$, $\left(\bar{A}^{\mathcal{R}}_{\theta_0,*}, *\right)_{\eta,r,a}$, $\left(*, \bar{A}^{\mathcal{L}}_{\theta_1,*}\right)_{\eta,r,a}$, and $\left(*, \bar{A}^{\mathcal{R}}_{\theta_1,*}\right)_{\eta,r,a}$ ($0 \leq \eta \leq 1$) for the limiting cases $\theta_0 = 0$ and $\theta_1 = 1$. This supplements the earlier papers of the authors, which only considered the case $0 < \theta_0 < \theta_1 < 1$. The proofs of most reiteration formulae are based on Holmstedt-type formulae. Applications to grand and small Lorentz spaces as well as to Lorentz–Karamata spaces are given.

## 1. Introduction

Let $\bar{A} := (A_0, A_1)$ be a compatible couple of (quasi-) Banach spaces such that $A_0 \cap A_1 \neq \{0\}$. The Peetre's $K$-functional on $A_0 + A_1$ is given by

$$K(t,f;\bar{A}) \equiv K(t,f) := \inf\left(\|f_0\|_{A_0} + t\|f_1\|_{A_1}\right) \quad (f = f_0 + f_1, f_i \in A_i, i=0,1, t>0). \tag{1}$$

The classical (Lions-Peetre) scale of interpolation spaces $\bar{A}_{\theta,q}$ is defined via the (quasi-)norms

$$\|f\|_{\theta,q} := \left\|u^{-\theta-1/q} K(u,f)\right\|_{q,(0,\infty)}, \tag{2}$$

where $\|*\|_{q,(a,b)}$ ($0<q\leq\infty$, $-\infty \leq a < b \leq \infty$) is the usual (quasi-)norm in the Lebesgue space $L_q$ on the interval $(a,b)$. For further information about basic properties of the $K$-functional and the real interpolation method we refer to e.g. [4].

The above definition requires $0<\theta<1$ or $\theta \in \{0,1\}$ for $q=\infty$, but the extreme cases $\theta \in \{0,1\}$ and $0 < q < \infty$ play an important role in certain questions in analysis. See for example [1, 2, 3, 5, 6, 12, 13, 15, 19, 20, 21] and references therein. One of the possibilities to consider these extreme cases is to involve an additional factor $a(u)$ in the formula (2):

$$\|f\|_{\theta,q;a} := \left\|u^{-\theta-1/q} a(u) K(u,f)\right\|_{q,(0,\infty)}.$$

For this factor various functions have been considered: logarithms [5, 6], broken-logarithmic functions [12], or more generally a slowly varying functions [1, 2, 13, 15, 20] (see Definition 8 below). In this way one gets a scale $\bar{A}_{\theta,q;a}$, where the extreme cases $\theta \in \{0,1\}$ for $q < \infty$ have sense.

An important property of this scale is that the reiteration spaces $\left(\bar{A}_{\theta_0,*}, \bar{A}_{\theta_1,*}\right)_{\eta,r,a}$ with $0 \leq \theta_0 < \theta_1 \leq 1$ and $0 < \eta < 1$ belong to the same scale. (Here and below "*" stands for further parameters.) However, for limiting cases $\eta \in \{0,1\}$ the reiteration spaces normally do not belong to this scale and new interpolation functors are needed to describe them. Thus, the limiting reiteration formulae lead in a natural way to a new limiting interpolation spaces $\bar{A}^{\mathcal{L}}_{\theta,*}$ and $\bar{A}^{\mathcal{R}}_{\theta,*}$. Following [13], we call them $\mathcal{R}$ und $\mathcal{L}$ spaces (see Definition 10 below). These limiting interpolation spaces occur in many fields of analysis. For example, it is shown [19, 3] that the so-called grand and small Lorentz spaces can be described in terms of the $\mathcal{R}$ and $\mathcal{L}$ spaces.

Naturally, the question arises about the description of the reiteration spaces $(*,*)_{\eta,r,a}$ ($0 \leq \eta \leq 1$) for couples where the first operand is $\bar{A}^{\mathcal{L}}_{\theta_0,*}$ or $\bar{A}^{\mathcal{R}}_{\theta_0,*}$ and/or the second operand is $\bar{A}^{\mathcal{L}}_{\theta_1,*}$ or $\bar{A}^{\mathcal{R}}_{\theta_1,*}$. In the previous papers [8, 9, 16, 17, 18], such reiteration formulae are established under condition $0 < \theta_0 < \theta_1 < 1$. Additionally, the couple $(\bar{A}^{\mathcal{R}}_{\theta,*}, \bar{A}^{\mathcal{L}}_{\theta,*})$ ($0 < \theta < 1$) is studied in the paper [11].

The principal aim of this paper is to examine the cases $\theta_0 = 0$ and $\theta_1 = 1$ and all possibilities $0 \leq \eta \leq 1$.

As in [8, 9, 11, 16, 17, 18], the present work finds its motivation in applications to the interpolation of grand and small Lorentz spaces (see Section 7 below). Other motivation was the paper [2], which gives a good technical tool to prove Holmstedt's formulae for the $K$-functional, which are needed to establish the reiteration formulae.

This paper is organized as follows. Section 3 contains some notations, definitions, and technical results. In Section 4 we collect necessary definitions and statements dealing with the real interpolation method involving slowly varying functions. The aim of Section 5 is to establish the Holmstedt's and reiteration formulae in the case when the second operand is either $\bar{A}^{\mathcal{L}}_{1,*}$ or $\bar{A}^{\mathcal{R}}_{1,*}$. In Section 6 we show how the case when the first operand is either $\bar{A}^{\mathcal{L}}_{0,*}$ or $\bar{A}^{\mathcal{R}}_{0,*}$ can be reduced to the previous case. Finally, in Section 7 we present some interpolation results for the grand and small Lorentz spaces as well as for Lorentz–Karamata spaces as applications of our general reiteration theorems.

## 2. Preliminaries

Throughout the paper, we write $X \subset Y$ for two (quasi-)normed spaces $X$ and $Y$ to indicate that $X$ is continuously embedded in $Y$. We write $X = Y$ if $X \subset Y$ and $Y \subset X$. For $f$ and $g$ being positive functions, we write $f \prec g$ if $f \leq Cg$, where the constant $C$ is independent on all significant quantities. Two positive functions $f$ and $g$ are considered equivalent ($f \approx g$) if $f \prec g$ and $g \prec f$. We write $f \uparrow$ ($f \downarrow$) if the function $f$ is non-decreasing (non-increasing). We adopt the conventions $1/\infty = 0$ and $1/0 = \infty$. The abbreviation $LHS(*)$ ($RHS(*)$) will be used for the left- (right-) hand side of the relation (*). By $\chi_{(a,b)}$ we denote the characteristic function on an interval $(a,b)$. $\|*\|_{q,(a,b)}$ is the usual (quasi-)norm in the Lebesgue space $L_q$ on the interval $(a,b)$ ($0 < q \leq \infty$, $0 \leq a < b \leq \infty$).

### 2.1. Slowly varying functions

We will need the definition and some basic properties of slowly varying functions. For further information about slowly varying functions, we refer to e.g. [13, 15, 20, 22].

**Definition 1.** We say that a positive Lebesgue measurable function $b$ is slowly varying on $(0,\infty)$, notation $b \in SV$, if, for each $\varepsilon > 0$, the function $t^{\varepsilon} b(t)$ is equivalent to an increasing function while the function $t^{-\varepsilon} b(t)$ is equivalent to a decreasing function.

**Lemma 2.** Let $b, b_1, b_2 \in SV$, $r>0$, $\alpha>0$, $0<q\leq\infty$, and $t \in (0,\infty)$.
 (i) Then $b^r \in SV$, $b\left(\frac{1}{t}\right) \in SV$, $b(t^r b_1(t)) \in SV$, and $b_1 b_2 \in SV$.
 (ii) If $f \approx g$, then $b \circ f \approx b \circ g$.
 (iii) $\left\|u^{\alpha-1/q} b(u)\right\|_{q,(0,t)} \approx t^{\alpha} b(t)$ and $\left\|u^{-\alpha-1/q} b(u)\right\|_{q,(t,\infty)} \approx t^{-\alpha} b(t)$.
 (iv) The functions $\left\|u^{-1/q} b(u)\right\|_{q,(0,t)}$ and $\left\|u^{-1/q} b(u)\right\|_{q,(t,\infty)}$ (if exist) belong to $SV$ and $b(t) \prec \left\|u^{-1/q} b(u)\right\|_{q,(0,t)}$ and $b(t) \prec \left\|u^{-1/q} b(u)\right\|_{q,(t,\infty)}$.

(v) If $\lambda \in (-\infty, \infty)$, then $\left\|u^{\lambda-1/q}b(u)\right\|_{q,(\frac{t}{2},t)} \approx \left\|u^{\lambda-1/q}b(u)\right\|_{q,(t,2t)} \approx t^\lambda b(t)$.

We will often use these properties without explicitly referencing every time. The following simple lemma will be also useful below.

**Lemma 3.** Let $a, b \in SV$ and $0<q, r\leq\infty$. The functions $\left\|s^{-1/r}b(s)\left\|u^{-1/q}a(u)\right\|_{q,(s,t)}\right\|_{r,(0,t)}$ and $\left\|s^{-1/r}b(s)\left\|u^{-1/q}a(u)\right\|_{q,(t,s)}\right\|_{r,(t,\infty)}$ (if exist) belong to $SV$ and for all $t \in (0, \infty)$
$$a(t)b(t) \prec a(t)\left\|s^{-1/r}b(s)\right\|_{r,(0,t)} \prec \left\|s^{-1/r}b(s)\left\|u^{-1/q}a(u)\right\|_{q,(s,t)}\right\|_{r,(0,t)}$$
and
$$a(t)b(t) \prec a(t)\left\|s^{-1/r}b(s)\right\|_{r,(t,\infty)} \prec \left\|s^{-1/r}b(s)\left\|u^{-1/q}a(u)\right\|_{q,(t,s)}\right\|_{r,(t,\infty)}.$$

*Proof.* We prove only the first statement. The second one can be proved similarly. We can assume that $\left\|s^{-1/r}b(s)\left\|u^{-1/q}a(u)\right\|_{q,(s,1)}\right\|_{r,(0,1)} < \infty$. Then, it is not hard to show that $\left\|s^{-1/r}b(s)\right\|_{r,(0,1)} < \infty$. Let $\varepsilon$ be an arbitrary positive number and the increasing functions $g$ and $h$ be such that $g(s) \approx s^{\varepsilon/2}b(s)$ and $h(u) \approx u^{\varepsilon/2}a(u)$. Put
$$G(t) = t^\varepsilon \left\|s^{-\varepsilon/2-1/r}g(s)\left\|u^{-\varepsilon/2-1/q}h(u)\right\|_{q,(t,s)}\right\|_{r,(t,\infty)}$$
$$\approx t^\varepsilon \left\|s^{-1/r}b(s)\left\|u^{-1/q}a(u)\right\|_{q,(t,s)}\right\|_{r,(t,\infty)}.$$
Let $0 < t_1 < t_2 < \infty$. By the change of variables $s = x\frac{t_1}{t_2}$ and $u = y\frac{t_1}{t_2}$, we get
$$G(t_1) = (t_1)^\varepsilon \left\|s^{-\varepsilon/2-1/r}g(s)\left\|u^{-\varepsilon/2-1/q}h(u)\right\|_{q,(t_1,s)}\right\|_{r,(t_1,\infty)}$$
$$= (t_1)^\varepsilon \left\|x^{-1/r}\left(x\frac{t_1}{t_2}\right)^{-\varepsilon/2} g\left(x\frac{t_1}{t_2}\right)\left\|u^{-1/q}u^{-\varepsilon/2}h(u)\right\|_{q,\left(t_1,x\frac{t_1}{t_2}\right)}\right\|_{r,(t_2,\infty)}$$
$$= (t_1)^\varepsilon \left\|x^{-1/r}\left(x\frac{t_1}{t_2}\right)^{-\varepsilon/2} g\left(x\frac{t_1}{t_2}\right)\left\|y^{-1/q}\left(y\frac{t_1}{t_2}\right)^{-\varepsilon/2} h\left(y\frac{t_1}{t_2}\right)\right\|_{q,(t_2,x)}\right\|_{r,(t_2,\infty)}$$
$$= (t_2)^\varepsilon \left\|x^{-\varepsilon/2-1/r}g\left(x\frac{t_1}{t_2}\right)\left\|y^{-\varepsilon/2-1/q}h\left(y\frac{t_1}{t_2}\right)\right\|_{q,(t_2,x)}\right\|_{r,(t_2,\infty)}$$
$$\leq (t_2)^\varepsilon \left\|x^{-\varepsilon/2-1/r}g(x)\left\|y^{-\varepsilon/2-1/q}h(y)\right\|_{q,(t_2,x)}\right\|_{r,(t_2,\infty)} = G(t_2).$$
The last inequality holds because the functions $g$ and $h$ are increasing. Thus, $t^\varepsilon \left\|s^{-1/r}b(s)\left\|u^{-1/q}a(u)\right\|_{q,(t,s)}\right\|_{r,(t,\infty)}$ is equivalent to an increasing function. Similarly can be shown that $t^{-\varepsilon} \left\|s^{-1/r}b(s)\left\|u^{-1/q}a(u)\right\|_{q,(t,s)}\right\|_{r,(t,\infty)}$ is equivalent to a decreasing function. This means that $\left\|s^{-1/r}b(s)\left\|u^{-1/q}a(u)\right\|_{q,(t,s)}\right\|_{r,(t,\infty)} \in SV$. Moreover, by Lemma 2,
$$\left\|s^{-1/r}b(s)\left\|u^{-1/q}a(u)\right\|_{q,(t,s)}\right\|_{r,(t,\infty)} \geq \left\|s^{-1/r}b(s)\left\|u^{-1/q}a(u)\right\|_{q,(t,s)}\right\|_{r,(2t,\infty)}$$
$$\geq \left\|s^{-1/r}b(s)\left\|u^{-1/q}a(u)\right\|_{q,(t,2t)}\right\|_{r,(2t,\infty)} \approx a(t)\left\|s^{-1/r}b(s)\right\|_{r,(t,\infty)} \succ a(t)b(t).$$
This completes the proof. □

**Remark 4.** (See [8, Remark 3].) Let $a, b \in SV$, $\lambda > 0$, and $\rho(t) = t^\lambda a(t)$.

(i) We may assume without loss of generality that the function $a$ is continuous.

(ii) Due to [1, Lemma 3.1], there exists a strongly increasing, differentiable function $\sigma(t) \approx \rho(t)$, such that $\sigma(+0) = 0$, $\sigma(\infty) = \infty$, and
$$\sigma'(t) \approx \sigma(t)\frac{1}{t}. \tag{3}$$

Obviously, the function $\sigma$ has inverse $\sigma^{(-1)}$ and it holds
$$\sigma^{(-1)'}(t) \approx \sigma^{(-1)}(t)\frac{1}{t}. \tag{4}$$

(iii) It is not difficult to show that for any compatible couple $\bar{A}$, it holds
$$K(\rho(t), f; \bar{A}) \approx K(\sigma(t), f; \bar{A}) \tag{5}$$
for all $f \in A_0 + A_1$ and $t \in (0, \infty)$.

(iv) Using [20, Proposition 2.2 (ii)], it can be shown that $b \circ \sigma \in SV$ and $b \circ \sigma^{(-1)} \in SV$.

(v) Let $0 \le \theta \le 1$, $0 < r \le \infty$. Consider the expression $\left\|\rho(t)^{-\theta} t^{-1/r} b(\rho(t)) \bar{K}(\rho(t), f)\right\|_{r,(0,\infty)}$. (5) implies that for all $f \in A_0 + A_1$
$$\left\|\rho(t)^{-\theta} t^{-1/r} a(\rho(t)) \bar{K}(\rho(t), f)\right\|_{r,(0,\infty)} \approx \left\|\sigma(t)^{-\theta} t^{-1/r} a(\sigma(t)) \bar{K}(\sigma(t), f)\right\|_{r,(0,\infty)}.$$
Using the change of variables $u = \sigma(t)$ and (3), we get
$$\left\|\sigma(t)^{-\theta} t^{-1/r} a(\sigma(t)) \bar{K}(\sigma(t), f)\right\|_{r,(0,\infty)} \approx \left\|u^{-\theta - 1/r} a(u) \bar{K}(u, f)\right\|_{r,(0,\infty)}.$$
Thus,
$$\left\|\rho(t)^{-\theta} t^{-1/r} b(\rho(t)) \bar{K}(\rho(t), f)\right\|_{r,(0,\infty)} \approx \left\|u^{-\theta - 1/r} a(u) \bar{K}(u, f)\right\|_{r,(0,\infty)}.$$

We will use the last property when making the change of variables.

**Lemma** 5. Let $0 < r, q, p \le \infty$, $a, b, c, \chi \in SV$, $\lambda > 0$ and $\rho(t) \approx t^\lambda \chi(t)$ be a strongly increasing, differentiable function such that $\rho(1) = 1$.

(i) If $\left\|s^{-1/r} c(s)\right\|_{r,(1,\infty)} < \infty$, $\left\|s^{-1/p} b(s)\right\|_{p,(1,\infty)} < \infty$ and $B(t) = \left\|s^{-1/p} b(s)\right\|_{p,(t,\infty)} c(\rho(t))$, then $\left\|s^{-1/r} B(s)\right\|_{r,(1,\infty)} < \infty$.

(ii) If $\left\|s^{-1/r} c(s)\right\|_{r,(0,1)} < \infty$, $\left\|\left\|s^{-1/p} b(s)\right\|u^{-1/q} a(u)\right\|_{q,(s,1)}\right\|_{p,(0,1)} < \infty$ and $B(t) = \left\|s^{-1/p} b(s)\right\|_{p,(0,t)} c(\rho(t))$, then $\left\|t^{-1/r} B(t) \left\|u^{-1/q} a(u)\right\|_{q,(t,1)}\right\|_{r,(0,1)} < \infty$.

(iii) If $\left\|s^{-1/r} c(s)\right\|_{r,(1,\infty)} < \infty$, $\left\|\left\|s^{-1/p} b(s)\right\|u^{-1/q} a(u)\right\|_{q,(1,s)}\right\|_{p,(1,\infty)} < \infty$ and $B(t) = \left\|s^{-1/p} b(s)\right\|_{p,(t,\infty)} c(\rho(t))$, then $\left\|t^{-1/r} B(t) \left\|u^{-1/q} a(u)\right\|_{q,(1,t)}\right\|_{r,(1,\infty)} < \infty$.

*Proof.* First, note that due Remark 4, such function $\rho(t)$ exists and $\rho(+0) = 0$, $\rho(\infty) = \infty$. Secondly, by the change of variables, we consider that
$$\left\|s^{-1/r} c(\rho(t))\right\|_{r,(0,1)} \approx \left\|s^{-1/r} c(s)\right\|_{r,(0,1)}$$
and
$$\left\|s^{-1/r} c(\rho(t))\right\|_{r,(1,\infty)} \approx \left\|s^{-1/r} c(s)\right\|_{r,(1,\infty)}.$$
Let us prove (i). We have
$$\left\|t^{-1/r} B(t)\right\|_{r,(1,\infty)} = \left\|t^{-1/r} \left\|s^{-1/p} b(s)\right\|_{p,(t,\infty)} c(\rho(t))\right\|_{r,(1,\infty)}$$
$$\le \left\|t^{-1/r} c(\rho(t))\right\|_{r,(1,\infty)} \left\|s^{-1/p} b(s)\right\|_{p,(1,\infty)} < \infty.$$
Now we prove (ii). In this case, we observe that by $0 < t < 1$ it holds
$$\left\|s^{-1/p} b(s)\right\|_{p,(0,t)} \left\|u^{-1/q} a(u)\right\|_{q,(t,1)} < \left\|\left\|s^{-1/p} b(s)\right\|u^{-1/q} a(u)\right\|_{q,(s,1)}\right\|_{p,(0,t)}$$

$$\leq \left\| s^{-1/p} b(s) \left\| u^{-1/q} a(u) \right\|_{q,(s,1)} \right\|_{p,(0,1)}.$$

Thus,
$$\left\| t^{-1/r} B(t) \left\| u^{-1/q} a(u) \right\|_{q,(t,1)} \right\|_{r,(0,1)}$$
$$= \left\| t^{-1/r} \left\| s^{-1/p} b(s) \right\|_{p,(0,t)} c(\rho(t)) \left\| u^{-1/q} a(u) \right\|_{q,(t,1)} \right\|_{r,(0,1)}$$
$$\prec \left\| t^{-1/r} c(\rho(t)) \right\|_{r,(0,1)} \left\| s^{-1/p} b(s) \left\| u^{-1/q} a(u) \right\|_{q,(s,1)} \right\|_{p,(0,1)} < \infty.$$

The assertion (iii) can be proved similarly. □

### 2.2. Hardy-type inequality

Using Fubini's theorem if the case $q = r$ or applying [20, Lemma 2.7] for $q < r$, we obtain the following Hardy-type inequality.

**Lemma 6.** Let $0 < q \leq r \leq \infty$, $b \in SV$, and $\mu > 0$. Then for all (Lebesgue-) measurable nonnegative functions $g$ on $(0,\infty)$

$$\left\| t^{-\mu-1/r} b(t) \left\| u^{-1/q} g(u) \right\|_{q,(0,t)} \right\|_{r,(0,\infty)} \prec \left\| t^{-\mu-1/r} b(t) g(t) \right\|_{r,(0,\infty)}$$

and

$$\left\| t^{\mu-1/r} b(t) \left\| u^{-1/q} g(u) \right\|_{q,(t,\infty)} \right\|_{r,(0,\infty)} \prec \left\| t^{\mu-1/r} b(t) g(t) \right\|_{r,(0,\infty)}.$$

## 3. Real interpolation method involving slowly varying functions: a short overview

Here we collect necessary definitions and statements dealing with the real interpolation method involving slowly varying functions. We consider a compatible couple of (quasi-) Banach spaces $\bar{A} = (A_0, A_1)$ such that $A_0 \cap A_1 \neq \{0\}$ and use the Peetre's K-functional on $A_0 + A_1$, given by (1).

**Definition 7**. As usual [4],
$$A_0 + \infty A_1 := \left\{ a \in A_0 + A_1 \colon \|a\|_{A_0 + \infty A_1} = \sup_{t > 0} K(t, a) = \lim_{t \to \infty} K(t, a) < \infty \right\}$$

and
$$A_1 + \infty A_0 := \left\{ a \in A_0 + A_1 \colon \|a\|_{A_1 + \infty A_0} = \sup_{t > 0} t^{-1} K(t, a) = \lim_{t \to 0} t^{-1} K(t, a) < \infty \right\}$$
are the Gagliardo completion of $A_0$ and $A_1$ in $A_0 + A_1$, respectively.

### 3.1. Standard interpolation spaces

**Definition 8.** [20]. Let $0 \leq \theta \leq 1$, $0 < q \leq \infty$ and $b \in SV$. We put
$$\bar{A}_{\theta,q;b} \equiv (A_0, A_1)_{\theta,q;b} := \left\{ f \in A_0 + A_1 \colon \|f\|_{\theta,q;b} = \left\| u^{-\theta-1/q} b(u) K(u, f) \right\|_{q,(0,\infty)} < \infty \right\}.$$

**Lemma 9.** [20, Proposition 2.5]. $\bar{A}_{\theta,q;b}$ is a (quasi-) Banach space. $A_0 \cap A_1 \subset \bar{A}_{\theta,q;b} \subset A_0 + A_1$ if and only if one of the following conditions is satisfied:
 (i) $0 < \theta < 1$,
 (ii) $\theta = 0$ and $\left\| u^{-1/q} b(u) \right\|_{q,(1,\infty)} < \infty$,
 (iii) $\theta = 1$ and $\left\| u^{-1/q} b(u) \right\|_{q,(0,1)} < \infty$.

Moreover, if none of these conditions holds, then $\bar{A}_{\theta,q;b} = \{0\}$.

We refer to the spaces $\bar{A}_{\theta,q;b}$ as standard interpolation spaces.

### 3.2. $\mathcal{R}$ and $\mathcal{L}$ limiting interpolation spaces

**Definition 10.** [20]. Let $0 \leq \theta \leq 1$, $0 < r, q \leq \infty$, and $a, b \in SV$. We put
$$\bar{A}^{\mathcal{L}}_{\theta,r,b,q,a} \equiv (A_0, A_1)^{\mathcal{L}}_{\theta,r,b,q,a} :=$$
$$\left\{ f \in A_0 + A_1 : \|f\|_{\mathcal{L};\theta,r,b,q,a} = \left\| t^{-1/r} b(t) \left\| u^{-\theta-1/q} a(u) K(u,f) \right\|_{q,(0,t)} \right\|_{r,(0,\infty)} < \infty \right\}$$
and
$$\bar{A}^{\mathcal{R}}_{\theta,r,b,q,a} \equiv (A_0, A_1)^{\mathcal{R}}_{\theta,r,b,q,a} :=$$
$$\left\{ f \in A_0 + A_1 : \|f\|_{\mathcal{R};\theta,r,b,q,a} = \left\| t^{-1/r} b(t) \left\| u^{-\theta-1/q} a(u) K(u,f) \right\|_{q,(t,\infty)} \right\|_{r,(0,\infty)} < \infty \right\}.$$

**Lemma 11.** [10]. Let $0 \leq \theta \leq 1$, $0 < r, q \leq \infty$, and $a, b \in SV$. Then $A_0 \cap A_1 \subset \bar{A}^{\mathcal{L}}_{\theta,r,b,q,a} \subset A_0 + A_1$ if and only if one of the following conditions is satisfied:
(i) $0 < \theta < 1$ and $\left\| t^{-1/r} b(t) \right\|_{r,(1,\infty)} < \infty$,
(ii) $\theta = 0$ and $\left\| t^{-1/r} b(t) \left\| u^{-1/q} a(u) \right\|_{q,(1,t)} \right\|_{r,(1,\infty)} < \infty$,
(iii) $\theta = 1$, $\left\| t^{-1/r} b(t) \right\|_{r,(1,\infty)} < \infty$, and $\left\| t^{-1/r} b(t) \left\| u^{-1/q} a(u) \right\|_{q,(0,t)} \right\|_{r,(0,1)} < \infty$.

Moreover, if none of these conditions holds, then $\bar{A}^{\mathcal{L}}_{\theta,r,b,q,a} = \{0\}$.

**Lemma 12.** [10]. Let $0 < r, q \leq \infty$, $0 \leq \theta \leq 1$, $a, b \in SV$. Then $A_0 \cap A_1 \subset \bar{A}^{\mathcal{R}}_{\theta,r,b,q,a} \subset A_0 + A_1$ if and only if one of the following conditions is satisfied:
(i) $0 < \theta < 1$ and $\left\| s^{-1/r} b(s) \right\|_{r,(0,1)} < \infty$.
(ii) $\theta = 0$, $\left\| s^{-1/r} b(s) \right\|_{r,(0,1)} < \infty$ and $\left\| s^{-1/r} b(s) \left\| u^{-1/q} a(u) \right\|_{q,(s,\infty)} \right\|_{r,(1,\infty)} < \infty$.
(iii) $\theta = 1$, $\left\| s^{-1/r} b(s) \left\| u^{-1/q} a(u) \right\|_{q,(s,1)} \right\|_{r,(0,1)} < \infty$.

Moreover, if none of these conditions holds, then $\bar{A}^{\mathcal{R}}_{\theta,r,b,q,a} = \{0\}$.

Similar definitions can be found for example in [1, 5, 6, 12, 13, 15] where the reader can find further properties of these spaces. These spaces arise naturally in reiteration formulae for the limiting cases. (See for example, [20, Theorems 3.2 and 4.3].) Following [13], we refer to these spaces as $\mathcal{L}$ and $\mathcal{R}$ spaces (or $\mathcal{L}$ and $\mathcal{R}$ limiting interpolation spaces). In view of the applications below, we are especially interested in $\mathcal{L}$ and $\mathcal{R}$ spaces if $a = 1$. The next lemma describes some of those in the cases $\theta = 0$ or $\theta = 1$.

**Lemma 13.** Let $0 < r \leq \infty$ and $b \in SV$.
(i) If $0 < q < \infty$ then $\bar{A}^{\mathcal{L}}_{1,r,b,q,1} = \{0\}$ and $\bar{A}^{\mathcal{R}}_{0,r,b,q,1} = \{0\}$.
(ii) If $\left\| t^{-1/r} b(t) \right\|_{r,(1,\infty)} < \infty$ then $\bar{A}^{\mathcal{L}}_{0,r,b,\infty,1} = \bar{A}_{0,r;b}$.
(iii) If $\left\| t^{-1/r} b(t) \right\|_{r,(0,1)} < \infty$ then $\bar{A}^{\mathcal{R}}_{1,r,b,\infty,1} = \bar{A}_{1,r;b}$.
(iv) If $\left\| t^{-1/r} b(t) \right\|_{r,(0,\infty)} < \infty$ then $\bar{A}^{\mathcal{L}}_{1,r,b,\infty,1} = A_1 + \infty A_0$.
(v) If $\left\| t^{-1/r} b(t) \right\|_{r,(0,\infty)} < \infty$ then $\bar{A}^{\mathcal{R}}_{0,r,b,\infty,1} = A_0 + \infty A_1$.

*Proof.* The statement (i) follows from Lemma 11 (iii) and Lemma 12 (ii). Let us prove (ii) the assertion (iii) can be proved similarly. Because $K(u,f)\uparrow$, we have
$$\|f\|_{\mathcal{L};0,r,b,\infty,1} = \left\| t^{-1/r} b(t) \| K(u,f) \|_{\infty,(0,t)} \right\|_{r,(0,\infty)} = \left\| t^{-1/r} b(t) K(t,f) \right\|_{r,(0,\infty)} = \|f\|_{0,r;b}.$$

Let us prove (iv). The statement (v) can be proved similarly. Because $u^{-1}K(u,f)\downarrow$, we see that

$$\|f\|_{\mathcal{L};1,r,b,\infty,1} = \left\|t^{-1/r}b(t)\|u^{-1}K(u,f)\|_{\infty,(0,t)}\right\|_{r,(0,\infty)}$$
$$= \left\|t^{-1/r}b(t)\right\|_{r,(0,\infty)} \sup_{u>0}(u^{-1}K(u,f)) \approx \|f\|_{A_1+\infty A_0}. \qquad \square$$

### 3.3. $\mathcal{RR}$, $\mathcal{RL}$, $\mathcal{LR}$, and $\mathcal{LL}$ extremal interpolation spaces

In the next definition, we introduce four further interpolation spaces. We follow [9, 10, 11, 16, 17, 18] where it has been shown that they appear in relation with the extreme reiteration results. We refer to these spaces as $\mathcal{LL}$, $\mathcal{LR}$, $\mathcal{RL}$, and $\mathcal{RR}$ extremal interpolation spaces.

**Definition 14.** Let $0 \leq \theta \leq 1$, $0 < s, q, r \leq \infty$, and $a, b, c \in SV$. The space $\bar{A}^{\mathcal{L},\mathcal{L}}_{\theta,s,c,r,b,q,a} \equiv (A_0, A_1)^{\mathcal{L},\mathcal{L}}_{\theta,s,c,r,b,q,a}$ is the set of all $f \in A_0 + A_1$ for which the (quasi-)norm

$$\|f\|_{\mathcal{L},\mathcal{L};\theta,s,c,r,b,q,a} = \left\|t^{-1/s}c(t)\left\|u^{-1/r}b(u)\|v^{-\theta-1/q}a(v)K(v,f)\|_{q,(0,u)}\right\|_{r,(0,t)}\right\|_{s,(0,\infty)}$$

is finite. The spaces $\bar{A}^{\mathcal{L},\mathcal{R}}_{\theta,s,c,r,b,q,a} \equiv (A_0, A_1)^{\mathcal{L},\mathcal{R}}_{\theta,s,c,r,b,q,a}$, $\bar{A}^{\mathcal{R},\mathcal{L}}_{\theta,s,c,r,b,q,a} \equiv (A_0, A_1)^{\mathcal{R},\mathcal{L}}_{\theta,s,c,r,b,q,a}$, and $\bar{A}^{\mathcal{R},\mathcal{R}}_{\theta,s,c,r,b,q,a} \equiv (A_0, A_1)^{\mathcal{R},\mathcal{R}}_{\theta,s,c,r,b,q,a}$ are defined analogously via the (quasi-)norms

$$\|f\|_{\mathcal{L},\mathcal{R};\theta,s,c,r,b,q,a} = \left\|t^{-1/s}c(t)\left\|u^{-1/r}b(u)\|v^{-\theta-1/q}a(v)K(v,f)\|_{q,(t,u)}\right\|_{r,(t,\infty)}\right\|_{s,(0,\infty)},$$

$$\|f\|_{\mathcal{R},\mathcal{L};\theta,s,c,r,b,q,a} = \left\|t^{-1/s}c(t)\left\|u^{-1/r}b(u)\|v^{-\theta-1/q}a(v)K(v,f)\|_{q,(u,t)}\right\|_{r,(0,t)}\right\|_{s,(0,\infty)},$$

and

$$\|f\|_{\mathcal{R},\mathcal{R};\theta,s,c,r,b,q,a} = \left\|t^{-1/s}c(t)\left\|u^{-1/r}b(u)\|v^{-\theta-1/q}a(v)K(v,f)\|_{q,(u,\infty)}\right\|_{r,(t,\infty)}\right\|_{s,(0,\infty)},$$

respectively.

For these spaces, it is also possible to formulate conditions under which they are trivial. For example, if $\left\|u^{-1/r}b(u)\right\|_{r,(1,\infty)} = \infty$ then $\bar{A}^{\mathcal{L},\mathcal{R}}_{\theta,s,c,r,b,q,a} = \{0\}$. We leave this to the reader. For other combinations of parameters, they are interpolation spaces for the couple $(A_0, A_1)$. Subsequently, we will assume we will assume that no spaces under consideration are trivial. The next lemma describes some extremal interpolation spaces if $q = \infty$, $a = 1$, and $\theta = 0$ or $\theta = 1$.

**Lemma 15.** Let $0 < s, r \leq \infty$, $b, c \in SV$. Then
   (i) $\bar{A}^{\mathcal{L},\mathcal{L}}_{0,s,c,r,b,\infty,1} = \bar{A}^{\mathcal{L}}_{0,s,c,r,b}$.
   (ii) $\bar{A}^{\mathcal{L},\mathcal{L}}_{1,s,c,r,b,\infty,1} = A_1 + \infty A_0$, provided $\left\|t^{-1/s}c(t)\|s^{-1/r}b(s)\|_{r,(0,t)}\right\|_{s,(0,\infty)} < \infty$.
   (iii) $\bar{A}^{\mathcal{L},\mathcal{R}}_{0,s,c,r,b,\infty,1} = \bar{A}^{\mathcal{R}}_{0,s,c,r,b}$.
   (iv) $\bar{A}^{\mathcal{L},\mathcal{R}}_{1,s,c,r,b,\infty,1} = \bar{A}_{1,s;B}$, where $B(t) := c(t)\|s^{-1/r}b(s)\|_{r,(t,\infty)}$.
   (v) $\bar{A}^{\mathcal{R},\mathcal{L}}_{0,s,c,r,b,\infty,1} = \bar{A}_{0,s;B}$, where $B(t) := c(t)\|s^{-1/r}b(s)\|_{r,(0,t)}$.
   (vi) $\bar{A}^{\mathcal{R},\mathcal{L}}_{1,s,c,r,b,\infty,1} = \bar{A}^{\mathcal{L}}_{1,s,c,r,b}$.
   (vii) $\bar{A}^{\mathcal{R},\mathcal{R}}_{0,s,c,r,b,\infty,1} = A_1 + \infty A_0$, provided
   $$\left\|t^{-1/s}c(t)\|s^{-1-1/r}b(s)K(s,f)\|_{r,(0,t)}\right\|_{s,(0,\infty)} < \infty.$$
   (viii) $\bar{A}^{\mathcal{R},\mathcal{R}}_{1,s,c,r,b,\infty,1} = \bar{A}^{\mathcal{R}}_{1,s,c,r,b}$.

*Proof.* We use the monotonicity of the $K$-functional and Definition 7.

(i)  $\|f\|_{\mathcal{L},\mathcal{L};0,s,c,r,b,\infty,1} = \left\|t^{-1/s}c(t)\left\|s^{-1/r}b(s)\|K(u,f)\|_{\infty,(0,s)}\right\|_{r,(0,t)}\right\|_{s,(0,\infty)}$
$= \left\|t^{-1/s}c(t)\left\|s^{-1/r}b(s)K(s,f)\right\|_{r,(0,t)}\right\|_{s,(0,\infty)} = \|f\|_{\mathcal{L};0,s,c,r,b}.$

(ii) $\|f\|_{\mathcal{L},\mathcal{L};1,s,c,r,b,\infty,1} = \left\|t^{-1/s}c(t)\left\|s^{-1/r}b(s)\|u^{-1}K(u,f)\|_{\infty,(0,s)}\right\|_{r,(0,t)}\right\|_{s,(0,\infty)}$
$= \left\|t^{-1/s}c(t)\left\|s^{-1/r}b(s)\right\|_{r,(0,t)}\right\|_{s,(0,\infty)} \|f\|_{A_1+\infty A_0}.$

(iii) $\|f\|_{\mathcal{L},\mathcal{R};0,s,c,r,b,\infty,1} = \left\|t^{-1/s}c(t)\left\|s^{-1/r}b(s)\|K(u,f)\|_{\infty,(t,s)}\right\|_{r,(t,\infty)}\right\|_{s,(0,\infty)}$
$= \left\|t^{-1/s}c(t)\left\|s^{-1/r}b(s)K(s,f)\right\|_{r,(t,\infty)}\right\|_{s,(0,\infty)} = \|f\|_{\mathcal{R};0,s,c,r,b}.$

(iv) $\|f\|_{\mathcal{L},\mathcal{R};1,s,c,r,b,\infty,1} = \left\|t^{-1/s}c(t)\left\|s^{-1/r}b(s)\|u^{-1}K(u,f)\|_{\infty,(t,s)}\right\|_{r,(t,\infty)}\right\|_{s,(0,\infty)}$
$= \left\|t^{-1-1/s}c(t)\left\|s^{-1/r}b(s)\right\|_{r,(t,\infty)}K(t,f)\right\|_{s,(0,\infty)} = \|f\|_{1,s;B}.$

The statements (v) – (viii) can be proved similarly. □

### 3.4. Some known formulae and auxiliary lemmas

Let $C \geq 1$, $t > 0$, and $\tilde{b}(t):= b(t^{-1})$. Under suitable conditions, the following formulae hold [4, Chap. 5, Proposition 1.2], [20], (Cf. [17, Lemma 2.12]):

$$K(t,f) \leq K(Ct,f) \leq CK(t,f),$$
$$K(t,f;A_0,A_1) = tK(t^{-1},f;A_1,A_0),$$
$$(A_0,A_1)_{\theta,q,b} = (A_1,A_0)_{1-\theta,q,\tilde{b}}, \tag{6}$$
$$(A_0,A_1)^{\mathcal{L}}_{\sigma,r,b,q,a} = (A_1,A_0)^{\mathcal{R}}_{1-\sigma,r,\tilde{b},q,\tilde{a}}, \tag{7}$$

$(A_0,A_1)^{\mathcal{L},\mathcal{L}}_{\theta,p,c,r,b,q,a} = (A_1,A_0)^{\mathcal{R},\mathcal{R}}_{1-\theta,p,\tilde{c},r,\tilde{b},q,\tilde{a}}$, $(A_0,A_1)^{\mathcal{L},\mathcal{R}}_{\theta,p,c,r,b,q,a} = (A_1,A_0)^{\mathcal{R},\mathcal{L}}_{1-\theta,p,\tilde{c},r,\tilde{b},q,\tilde{a}}$. (8)

**Lemma 16.** [9]. Let $0 < r \leq q \leq \infty$ and $a \in SV$. Then, for all $f \in A_0 + A_1$ and $t>0$
$$\left\|u^{-\theta-1/q}a(u)K(u,f)\right\|_{q,(t,\infty)} \prec \left\|u^{-\theta-1/r}a(u)K(u,f)\right\|_{r,(t,\infty)}, \text{ provided } \theta > 0$$
and
$$\left\|u^{-\theta-1/q}a(u)K(u,f)\right\|_{q,(0,t)} \prec \left\|u^{-\theta-1/r}a(u)K(u,f)\right\|_{r,(0,t)}, \text{ provided } \theta < 1.$$
Moreover, it holds
$$\left\|u^{-1/r}a(u)\right\|_{r,(t,\infty)}K(t,f) \prec \left\|u^{-1/r}a(u)K(u,f)\right\|_{r,(t,\infty)},$$
$$t^{-1}\left\|u^{-1/r}a(u)\right\|_{r,(0,t)}K(t,f) \prec \left\|u^{-1-1/r}a(u)K(u,f)\right\|_{r,(0,t)}.$$

**Corollary 17**. (Cf. [17, (2.8), (2.9)].) Let $0 < p, r \leq \infty$, $a, b \in SV$.

(i) If $0 \leq \theta < 1$ then
$$\bar{A}^{\mathcal{L}}_{\theta,p,b,r,a} \subset \bar{A}_{\theta,p;ab}.$$

(ii) If $0 < \theta \leq 1$ then
$$\bar{A}^{\mathcal{R}}_{\theta,p,b,r,a} \subset \bar{A}_{\theta,p;ab}.$$

(iii) If $\left\|u^{-1/r}a(u)\right\|_{r,(0,1)} < \infty$ and $B(t) = \left\|u^{-1/r}a(u)\right\|_{r,(0,t)}b(t)$, then
$$\bar{A}^{\mathcal{L}}_{1,p,b,r,a} \subset \bar{A}_{1,p;B}.$$

(iv) If $\left\|u^{-1/r}a(u)\right\|_{r,(1,\infty)} < \infty$ and $B(t) = \left\|u^{-1/r}a(u)\right\|_{r,(t,\infty)}b(t)$, then
$$\bar{A}^{\mathcal{R}}_{0,p,b,r,a} \subset \bar{A}_{0,p;B}.$$

*Proof.* Let us prove the first statement. The second one can be proved similarly. Because $\theta < 1$, by Lemma 16, we conclude that for all $f \in A_0 + A_1$ and $t>0$

$$t^{-\theta}a(t)K(t,f) \prec \left\|u^{-\theta-1/r}a(u)K(u,f)\right\|_{r,(0,t)}.$$

Hence,

$$\|f\|_{\theta,p;ab} = \left\|t^{-\theta-1/p}a(t)b(t)K(t,f)\right\|_{p,(0,\infty)}$$

$$\prec \left\|t^{-1/p}b(t)\left\|u^{-\theta-1/r}a(u)K(u,f)\right\|_{r,(0,t)}\right\|_{p,(0,\infty)} = \|f\|_{\mathcal{L};\theta,p,b,r,a}.$$

Now we prove the third statement. The fourth one can be proved similarly. By Lemma 16, we have

$$\|f\|_{1,p;B} = \left\|t^{-1-1/p}\left\|u^{-1/r}a(u)\right\|_{r,(0,t)}b(t)K(t,f)\right\|_{p,(0,\infty)}$$

$$\prec \left\|t^{-1/p}b(t)\left\|u^{-1-1/r}a(u)K(u,f)\right\|_{r,(0,t)}\right\|_{p,(0,\infty)} = \|f\|_{\mathcal{L};1,p,b,r,a}.$$

Thus, $\bar{A}^{\mathcal{L}}_{1,p,b,ra} \subset \bar{A}_{1,p;B}$. □

The next lemma is an analog to [16, Lemma 2.12], [17, Lemma 2.10], and [18, Lemma 2.10] and can be proved similarly.

**Lemma 18.** Let $0 \leq \theta \leq 1$, $0 < q,r \leq \infty$, $a,b \in SV$. Then, for all $f \in A_0 + A_1$ and $u>0$,

$$u^{-\theta}a(u)\left\|s^{-1/r}b(s)\right\|_{r,(0,u)}K(u,f) \prec \left\|t^{-1/r}b(t)\left\|s^{-\theta-\frac{1}{q}}a(s)K(s,f)\right\|_{q,(t,u)}\right\|_{r,(0,u)}, \quad (9)$$

$$u^{-\theta}a(u)\left\|s^{-1/r}b(s)\right\|_{r,(u,\infty)}K(u,f) \prec \left\|t^{-1/r}b(t)\left\|s^{-\theta-\frac{1}{q}}a(s)K(s,f)\right\|_{q,(u,t)}\right\|_{r,(u,\infty)}, \quad (10)$$

$$u^{-\theta}a(u)b(u)K(u,f) \prec \left\|t^{-1/r}b(t)\left\|s^{-\theta-\frac{1}{q}}a(s)K(s,f)\right\|_{q,(t,\infty)}\right\|_{r,(u,\infty)}, \quad (11)$$

$$u^{-\theta}a(u)b(u)K(u,f) \prec \left\|t^{-1/r}b(t)\left\|s^{-\theta-\frac{1}{q}}a(s)K(s,f)\right\|_{q,(0,t)}\right\|_{r,(0,u)}. \quad (12)$$

**Corollary** 19. Let $0 \leq \theta \leq 1$, $0 < p,q,r \leq \infty$, $a,b,c \in SV$. Then
(i) $\bar{A}^{\mathcal{R}}_{\theta,r,b,q,a} \subset \bar{A}_{\theta,\infty;B}$, where $B(u) = a(u)\left\|s^{-1/r}b(s)\right\|_{r,(0,u)}$,
(ii) $\bar{A}^{\mathcal{L}}_{\theta,r,b,q,a} \subset \bar{A}_{\theta,\infty;B}$, where $B(u) := a(u)\left\|s^{-1/r}b(s)\right\|_{r,(u,\infty)}$,
(iii) $\bar{A}^{\mathcal{R},\mathcal{R}}_{\theta,p,c,r,b,q,a} \subset \bar{A}_{\theta,p;abc}$,
(iv) $\bar{A}^{\mathcal{L},\mathcal{L}}_{\theta,p,c,r,b,q,a} \subset \bar{A}_{\theta,p;abc}$,
(v) $\bar{A}^{\mathcal{R},\mathcal{L}}_{\theta,p,c,r,b,q,a} \subset \bar{A}_{\theta,p;B}$, where $B(u) = c(u)a(u)\left\|s^{-1/r}b(s)\right\|_{r,(0,u)}$,
(vi) $\bar{A}^{\mathcal{L},\mathcal{R}}_{\theta,p,c,r,b,q,a} \subset \bar{A}_{\theta,p;B}$, where $B(u) = c(u)a(u)\left\|s^{-1/r}b(s)\right\|_{r,(u,\infty)}$.

*Proof.* First, we prove the embedding (i). By (9), we have

$$\|f\|_{\theta,\infty;B} = \sup_{0<u<\infty}\left(u^{-\theta}a(u)\left\|s^{-1/r}b(s)\right\|_{r,(0,u)}K(u,f)\right)$$

$$\prec \sup_{0<u<\infty}\left(\left\|t^{-1/r}b(t)\left\|s^{-\theta-\frac{1}{q}}a(s)K(s,f)\right\|_{q,(t,u)}\right\|_{r,(0,u)}\right)$$

$$= \left\|t^{-1/r}b(t)\left\|s^{-\theta-\frac{1}{q}}a(s)K(s,f)\right\|_{q,(t,\infty)}\right\|_{r,(0,\infty)} = \|f\|_{\mathcal{R};\theta,r,b,q,a}.$$

Now we prove the embedding (iii). By (11), we have

$$\|f\|_{\theta,p;abc} = \left\|u^{-\theta-1/p}a(u)b(u)c(u)K(u,f)\right\|_{p,(0,\infty)}$$

$$\prec \left\|u^{-1/p}c(u)\left\|t^{-1/r}b(t)\left\|s^{-\theta-\frac{1}{q}}a(s)K(s,f)\right\|_{q,(t,\infty)}\right\|_{r,(u,\infty)}\right\|_{p,(0,\infty)} = \|f\|_{\mathcal{R},\mathcal{R};\theta,p,c,r,b,q,a}.$$

Let us prove the embedding (v). By (9), we have
$$\|f\|_{\theta,p;B} = \left\|u^{-\theta-1/p}c(u)a(u)\left\|s^{-1/r}b(s)\right\|_{r,(0,u)}K(u,f)\right\|_{p,(0,\infty)}$$
$$\prec \left\|u^{-1/p}c(u)\left\|t^{-1/r}b(t)\left\|s^{-\theta-\frac{1}{q}}a(s)K(s,f)\right\|_{q,(t,u)}\right\|_{r,(0,u)}\right\|_{p,(0,\infty)} = \|f\|_{\mathcal{R},\mathcal{L};\theta,p,c,r,b,q,a}.$$

The embeddings (ii), (iv), and (vi) can be proved similarly. □

**Lemma 20.** Let $0 < \theta < 1$, $0 < p \leq \infty$, $0 < r < q \leq \infty$, and $a, b \in SV$. Then, for all $f \in A_0 + A_1$ and $t>0$,
$$\left\|s^{-1/q}b(s)\left\|u^{-\theta-1/p}a(u)K(u,f)\right\|_{p,(0,s)}\right\|_{q,(0,t)} \prec$$
$$\prec \left\|s^{-1/r}b(s)\left\|u^{-\theta-1/p}a(u)K(u,f)\right\|_{p,(0,s)}\right\|_{r,(0,t)} \quad (13)$$
and
$$\left\|s^{-1/q}b(s)\left\|u^{-\theta-1/p}a(u)K(u,f)\right\|_{p,(s,\infty)}\right\|_{q,(t,\infty)} \prec$$
$$\prec \left\|s^{-1/r}b(s)\left\|u^{-\theta-1/p}a(u)K(u,f)\right\|_{p,(s,\infty)}\right\|_{r,(t\infty)}.$$

*Proof.* We prove the first estimate. The second one can be proved analogously. The function $s^{1-\theta}a(s)$ is equivalent to an increasing function. Moreover, it is known [20, Theorem 3.1] that for all $f \in A_0 + A_1$ and $t>0$
$$K(t^{1-\theta}a(t), f; \bar{A}_{\theta,p,a}, A_1) \approx \left\|u^{-\theta-1/p}a(u)K(u,f)\right\|_{p,(0,t)}.$$
Therefore, the function $\frac{1}{t^{1-\theta}a(t)}\left\|u^{-\theta-1/p}a(u)K(u,f)\right\|_{p,(0,t)}$ is equivalent to a decreasing function. Hence, for $v < t$, we get
$$b(v)\left\|u^{-\theta-1/p}a(u)K(u,f)\right\|_{p,(0,v)}$$
$$\approx \left\|s^{1-\theta-1/r}a(s)b(s)\right\|_{r,(0,v)} \frac{1}{v^{1-\theta}a(v)}\left\|u^{-\theta-1/p}a(u)K(u,f)\right\|_{p,(0,v)}$$
$$\prec \left\|s^{-1/r}b(s)\left\|u^{-\theta-1/p}a(u)K(u,f)\right\|_{p,(0,s)}\right\|_{r,(0,v)}$$
$$\leq \left\|s^{-1/r}b(s)\left\|u^{-\theta-1/p}a(u)K(u,f)\right\|_{p,(0,s)}\right\|_{r,(0,t)}.$$

Hence, (13) is proved for $q = \infty$:
$$\sup_{0<v<t} b(v)\left\|u^{-\theta-1/p}a(u)K(u,f)\right\|_{p,(0,v)} \prec \left\|s^{-1/r}b(s)\left\|u^{-\theta-1/p}a(u)K(u,f)\right\|_{p,(0,s)}\right\|_{r,(0,t)}.$$

Denote $g(s) = b(s)\left\|u^{-\theta-1/p}a(u)K(u,f)\right\|_{p,(0,s)}$. If $q < \infty$ using the last estimate, we get
$$LHS(13) = \left\{\int_0^t (g(s))^q \frac{ds}{s}\right\}^{\frac{1}{q}} \leq \left(\sup_{0<s<t} g(s)\right)^{\frac{q-r}{q}} \left\{\int_0^t (g(s))^r \frac{ds}{s}\right\}^{\frac{1}{q}}$$
$$\prec \left\|s^{-1/r}b(s)\left\|u^{-\theta-1/p}a(u)K(u,f)\right\|_{p,(0,s)}\right\|_{r,(0,t)}.$$

This completes the proof. □

**Lemma 21.** (Cf. [14, Theorems 3.6 and 3.7].) Let $0 < q, r \leq \infty$, $a, b \in SV$, and $\beta > 0$. Then, for all $f \in A_0 + A_1$,
$$\left\|t^{-\beta-1/r}b(t)\left\|u^{\alpha-1/q}a(u)K(u,f)\right\|_{q,(0,t)}\right\|_{r,(0,\infty)} \approx \left\|t^{\alpha-\beta-1/r}a(t)b(t)K(t,f)\right\|_{r,(0,\infty)},$$
provided that $\alpha > -1$, and
$$\left\|t^{\beta-1/r}b(t)\left\|u^{\alpha-1/q}a(u)K(u,f)\right\|_{q,(t,\infty)}\right\|_{r,(0,\infty)} \approx \left\|t^{\alpha+\beta-1/r}a(t)b(t)K(t,f)\right\|_{r,(0,\infty)},$$

provided that $\alpha < 0$.

*Proof.* We prove the first equivalence. The second one can be prove analogously. We have
$$\left\|u^{\alpha-1/q}a(u)K(u,f)\right\|_{q,(0,t)} \geq \left\|u^{\alpha-1/q}a(u)K(u,f)\right\|_{q,\left(\frac{t}{2},t\right)} \geq \left\|u^{\alpha-1/q}a(u)\right\|_{q,\left(\frac{t}{2},t\right)}K\left(\frac{t}{2},f\right)$$
$$\approx t^{\alpha}a(t)K(t,f).$$
Hence,
$$\left\|t^{-\beta-1/r}b(t)\left\|u^{\alpha-1/q}a(u)K(u,f)\right\|_{q,(0,t)}\right\|_{r,(0,\infty)} \succ \left\|t^{\alpha-\beta-1/r}a(t)b(t)K(t,f)\right\|_{r,(0,\infty)}.$$
Let us prove the converse inequality.

**Case $q \leq r$.** Because $\beta > 0$, by Lemma 6, we get
$$\left\|t^{-\beta-1/r}b(t)\left\|u^{\alpha-1/q}a(u)K(u,f)\right\|_{q,(0,t)}\right\|_{r,(0,\infty)} \prec \left\|t^{\alpha-\beta-1/r}a(t)b(t)K(t,f)\right\|_{r,(0,\infty)}.$$
In particular,
$$\left\|t^{-\beta-1/r}b(t)\left\|u^{\alpha-1/r}a(u)K(u,f)\right\|_{r,(0,t)}\right\|_{r,(0,\infty)} \prec \left\|t^{\alpha-\beta-1/r}a(t)b(t)K(t,f)\right\|_{r,(0,\infty)}.$$
**Case $q > r$.** Because $\alpha > -1$, using Lemma 16 and the last estimate, we get
$$\left\|t^{-\beta-1/r}b(t)\left\|u^{\alpha-1/q}a(u)K(u,f)\right\|_{q,(0,t)}\right\|_{r,(0,\infty)}$$
$$\prec \left\|t^{-\beta-1/r}b(t)\left\|u^{\alpha-1/r}a(u)K(u,f)\right\|_{r,(0,t)}\right\|_{r,(0,\infty)} \prec \left\|t^{\alpha-\beta-1/r}a(t)b(t)K(t,f)\right\|_{r,(0,\infty)}.$$
This completes the proof. □

**Lemma 22.** Let $0 \leq \theta_0 < \theta_1 \leq 1$, $0 < p, q, r \leq \infty$, and $a, b, c \in SV$, then, for all $f \in A_0 + A_1$,
$$\left\|t^{\theta_0-\theta_1-1/p}c(t)\left\|s^{-1/r}b(s)\left\|u^{-\theta_0-1/q}a(u)K(u,f)\right\|_{q,(0,s)}\right\|_{r,(0,t)}\right\|_{p,(0,\infty)}$$
$$\approx \left\|t^{-\theta_1-1/p}a(t)b(t)c(t)K(t,f)\right\|_{p,(0,\infty)} = \|f\|_{\theta_1,p,abc} \qquad (14)$$
and
$$\left\|t^{\theta_1-\theta_0-1/p}c(t)\left\|s^{-1/r}b(s)\left\|u^{-\theta_1-1/q}a(u)K(u,f)\right\|_{q,(s,\infty)}\right\|_{r,(t,\infty)}\right\|_{p,(0,\infty)}$$
$$\approx \left\|t^{-\theta_0-1/p}a(t)b(t)c(t)K(t,f)\right\|_{p,(0,\infty)} = \|f\|_{\theta_0,p,abc}.$$

*Proof.* We prove the first equivalence. The second one can be proved analogously or using the first equivalence and symmetry argument. For fixed $r$, denote $I(r) = LHS(14)$.

**Case $r \leq p$.** Because $\theta_0 - \theta_1 < 0$ and $\theta_0 < 1$, by **Lemma 6** (taking $g(s) = b(s)\left\|u^{-\theta_0-1/q}a(u)K(u,f)\right\|_{q,(0,s)}$) and by **Lemma 21**, we get
$$I(r) = \left\|t^{\theta_0-\theta_1-1/p}c(t)\left\|s^{-1/r}g(s)\right\|_{r,(0,t)}\right\|_{p,(0,\infty)} \prec \left\|t^{\theta_0-\theta_1-1/p}c(t)g(t)\right\|_{p,(0,\infty)}$$
$$= \left\|t^{\theta_0-\theta_1-1/p}c(t)b(t)\left\|u^{-\theta_0-1/q}a(u)K(u,f)\right\|_{q,(0,t)}\right\|_{p,(0,\infty)}$$
$$\approx \left\|t^{-\theta_1-1/p}a(t)b(t)c(t)K(t,f)\right\|_{p,(0,\infty)} = RHS(14).$$
In particular, it is shown that
$$I(p) \prec RHS(14). \qquad (15)$$
**Case $r > p$.** By (13) and (15), we get
$$I(r) \prec \left\|t^{\theta_0-\theta_1-1/p}c(t)\left\|s^{-1/p}b(s)\left\|u^{-\theta_0-1/q}a(u)K(u,f)\right\|_{q,(0,s)}\right\|_{p,(0,t)}\right\|_{p,(0,\infty)}$$
$$= I(p) \prec RHS(14).$$
Let us prove the inverse estimate. Using (12), we conclude that
$$RHS(14) = \left\|t^{\theta_0-\theta_1-1/p}c(t)t^{-\theta_0}a(t)b(t)K(t,f)\right\|_{p,(0,\infty)}$$

$$\prec \left\| t^{\theta_0-\theta_1-1/p}c(t)\left\|s^{-1/r}b(s)\|u^{-\theta_0-1/q}a(u)K(u,f)\|_{q,(0,s)}\right\|_{r,(0,t)}\right\|_{p,(0,\infty)}.$$

This completes the proof. □

**Lemma 23.** Let $0 \le \theta_0 < \theta_1 \le 1$, $0 < q, r \le \infty$, and $a, b \in SV$.

(i) If $\left\|s^{-1/r}b(s)\right\|_{r,(0,1)} < \infty$, then, for all $f \in A_0 + A_1$,

$$\left\| t^{\theta_0-\theta_1-1/p}c(t)\left\|s^{-1/r}b(s)\|u^{-\theta_0-1/q}a(u)K(u,f)\|_{q,(s,t)}\right\|_{r,(0,t)}\right\|_{p,(0,\infty)}$$
$$\approx \left\| t^{-\theta_1-1/p}a(t)\|s^{-1/r}b(s)\|_{r,(0,t)}c(t)K(t,f)\right\|_{p,(0,\infty)} = \|f\|_{\theta_1,p,a^\#}, \quad (16)$$

where $a^\#(t) = a(t)\|s^{-1/r}b(s)\|_{r,(0,t)}c(t)$.

(ii) If $\left\|s^{-1/r}b(s)\right\|_{r,(1,\infty)} < \infty$, then, for all $f \in A_0 + A_1$,

$$\left\| t^{\theta_1-\theta_0-1/p}c(t)\left\|s^{-1/r}b(s)\|u^{-\theta_1-1/q}a(u)K(u,f)\|_{q,(t,s)}\right\|_{r,(t,\infty)}\right\|_{p,(0,\infty)}$$
$$\approx \left\| t^{-\theta_0-1/p}a(t)\|s^{-1/r}b(s)\|_{r,(t,\infty)}c(t)K(t,f)\right\|_{p,(0,\infty)} = \|f\|_{\theta_0,p,a^\#},$$

where $a^\#(t) = a(t)\|s^{-1/r}b(s)\|_{r,(t,\infty)}c(t)$.

*Proof.* We prove the first equivalence. The second one can be proved analogously or using symmetry argument. By (9), we have

$$t^{-\theta_0}a(t)\|s^{-1/r}b(s)\|_{r,(0,t)}K(t,f) \prec \left\|s^{-1/r}b(s)\|u^{-\theta_0-1/q}a(u)K(u,f)\|_{q,(s,t)}\right\|_{r,(0,t)}.$$

Hence,

$$\text{RHS}(16) \prec \left\| t^{\theta_0-\theta_1-1/p}c(s)\left\|s^{-1/r}b(s)\|u^{-\theta_0-1/q}a(u)K(u,f)\|_{q,(s,t)}\right\|_{r,(0,t)}\right\|_{p,(0,\infty)}$$
$$= \text{LHS}(16).$$

Let us prove the inverse estimate. By **Lemma 21**, we conclude

$$\text{LHS}(16) \le \left\| t^{\theta_0-\theta_1-1/p}c(t)\|s^{-1/r}b(s)\|_{r,(0,t)}\|u^{-\theta_0-1/q}a(u)K(u,f)\|_{q,(0,t)}\right\|_{p,(0,\infty)}$$
$$\approx \left\| t^{-\theta_1-1/p}c(t)\|s^{-1/r}b(s)\|_{r,(0,t)}a(t)K(t,f)\right\|_{p,(0,\infty)} = \text{RHS}(16).$$

This completes the proof. □

## 4. The limiting case $\theta_1 = 1$

In this section we establish some limiting reiteration formulae for the couples of the form $(*, \bar{A}^{\mathcal{L}}_{1,*})$ and $(*, \bar{A}^{\mathcal{R}}_{1,*})$. Here we need the relevant Holmstedt-type formulae. The proofs of all Holmstedt-type formulae below are based on the paper [2]. Therefore, we adopt the notation from it. Recall that $K(t, f)$ stands for $K(t, f; A_0, A_1)$ and that we will often be using the basic properties of slowly varying functions formulated in Lemma 2.

### 4.1. Couples of the form $(*, \bar{A}^{\mathcal{L}}_{1,*})$

In this subsection we consider the space $\bar{A}^{\mathcal{L}}_{1,r_1,b_1,q_1,a_1}$ as the second operand. Hence, in view of **Lemma 11**, in all theorems of this subsection, we assume that $0 < q_1, r_1 \le \infty$, $a_1, b_1 \in SV$,

$$\left\|s^{-1/r_1}b_1(s)\right\|_{r_1,(1,\infty)} < \infty, \text{ and } \left\|s^{-1/r_1}b_1(s)\|u^{-1/q_1}a_1(u)\|_{q_1,(0,s)}\right\|_{r_1,(0,1)} < \infty.$$

Additionally, we put

$$\chi_1(t) = \left\|s^{-1/r_1}b_1(s)\left\|u^{-1/q_1}a_1(u)\right\|_{q_1,(0,s)}\right\|_{r_1,(0,t)} +$$
$$\left\|s^{-1/r_1}b_1(s)\right\|_{r_1,(t,\infty)}\left\|u^{-1/q_1}a_1(u)\right\|_{q_1,(0,t)}.$$

Note that $\chi_1 \in SV$ and
$$\left\|s^{-1/r_1}b_1(s)\right\|_{r_1,(t,\infty)} a_1(t) \prec \left\|s^{-1/r_1}b_1(s)\right\|_{r_1,(t,\infty)}\left\|u^{-1/q_1}a_1(u)\right\|_{q_1,(0,t)} \prec \chi_1(t). \qquad (17)$$

In proving Holmstedt-type formulae in this subsection, by $\Phi_1$ we denote the function space corresponding to $\bar{A}^{\mathcal{L}}_{1,r_1,b_1,q_1,a_1}$:
$$\|F\|_{\Phi_1} = \left\|s^{-1/r_1}b_1(s)\left\|u^{-1-1/q_1}a_1(u)F(u)\right\|_{q_1,(0,s)}\right\|_{r_1,(0,\infty)}.$$

Following [2], we consider the next auxiliary functions $J(t, f)$, $g_1(t)$, and $h_1(t)$.

$$J(t,f) := \|\chi_{(t,\infty)}(*)K(*,f)\|_{\Phi_1} = \left\|s^{-1/r_1}b_1(s)\left\|u^{-1-1/q_1}a_1(u)K(u,f)\right\|_{q_1,(t,s)}\right\|_{r_1,(t,\infty)}.$$

$$g_1(t) := t\|\chi_{(t,\infty)}(*)\|_{\Phi_1} = t\left\|s^{-1/r_1}b_1(s)\left\|u^{-1-1/q_1}a_1(u)\right\|_{q_1,(t,s)}\right\|_{r_1,(t,\infty)}$$
$$\leq t\left\|s^{-1/r_1}b_1(s)\left\|u^{-1-1/q_1}a_1(u)\right\|_{q_1,(t,\infty)}\right\|_{r_1,(t,\infty)} \approx a_1(t)\left\|s^{-1/r_1}b_1(s)\right\|_{r_1,(t,\infty)}$$

and

$$h_1(t) := \|* \chi_{(0,t)}(*)\|_{\Phi_1} = \left\|s^{-1/r_1}b_1(s)\left\|u^{-1/q_1}a_1(u)\right\|_{q_1,(0,\min(t,s))}\right\|_{r_1,(0,\infty)}$$
$$\approx \left\|s^{-1/r_1}b_1(s)\left\|u^{-1/q_1}a_1(u)\right\|_{q_1,(0,s)}\right\|_{r_1,(0,t)} + \left\|s^{-1/r_1}b_1(s)\right\|_{r_1,(t,\infty)}\left\|u^{-1/q_1}a_1(u)\right\|_{q_1,(0,t)}$$
$$= \chi_1(t).$$

Because of (17), we conclude that
$$g_1(t) \prec h_1(t) \approx g_1(t) + h_1(t) \approx \chi_1(t). \qquad (18)$$

**Theorem 24.** Put $\rho(t) = \frac{t}{\chi_1(t)}$. For all $f \in A_0 + \bar{A}^{\mathcal{L}}_{1,r_1,b_1,q_1,a_1}$ and $t > 0$,
$$K\left(\rho(t), f; A_0, \bar{A}^{\mathcal{L}}_{1,r_1,b_1,q_1,a_1}\right)$$
$$\approx K(t,f) + \rho(t)\left\|s^{-1/r_1}b_1(s)\left\|u^{-1-1/q_1}a_1(u)K(u,f)\right\|_{q_1,(t,s)}\right\|_{r_1,(t,\infty)}.$$

*Proof.* By (18), we arrive at $\rho(t) \approx \frac{t}{g_1(t)+h_1(t)}$. Theorem 3 of [2] completes the proof. □

**Theorem 25.** Let $0 < r \leq \infty$, and $a \in SV$. Put $\rho(t) = \frac{t}{\chi_1(t)}$. For given $0 \leq \theta \leq 1$, denote $a^{\#} = \chi_1^{\theta} a \circ \rho$.

(i) If $0 < \theta < 1$. Then
$$\left(A_0, \bar{A}^{\mathcal{L}}_{1,r_1,b_1,q_1,a_1}\right)_{\theta,r,a} = \bar{A}_{\theta,r,a^{\#}}.$$

(ii) If $\left\|s^{-1/r}a(s)\right\|_{r,(1,\infty)} < \infty$, then
$$\left(A_0, \bar{A}^{\mathcal{L}}_{1,r_1,b_1,q_1,a_1}\right)_{0,r,a} = \bar{A}_{0,r,a\circ\rho}.$$

(iii) If $\left\|s^{-1/r}a(s)\right\|_{r,(0,1)} < \infty$, then
$$\left(A_0, \bar{A}^{\mathcal{L}}_{1,r_1,b_1,q_1,a_1}\right)_{1,r,a} = \bar{A}_{1,r,a^{\#}} \cap \bar{A}^{\mathcal{L,R}}_{1,r,a\circ\rho,r_1,b_1,q_1,a_1}.$$

*Proof.* Denote $\bar{K}(t,f) = K\left(t, f; A_0, \bar{A}^{\mathcal{L}}_{1,r_1,b_1,q_1,a_1}\right)$ and $Z = \left(A_0, \bar{A}^{\mathcal{L}}_{1,r_1,b_1,q_1,a_1}\right)_{\theta,r,a}$ $(0 \leq \theta \leq 1)$. By Remark 4 (v) and Theorem 24 we can write
$$\|f\|_Z = \left\|u^{-\theta-1/r}a(u)\bar{K}(u,f)\right\|_{r,(0,\infty)}$$
$$\approx \left\|\rho(t)^{-\theta}t^{-1/r}a(\rho(t))\bar{K}(\rho(t),f)\right\|_{r,(0,\infty)} \approx I_1 + I_2,$$

where
$$I_1 := \|\rho(t)^{-\theta} t^{-1/r} a(\rho(t)) K(t,f)\|_{r,(0,\infty)}$$
$$= \|t^{-\theta-1/r} \chi_1(t)^\theta a(\rho(t)) K(t,f)\|_{r,(0,\infty)} = \|f\|_{\theta,r;a^\#}$$

and
$$I_2 := \left\|\rho(t)^{1-\theta} t^{-1/r} a(\rho(t)) \left\|s^{-1/r_1} b_1(s) \|u^{-1-1/q_1} a_1(u) K(u,f)\|_{q_1,(t,s)}\right\|_{r_1,(t,\infty)}\right\|_{r,(0,\infty)}$$
$$= \left\|t^{1-\theta-1/r} \chi_1(t)^{\theta-1} a(\rho(t)) \left\|s^{-1/r_1} b_1(s) \|u^{-1-1/q_1} a_1(u) K(u,f)\|_{q_1,(t,s)}\right\|_{r_1,(t,\infty)}\right\|_{r,(0,\infty)}.$$

If $\theta = 1$, we have
$$I_2 = \left\|t^{-1/r} a(\rho(t)) \left\|s^{-1/r_1} b_1(s) \|u^{-1-1/q_1} a_1(u) K(u,f)\|_{q_1,(t,s)}\right\|_{r_1,(t,\infty)}\right\|_{r,(0,\infty)}$$
$$= \|f\|_{\mathcal{L},\mathcal{R};1,r,a\circ\rho,r_1,b_1,q_1,a_1}.$$

Thus, the statement (iii) is proved. Consider the case $0 \le \theta < 1$. To proof the statements (i) and (ii) it is enough to show that $I_2 \prec I_1$. By Lemma 23 and (17), we get
$$I_2 \approx \left\|t^{-\theta-1/r} \chi_1(t)^{\theta-1} a(\rho(t)) \|s^{-1/r_1} b_1(s)\|_{r_1,(t,\infty)} a_1(t) K(t,f)\right\|_{r,(0,\infty)}$$
$$\prec \|t^{-\theta-1/r} \chi_1(t)^\theta a(\rho(t)) K(t,f)\|_{r,(0,\infty)} \approx I_1.$$

This completes the proof. $\square$

**Theorem 26.** Let $0 < r_0 \le \infty$, $b_0 \in SV$, and $\|s^{-1/r_0} b_0(s)\|_{r_0,(1,\infty)} < \infty$. Put $\chi_0(t) = \|s^{-1/r_0} b_0(s)\|_{r_0,(t,\infty)}$ and $\rho(t) = t\frac{\chi_0(t)}{\chi_1(t)}$. Then, for all $f \in A_{0,r_0,b_0} + \bar{A}^{\mathcal{L}}_{1,r_1,b_1,q_1,a_1}$ and $t > 0$,
$$K(\rho(t), f; A_{0,r_0,b_0}, \bar{A}^{\mathcal{L}}_{1,r_1,b_1,q_1,a_1}) \approx \|u^{-1/r_0} b_0(u) K(u,f)\|_{r_0,(0,t)}$$
$$+ \rho(t) \left\|s^{-1/r_1} b_1(s) \|u^{-1-1/q_1} a_1(u) K(u,f)\|_{q_1,(t,s)}\right\|_{r_1,(t,\infty)} + \chi_0(t) K(t,f). \quad (19)$$

*Proof.* Let $\Phi_0$ be the function space corresponding to $\bar{A}_{0,r_0,b_0}$:
$$\|F\|_{\Phi_0} = \|u^{-1/r_0} b_0(u) F(u)\|_{r_0,(0,\infty)}.$$

Following [2], we consider the functions
$$I(t,f) := \|\chi_{(0,t)}(*) K(*,f)\|_{\Phi_0} = \|u^{-1/r_0} b_0(u) K(u,f)\|_{r_0,(0,t)},$$
$$g_0(t) := t\|\chi_{(t,\infty)}(*)\|_{\Phi_0} = t\|u^{-1/r_0} b_0(u)\|_{r_0,(t,\infty)} = t\chi_0(t),$$
and
$$h_0(t) := \|* \chi_{(0,t)}(*)\|_{\Phi_0} = \|u^{1-1/r_0} b_0(u)\|_{r_0,(0,t)} \approx tb_0(t).$$

Thus,
$$h_0(t) \prec g_0(t) \approx g_0(t) + h_0(t) \approx t\chi_0(t).$$

Now we estimate the auxiliary functions $g(t)$ and $h(t)$.
$$\frac{1}{g(t)} := \left\|\chi_{(t,\infty)}(u) \frac{u}{g_0(u) + h_0(u)}\right\|_{\Phi_1}$$
$$= \left\|s^{-1/r_1} b_1(s) \left\|u^{-1-1/q_1} a_1(u) \chi_{(t,\infty)}(u) \frac{u}{g_0(u) + h_0(u)}\right\|_{q_1,(0,s)}\right\|_{r_1,(0,\infty)}$$
$$\approx \left\|s^{-1/r_1} b_1(s) \left\|u^{-1-1/q_1} \frac{a_1(u)}{\chi_0(u)}\right\|_{q_1,(t,s)}\right\|_{r_1,(t,\infty)}$$
$$\le \left\|u^{-1-1/q_1} \frac{a_1(u)}{\chi_0(u)}\right\|_{q_1,(t,\infty)} \|s^{-1/r_1} b_1(s)\|_{r_1,(t,\infty)} \approx t^{-1} \frac{a_1(t)}{\chi_0(t)} \|s^{-1/r_1} b_1(s)\|_{r_1,(t,\infty)}.$$

By (17), we have
$$\frac{1}{g(t)} \prec t^{-1}\frac{\chi_1(t)}{\chi_0(t)} = \frac{1}{\rho(t)}.$$
$$h(t) := \left\|\chi_{(0,t)}(u)\frac{u}{g_1(u)+h_1(u)}\right\|_{\Phi_0} \approx \left\|u^{1-1/r_0}\frac{b_0(u)}{\chi_1(u)}\right\|_{r_0,(0,t)} \approx t\frac{b_0(t)}{\chi_1(t)} \prec t\frac{\chi_0(t)}{\chi_1(t)} = \rho(t).$$
Here we have used the estimate (18). Hence, we arrive at
$$h(t) \approx t\frac{b_0(t)}{\chi_1(t)} \prec \rho(t) \prec g(t).$$
Using [2, Theorem 4, Case 1], we get
$$K\bigl(\rho(t), f; A_{0,r_0,b_0}, \bar{A}^{\mathcal{L}}_{1,r_1,b_1,q_1,a_1}\bigr) \approx \left\|u^{-1/r_0}b_0(u)K(u,f)\right\|_{r_0,(0,t)}$$
$$+ \rho(t)\left\|s^{-1/r_1}b_1(s)\left\|u^{-1-1/q_1}a_1(u)K(u,f)\right\|_{q_1,(t,s)}\right\|_{r_1,(t,\infty)}$$
$$+ \{t\chi_0(t) + \rho(t)\chi_1(t)\}t^{-1}K(t,f).$$
Note that
$$\{t\chi_0(t) + \rho(t)\chi_1(t)\}t^{-1} = 2\chi_0(t).$$
Thus, we obtain formula (19). □

**Theorem 27.** Let $0 < r, r_0 \leq \infty$, $a, b_0 \in SV$, and $\left\|u^{-1/r_0}b_0(u)\right\|_{r_0,(1,\infty)} < \infty$. Put $\chi_0(t) := \left\|u^{-1/r_0}b_0(u)\right\|_{r_0,(t,\infty)}$ and $\rho(t) = t\frac{\chi_0(t)}{\chi_1(t)}$. For given $0 \leq \theta \leq 1$, denote $a^\# = \chi_0^{(1-\theta)}\chi_1^{\theta}a \circ \rho$.

(i) If $0 < \theta < 1$, then
$$\bigl(\bar{A}_{0,r_0,b_0}, \bar{A}^{\mathcal{L}}_{1,r_1,b_1,q_1,a_1}\bigr)_{\theta,r,a} = \bar{A}_{\theta,r,a^\#}.$$

(ii) If $\left\|s^{-1/r}a(s)\right\|_{r,(1,\infty)} < \infty$, then
$$\bigl(\bar{A}_{0,r_0,b_0}, \bar{A}^{\mathcal{L}}_{1,r_1,b_1,q_1,a_1}\bigr)_{0,r,a} = \bar{A}_{0,r,a^\#} \cap \bar{A}^{\mathcal{L}}_{0,r,a\circ\rho,r_0,b_0}.$$

(iii) If $\left\|s^{-1/r}a(s)\right\|_{r,(0,1)} < \infty$, then
$$\bigl(\bar{A}_{0,r_0,b_0}, \bar{A}^{\mathcal{L}}_{1,r_1,b_1,q_1,a_1}\bigr)_{1,r,a} = \bar{A}_{1,r,a^\#} \cap \bar{A}^{\mathcal{L},\mathcal{R}}_{1,r,a\circ\rho,r_1,b_1,q_1,a_1}.$$

*Proof.* Denote $\bar{K}(t,f) = K\bigl(t,f; \bar{A}_{0,r_0,b_0}, \bar{A}^{\mathcal{L}}_{1,r_1,b_1,q_1,a_1}\bigr)$, and $Z = \bigl(\bar{A}_{0,r_0,b_0}, \bar{A}^{\mathcal{L}}_{1,r_1,b_1,q_1,a_1}\bigr)_{\theta,r,a}$ ($0 \leq \theta \leq 1$). By Remark 4 (v) and (19) we can write
$$\|f\|_Z \approx \left\|\rho(t)^{-\theta}t^{-1/r}a(\rho(t))\bar{K}(\rho(t),f)\right\|_{r,(0,\infty)}$$
$$= \left\|t^{-\theta-1/r}a(\rho(t))\chi_0(t)^{-\theta}\chi_1(t)^{\theta}\bar{K}(\rho(t),f)\right\|_{r,(0,\infty)} \approx I_1(\theta) + I_2(\theta) + I_3(\theta),$$
where
$$I_1(\theta) := \left\|t^{-\theta-1/r}a(\rho(t))\chi_0(t)^{-\theta}\chi_1(t)^{\theta}\left\|u^{-1/r_0}b_0(u)K(u,f)\right\|_{r_0,(0,t)}\right\|_{r,(0,\infty)},$$
$$I_2(\theta) := \left\|\rho(t)^{1-\theta}t^{-1/r}a(\rho(t))\left\|s^{-1/r_1}b_1(s)\left\|u^{-1-1/q_1}a_1(u)K(u,f)\right\|_{q_1,(t,s)}\right\|_{r_1,(t,\infty)}\right\|_{r,(0,\infty)}$$
$$=$$
$$\left\|t^{1-\theta-\frac{1}{r}}a(\rho(t))\chi_0(t)^{1-\theta}\chi_1(t)^{\theta-1}\left\|s^{-\frac{1}{r_1}}b_1(s)\left\|u^{-1-\frac{1}{q_1}}a_1(u)K(u,f)\right\|_{q_1,(t,s)}\right\|_{r_1,(t,\infty)}\right\|_{r,(0,\infty)},$$
and
$$I_3(\theta) := \left\|t^{-\theta-1/r}a(\rho(t))\chi_0(t)^{-\theta}\chi_1(t)^{\theta}\chi_0(t)K(t,f)\right\|_{r,(0,\infty)}$$
$$= \left\|t^{-\theta-1/r}a(\rho(t))\chi_0(t)^{1-\theta}\chi_1(t)^{\theta}K(t,f)\right\|_{r,(0,\infty)} = \|f\|_{\theta,r;a^\#}.$$
Consider $I_1$. If $0 < \theta \leq 1$, using the first equivalence from Lemma 21, we conclude that
$$I_1(\theta) \approx \left\|t^{-\theta-1/r}a(\rho(t))\chi_0(t)^{-\theta}\chi_1(t)^{\theta}b_0(t)K(t,f)\right\|_{r,(0,\infty)}$$
$$\prec \left\|t^{-\theta-1/r}a(\rho(t))\chi_0(t)^{1-\theta}\chi_1(t)^{\theta}K(t,f)\right\|_{r,(0,\infty)} = I_3(\theta).$$

Hence, to prove the statement (i), it is enough to show that $I_2(\theta) \prec I_3(\theta)$ if $0 < \theta < 1$. If $\theta = 0$, we obtain
$$I_1(0) = \left\|t^{-1/r}a(\rho(t))\left\|u^{-1/r_0}b_0(u)K(u,f)\right\|_{r_0,(0,t)}\right\|_{r,(0,\infty)} = \|f\|_{\mathcal{L};0,r,a\circ\rho,r_0,b_0}.$$
Hence, to prove the statement (ii), it is enough to show that $I_2(0) \prec I_3(0)$ or $I_2(0) \prec I_1(0)$.
Consider $I_2$. We have
$$I_2(1) = \left\|t^{-1/r}a(\rho(t))\left\|s^{-1/r_1}b_1(s)\left\|u^{-1-1/q_1}a_1(u)K(u,f)\right\|_{q_1,(t,s)}\right\|_{r_1,(t,\infty)}\right\|_{r,(0,\infty)}$$
$$= \|f\|_{\mathcal{L},\mathcal{R};1,r,a\circ\rho,r_1,b_1,q_1,a_1}.$$
Hence, the statement (iii) is proved. Let $0 \leq \theta < 1$. The proof of the statements (i) and (ii) is completed when we show that $I_2(\theta) \prec I_3(\theta)$. From Lemma 23 (ii) and (17), it follows
$$I_2(\theta) \approx \left\|t^{-\theta-1/r}a(\rho(t))\chi_0(t)^{1-\theta}\chi_1(t)^{\theta-1}\left\|s^{-1/r_1}b_1(s)\right\|_{r_1,(t,\infty)}a_1(t)K(t,f)\right\|_{r,(0,\infty)}$$
$$\prec \left\|t^{-\theta-1/r}\chi_0(t)^{1-\theta}\chi_1(t)^{\theta}a(\rho(t))K(t,f)\right\|_{r,(0,\infty)} \approx I_3(\theta).$$
This completes the proof. □

**Theorem 28.** Let $0 < \theta_0 < 1$, $0 < r_0 \leq \infty$, and $b_0 \in SV$. Put $\rho(t) = t^{1-\theta_0}\frac{b_0(t)}{\chi_1(t)}$. Then, for all $f \in \bar{A}_{\theta_0,r_0,b_0} + \bar{A}^{\mathcal{L}}_{1,r_1,b_1,q_1,a_1}$ and $t > 0$,
$$K(\rho(t),f;\bar{A}_{\theta_0,r_0,b_0},\bar{A}^{\mathcal{L}}_{1,r_1,b_1,q_1,a_1}) \approx \left\|u^{-\theta_0-1/r_0}b_0(u)K(u,f)\right\|_{r_0,(0,t)} +$$
$$\rho(t)\left\|s^{-1/r_1}b_1(s)\left\|u^{-1-1/q_1}a_1(u)K(u,f)\right\|_{q_1,(t,s)}\right\|_{r_1,(t,\infty)}.$$

*Proof.* Let $\Phi_0$ be the function space corresponding to $\bar{A}_{\theta_0,r_0,b_0}$:
$$\|F\|_{\Phi_0} = \left\|u^{-\theta_0-1/r_0}b_0(u)F(u)\right\|_{r_0,(0,\infty)}.$$
From the proof of [8, Theorem 18], it is known that
$$I(t,f) = \left\|u^{-\theta_0-1/r_0}b_0(u)K(u,f)\right\|_{r_0,(0,t)},$$
and
$$h_0(t) \approx g_0(t) \approx g_0(t) + h_0(t) \approx t^{1-\theta_0}b_0(t).$$
Now, using (17) and (18), we estimate the functions $g(t)$ and $h(t)$.
$$\frac{1}{g(t)} := \left\|\chi_{(t,\infty)}(u)\frac{u}{g_0(u)+h_0(u)}\right\|_{\Phi_1}$$
$$\approx \left\|s^{-1/r_1}b_1(s)\left\|u^{\theta_0-1-1/q_1}\frac{a_1(u)}{b_0(u)}\right\|_{q_1,(t,s)}\right\|_{r_1,(t,\infty)}$$
$$\leq \left\|u^{\theta_0-1-1/q_1}\frac{a_1(u)}{b_0(u)}\right\|_{q_1,(t,\infty)}\left\|s^{-1/r_1}b_1(s)\right\|_{r_1,(t,\infty)}$$
$$\approx t^{\theta_0-1}\frac{a_1(t)}{b_0(t)}\left\|s^{-1/r_1}b_1(s)\right\|_{r_1,(t,\infty)} \prec t^{\theta_0-1}\frac{\chi_1(t)}{b_0(t)} = \frac{1}{\rho(t)}.$$
$$h(t) := \left\|\chi_{(0,t)}(u)\frac{u}{g_1(u)+h_1(u)}\right\|_{\Phi_0} = \left\|u^{1-\theta_0-1/r_0}\frac{b_0(u)}{g_1(u)+h_1(u)}\right\|_{r_0,(0,t)}$$
$$\approx \left\|u^{1-\theta_0-1/r_0}\frac{b_0(u)}{\chi_1(u)}\right\|_{r_0,(0,t)} \approx t^{1-\theta_0}\frac{b_0(t)}{\chi_1(t)} = \rho(t).$$
So,
$$h(t) \approx \rho(t) \prec g(t).$$
Finally, note that
$$\frac{h_0(t)}{h_1(t)} \approx \frac{t^{1-\theta_0}b_0(t)}{\chi_1(t)} = \rho(t).$$

Theorem 4 (Case 4) of [2] completes the proof. □

**Theorem** 29. Let $0 < \theta_0 < 1$, $0 < r, r_0 \leq \infty$, and $a, b_0 \in SV$. Put $\rho(t) = t^{1-\theta_0}\frac{b_0(t)}{\chi_1(t)}$. For given $0 \leq \theta \leq 1$, denote $a^{\#} = b_0^{1-\theta}\chi_1^{\theta}a\circ\rho$ and $\eta = (1-\theta)\theta_0 + \theta$.

(i) If $0 < \theta < 1$, then
$$\left(\bar{A}_{\theta_0,r_0,b_0}, \bar{A}^{\mathcal{L}}_{1,r_1,b_1,q_1,a_1}\right)_{\theta,r,a} = \bar{A}_{\eta,r,a^{\#}},$$
where $\eta = (1-\theta)\theta_0 + \theta$.

(ii) If $\left\|s^{-1/r}a(s)\right\|_{r,(1,\infty)} < \infty$, then
$$\left(\bar{A}_{\theta_0,r_0,b_0}, \bar{A}^{\mathcal{L}}_{1,r_1,b_1,q_1,a_1}\right)_{0,r,a} = \bar{A}^{\mathcal{L}}_{\theta_0,r,a\circ\rho,r_0,b_0}.$$

(iii) If $\left\|s^{-1/r}a(s)\right\|_{r,(0,1)} < \infty$, then
$$\left(\bar{A}_{\theta_0,r_0,b_0}, \bar{A}^{\mathcal{L}}_{1,r_1,b_1,q_1,a_1}\right)_{1,r,a} = \bar{A}_{1,r,a^{\#}} \cap \bar{A}^{\mathcal{L},\mathcal{R}}_{1,r,a\circ\rho,r_1,b_1,q_1,a_1}.$$

*Proof.* Denote $\bar{K}(t,f) = K(t,f;\bar{A}_{\theta_0,r_0,b_0}, \bar{A}^{\mathcal{L}}_{1,r_1,b_1,q_1,a_1})$ and $Z = \left(\bar{A}_{\theta_0,r_0,b_0}, \bar{A}^{\mathcal{L}}_{1,r_1,b_1,q_1,a_1}\right)_{\theta,r,a}$ ($0 \leq \theta \leq 1$). We write for shortness
$$e(t) = b_0(t)^{-\theta}\chi_1(t)^{\theta}a(\rho(t)) = a^{\#}(t)b_0(t)^{-1}.$$

By Remark 4 (v) and Theorem 28 we can write
$$\|f\|_Z \approx \left\|\rho(t)^{-\theta}t^{-1/r}a(\rho(t))\bar{K}(\rho(t),f)\right\|_{r,(0,\infty)}$$
$$= \left\|t^{-\theta(1-\theta_0)-1/r}e(t)\bar{K}(\rho(t),f)\right\|_{r,(0,\infty)} \approx I_1(\theta) + I_2(\theta),$$

where
$$I_1(\theta) := \left\|t^{-\theta(1-\theta_0)-1/r}e(t)\left\|u^{-\theta_0-1/r_0}b_0(u)K(u,f)\right\|_{r_0,(0,t)}\right\|_{r,(0,\infty)}$$

and
$$I_2(\theta) := \left\|\rho(t)^{-\theta}t^{-\frac{1}{r}}a(\rho(t))\rho(t)\left\|s^{-\frac{1}{r_1}}b_1(s)\left\|u^{-1-\frac{1}{q_1}}a_1(u)K(u,f)\right\|_{q_1,(t,s)}\right\|_{r_1,(t,\infty)}\right\|_{r,(0,\infty)}$$
$$=$$
$$\left\|t^{(1-\theta)(1-\theta_0)-\frac{1}{r}}\left(\frac{b_0(t)}{\chi_1(t)}\right)^{1-\theta}a(\rho(t))\left\|s^{-\frac{1}{r_1}}b_1(s)\left\|u^{-1-\frac{1}{q_1}}a_1(u)K(u,f)\right\|_{q_1,(t,s)}\right\|_{r_1,(t,\infty)}\right\|_{r,(0,\infty)}.$$

Consider the case $\theta = 0$. We have
$$I_1(0) = \left\|t^{-1/r}a(\rho(t))\left\|u^{-\theta_0-1/r_0}b_0(u)K(u,f)\right\|_{r_0,(0,t)}\right\|_{r,(0,\infty)} = \|f\|_{\mathcal{L};\theta_0,r,a\circ\rho,r_0,b_0}.$$

By Lemma 23 (ii) and (17), we observe
$$I_2(0) = \left\|t^{1-\theta_0-\frac{1}{r}}\frac{b_0(t)}{\chi_1(t)}a(\rho(t))\left\|s^{-\frac{1}{r_1}}b_1(s)\left\|u^{-1-\frac{1}{q_1}}a_1(u)K(u,f)\right\|_{q_1,(t,s)}\right\|_{r_1,(t,\infty)}\right\|_{r,(0,\infty)}$$
$$\approx \left\|t^{-\theta_0-1/r}\frac{b_0(t)}{\chi_1(t)}a(\rho(t))\|s^{-1/r_1}b_1(s)\|_{r_1,(t,\infty)}a_1(t)K(t,f)\right\|_{r,(0,\infty)}$$
$$\prec \left\|t^{-\theta_0-1/r}a(\rho(t))b_0(t)K(t,f)\right\|_{r,(0,\infty)} = \|f\|_{\theta_0,r;b_0a\circ\rho}.$$

By Corollary 17 (i), we know that $\bar{A}^{\mathcal{L}}_{\theta_0,r,a\circ\rho,r_0,b_0} \subset \bar{A}_{\theta_0,r,b_0a\circ\rho}$. This means, that
$$I_2(0) \prec \|f\|_{\theta_0,r;b_0a\circ\rho} \prec \|f\|_{\mathcal{L};\theta_0,r,a\circ\rho,r_0,b_0} = I_1(0).$$

Thus, the statement (ii) is proved. Consider $I_1(\theta)$ for $0 < \theta \leq 1$. By Lemma 21, we have
$$I_1(\theta) \approx \left\|t^{-\theta(1-\theta_0)-\theta_0-1/r}e(t)b_0(t)K(t,f)\right\|_{r,(0,\infty)}$$
$$= \left\|t^{-\theta(1-\theta_0)-\theta_0-1/r}a^{\#}(t)K(t,f)\right\|_{r,(0,\infty)}.$$

Hence,

$$I_1(\theta) \approx \|f\|_{\eta,r;a^\#} \qquad \text{provided } 0 < \theta \leq 1, \tag{20}$$

Consider $I_2(\theta)$ for $0 < \theta < 1$. Note that $(1-\theta)(1-\theta_0) = 1-\eta$. By Lemma 23 (ii) and (17), we get

$$I_2(\theta) \approx \left\| t^{-\eta-1/r} b_0(t)^{(1-\theta)} \chi_1(t)^{\theta-1} a(\rho(t)) \left\| s^{-1/r_1} b_1(s) \right\|_{r_1,(t,\infty)} a_1(t) K(t,f) \right\|_{r,(0,\infty)}$$
$$\prec \left\| t^{-\eta-1/r} b_0(t)^{(1-\theta)} \chi_1(t)^{\theta} a(\rho(t)) K(t,f) \right\|_{r,(0,\infty)} = \|f\|_{\eta,r;a^\#}.$$

So, taking into account (20), we arrive at $I_1(\theta) + I_2(\theta) \approx \|f\|_{\eta,r;a^\#}$. Thus, the statement (i) is proved. For $I_2(1)$, We have

$$I_2(1) = \left\| t^{-1/r} a(\rho(t)) \left\| s^{-1/r_1} b_1(s) \left\| u^{-1-1/q_1} a_1(u) K(u,f) \right\|_{q_1,(t,s)} \right\|_{r_1,(t,\infty)} \right\|_{r,(0,\infty)}$$
$$= \|f\|_{\mathcal{L},\mathcal{R};1,r,a^\circ\rho,r_1,b_1,q_1,a_1}.$$

Together with (20), this means that the statement (iii) is proved. $\square$

**Theorem 30.** Let $0 < r_0, q_0 \leq \infty$, $a_0, b_0 \in SV$, and

$$\left\| s^{-1/r_0} b_0(s) \left\| u^{-1/q_0} a_0(u) \right\|_{q_0,(1,s)} \right\|_{r_0,(1,\infty)} < \infty.$$

Put

$$\chi_0(t) = \left\| s^{-1/r_0} b_0(s) \left\| u^{-1/q_0} a_0(u) \right\|_{q_0,(t,s)} \right\|_{r_0,(t,\infty)}$$

and $\rho(t) = t \frac{\chi_0(t)}{\chi_1(t)}$. Then, for all $f \in \bar{A}^{\mathcal{L}}_{0,r_0,b_0,q_0,a_0} + \bar{A}^{\mathcal{L}}_{1,r_1,b_1,q_1,a_1}$ and $t > 0$,

$$K\left(\rho(t), f; \bar{A}^{\mathcal{L}}_{0,r_0,b_0,q_0,a_0}, \bar{A}^{\mathcal{L}}_{1,r_1,b_1,q_1,a_1}\right)$$
$$\approx \left\| s^{-1/r_0} b_0(s) \left\| u^{-1/q_0} a_0(u) K(u,f) \right\|_{q_0,(0,s)} \right\|_{r_0,(0,t)}$$
$$+ \left\| s^{-1/r_0} b_0(s) \right\|_{r_0,(t,\infty)} \left\| u^{-1/q_0} a_0(u) K(u,f) \right\|_{q_0,(0,t)}$$
$$+ \rho(t) \left\| s^{-1/r_1} b_1(s) \left\| u^{-1-1/q_1} a_1(u) K(u,f) \right\|_{q_1,(t,s)} \right\|_{r_1,(t,\infty)} + \chi_0(t) K(t,f). \tag{21}$$

*Proof.* Let $\Phi_0$ be the function space corresponding to $\bar{A}^{\mathcal{L}}_{0,r_0,b_0,q_0,a_0}$:

$$\|F\|_{\Phi_0} = \left\| s^{-1/r_0} b_0(s) \left\| u^{-1/q_0} a_0(u) F(u) \right\|_{q_0,(0,s)} \right\|_{r_0,(0,\infty)}.$$

Let us calculate the functions $I(t,f)$, $g_0(t)$, and $h_0(t)$.

$$I(t,f) := \left\| \chi_{(0,t)}(*) K(*,f) \right\|_{\Phi_0}$$
$$= \left\| s^{-1/r_0} b_0(s) \left\| u^{-1/q_0} a_0(u) \chi_{(0,t)}(u) K(u,f) \right\|_{q_0,(0,s)} \right\|_{r_0,(0,\infty)}$$
$$\approx \left\| s^{-1/r_0} b_0(s) \left\| u^{-1/q_0} a_0(u) K(u,f) \right\|_{q_0,(0,s)} \right\|_{r_0,(0,t)} +$$
$$\left\| s^{-1/r_0} b_0(s) \right\|_{r_0,(t,\infty)} \left\| u^{-1/q_0} a_0(u) K(u,f) \right\|_{q_0,(0,t)}.$$

$$g_0(t) := t \left\| \chi_{(t,\infty)}(*) \right\|_{\Phi_0} = t \left\| s^{-1/r_0} b_0(s) \left\| u^{-1/q_0} a_0(u) \right\|_{q_0,(t,s)} \right\|_{r_0,(t,\infty)} = t\chi_0(t).$$

$$h_0(t) := \left\| * \chi_{(0,t)}(*) \right\|_{\Phi_0} = \left\| s^{-1/r_0} b_0(s) \left\| u^{1-1/q_0} a_0(u) \right\|_{q_0,(0,\min(t,s))} \right\|_{r_0,(0,\infty)}$$
$$\approx \left\| s^{-1/r_0} b_0(s) \left\| u^{1-1/q_0} a_0(u) \right\|_{q_0,(0,s)} \right\|_{r_0,(0,t)} + \left\| s^{-1/r_0} b_0(s) \left\| u^{1-1/q_0} a_0(u) \right\|_{q_0,(0,t)} \right\|_{r_0,(t,\infty)}$$
$$\approx \left\| s^{1-1/r_0} b_0(s) a_0(s) \right\|_{r_0,(0,t)} + \left\| s^{-1/r_0} b_0(s) \right\|_{r_0,(t,\infty)} t a_0(t)$$
$$\approx t b_0(t) a_0(t) + \left\| s^{-1/r_0} b_0(s) \right\|_{r_0,(t,\infty)} t a_0(t) \approx t a_0(t) \left\| s^{-1/r_0} b_0(s) \right\|_{r_0,(t,\infty)}.$$

Note that

$$\chi_0(t) > \left\| \left\| s^{-1/r_0} b_0(s) \right\| u^{-1/q_0} a_0(u) \right\|_{q_0,(t,s)} \right\|_{r_0,(2t,\infty)}$$
$$> \left\| s^{-1/r_0} b_0(s) \right\|_{r_0,(2t,\infty)} \left\| u^{-1/q_0} a_0(u) \right\|_{q_0,(t,2t)} \approx a_0(t) \left\| s^{-1/r_0} b_0(s) \right\|_{r_0,(t,\infty)}. \quad (22)$$

Hence,
$$h_0(t) \prec g_0(t) \approx g_0(t) + h_0(t) \approx t\chi_0(t).$$

Now we estimate the functions $g(t)$ and $h(t)$.
$$\frac{1}{g(t)} := \left\| \chi_{(t,\infty)}(u) \frac{u}{g_0(u) + h_0(u)} \right\|_{\Phi_1}$$
$$= \left\| \left\| s^{-1/r_1} b_1(s) \left\| u^{-1-1/q_1} a_1(u) \frac{u}{g_0(u) + h_0(u)} \right\|_{q_1,(t,s)} \right\|_{r_1,(t,\infty)}$$
$$\approx \left\| \left\| s^{-1/r_1} b_1(s) \left\| u^{-1-1/q_1} \frac{a_1(u)}{\chi_0(u)} \right\|_{q_1,(t,s)} \right\|_{r_1,(t,\infty)}$$
$$\prec \left\| \left\| s^{-1/r_1} b_1(s) \left\| u^{-1-1/q_1} \frac{a_1(u)}{\chi_0(u)} \right\|_{q_1,(t,\infty)} \right\|_{r_1,(t,\infty)} \approx t^{-1} \frac{a_1(t)}{\chi_0(t)} \left\| s^{-1/r_1} b_1(s) \right\|_{r_1,(t,\infty)}.$$

Hence by (17), $\rho(t) \prec g(t)$. Using (18), we get
$$h(t) := \left\| \chi_{(0,t)}(u) \frac{u}{g_1(u) + h_1(u)} \right\|_{\Phi_0}$$
$$= \left\| \left\| s^{-1/r_0} b_0(s) \left\| u^{-1/q_0} a_0(u) \chi_{(0,t)}(u) \frac{u}{g_1(u) + h_1(u)} \right\|_{q_0,(0,s)} \right\|_{r_0,(0,\infty)}$$
$$\approx \left\| \left\| s^{-1/r_0} b_0(s) \left\| u^{1-1/q_0} \frac{a_0(u)}{\chi_1(u)} \right\|_{q_0,(0,\min(t,s))} \right\|_{r_0,(0,\infty)}$$
$$\approx \left\| \left\| s^{-1/r_0} b_0(s) \left\| u^{1-1/q_0} \frac{a_0(u)}{\chi_1(u)} \right\|_{q_0,(0,s)} \right\|_{r_0,(0,t)}$$
$$+ \left\| \left\| s^{-1/r_0} b_0(s) \left\| u^{1-1/q_0} \frac{a_0(u)}{\chi_1(u)} \right\|_{q_0,(0,t)} \right\|_{r_0,(t,\infty)}$$
$$\approx \left\| s^{1-1/r_0} b_0(s) \frac{a_0(s)}{\chi_1(s)} \right\|_{r_0,(0,t)} + t \frac{a_0(t)}{\chi_1(t)} \left\| s^{-1/r_0} b_0(s) \right\|_{r_0,(t,\infty)}$$
$$\approx tb_0(t) \frac{a_0(t)}{\chi_1(t)} + t \frac{a_0(t)}{\chi_1(t)} \left\| s^{-1/r_0} b_0(s) \right\|_{r_0,(t,\infty)} \approx t \frac{a_0(t)}{\chi_1(t)} \left\| s^{-1/r_0} b_0(s) \right\|_{r_0,(t,\infty)}.$$

Because of (22), we get $h(t) \prec \rho(t)$. In addition, we have
$$\frac{g_0(t)}{g_1(t) + h_1(t)} \approx \frac{t\chi_0(t)}{\chi_1(t)} = \rho(t) \quad \text{and} \quad \frac{g_0(t) + h_0(t)}{h_1(t)} \approx \frac{t\chi_0(t)}{\chi_1(t)} = \rho(t).$$

Theorem 4 (Case 1) of [2] completes the proof. □

**Theorem** 31. Let $0 < r, r_0, q_0 \leq \infty$, $a, a_0, b_0 \in SV$, and
$$\left\| \left\| s^{-1/r_0} b_0(s) \right\| u^{-1/q_0} a_0(u) \right\|_{q_0,(1,s)} \right\|_{r_0,(1,\infty)} < \infty.$$

Put $\chi_0(t) = \left\| \left\| s^{-1/r_0} b_0(s) \right\| u^{-1/q_0} a_0(u) \right\|_{q_0,(t,s)} \right\|_{r_0,(t,\infty)}$ and $\rho(t) = t \frac{\chi_0(t)}{\chi_1(t)}$. For given $0 \leq \theta \leq 1$, denote $a^{\#} = \chi_0^{(1-\theta)} \chi_1^{\theta} a \circ \rho$.

(i) If $0 < \theta < 1$, then
$$\left( \bar{A}^{\mathcal{L}}_{0,r_0,b_0,q_0,a_0}, \bar{A}^{\mathcal{L}}_{1,r_1,b_1,q_1,a_1} \right)_{\theta,r,a} = \bar{A}_{\theta,r,a^{\#}}.$$

(ii) If $\left\|s^{-1/r}a(s)\right\|_{r,(1,\infty)} < \infty$, then
$$\left(\bar{A}^{\mathcal{L}}_{0,r_0,b_0,q_0,a_0}, \bar{A}^{\mathcal{L}}_{1,r_1,b_1,q_1,a_1}\right)_{0,r,a} = \bar{A}_{0,r,a^\#} \cap \bar{A}^{\mathcal{L}}_{0,r,B_0,q_0,a_0} \cap \bar{A}^{\mathcal{L},\mathcal{L}}_{0,r,a^\circ\rho,r_0,b_0,q_0,a_0},$$
where $B_0(t) = \left\|s^{-1/r_0}b_0(s)\right\|_{r_0,(t,\infty)} a(\rho(t))$.

(iii) If $\left\|s^{-1/r}a(s)\right\|_{r,(0,1)} < \infty$, then
$$\left(\bar{A}^{\mathcal{L}}_{0,r_0,b_0,q_0,a_0}, \bar{A}^{\mathcal{L}}_{1,r_1,b_1,q_1,a_1}\right)_{1,r,a} = \bar{A}_{1,r,a^\#} \cap \bar{A}^{\mathcal{L},\mathcal{R}}_{1,r,a^\circ\rho,r_1,b_1,q_1,a_1}.$$

*Proof.* First, note that due to Lemma 5 (iii) and Lemma 11, the space $\bar{A}^{\mathcal{L}}_{0,r,B_0,q_0,a_0}$ is an intermediate space for the couple $(A_0, A_1)$. Denote $\bar{K}(t,f) = K\left(t,f;\bar{A}^{\mathcal{L}}_{0,r_0,b_0,q_0,a_0}, \bar{A}^{\mathcal{L}}_{1,r_1,b_1,q_1,a_1}\right)$ and $Z = \left(\bar{A}^{\mathcal{L}}_{0,r_0,b_0,q_0,a_0}, \bar{A}^{\mathcal{L}}_{1,r_1,b_1,q_1,a_1}\right)_{\theta,r,a}$ $(0 \le \theta \le 1)$. Using the change of variables and Theorem 30, we can write

$$\|f\|_Z \approx \left\|\rho(t)^{-\theta}t^{-\frac{1}{r}}a(\rho(t))\bar{K}(\rho(t),f)\right\|_{r,(0,\infty)}$$
$$= \left\|t^{-\theta-\frac{1}{r}}\left(\frac{\chi_1(t)}{\chi_0(t)}\right)^\theta a(\rho(t))\bar{K}(\rho(t),f)\right\|_{r,(0,\infty)} \approx I_1(\theta) + I_2(\theta) + I_3(\theta) + I_4(\theta),$$

where

$$I_1(\theta) := \left\|t^{-\theta-\frac{1}{r}}\left(\frac{\chi_1(t)}{\chi_0(t)}\right)^\theta a(\rho(t))\left\|s^{-1/r_0}b_0(s)\left\|u^{-1/q_0}a_0(u)K(u,f)\right\|_{q_0,(0,s)}\right\|_{r_0,(0,t)}\right\|_{r,(0,\infty)},$$

$$I_2(\theta) := \left\|t^{-\theta-\frac{1}{r}}\left(\frac{\chi_1(t)}{\chi_0(t)}\right)^\theta a(\rho(t))\left\|s^{-1/r_0}b_0(s)\right\|_{r_0,(t,\infty)}\left\|u^{-1/q_0}a_0(u)K(u,f)\right\|_{q_0,(0,t)}\right\|_{r,(0,\infty)},$$

$I_3(\theta)$
$$:= \left\|t^{-\theta-\frac{1}{r}}\left(\frac{\chi_1(t)}{\chi_0(t)}\right)^\theta a(\rho(t))\rho(t)\left\|s^{-\frac{1}{r_1}}b_1(s)\left\|u^{-1-\frac{1}{q_1}}a_1(u)K(u,f)\right\|_{q_1,(t,s)}\right\|_{r_1,(t,\infty)}\right\|_{r,(0,\infty)}$$
$$= \left\|t^{1-\theta-\frac{1}{r}}\left(\frac{\chi_1(t)}{\chi_0(t)}\right)^{\theta-1} a(\rho(t))\left\|s^{-\frac{1}{r_1}}b_1(s)\left\|u^{-1-\frac{1}{q_1}}a_1(u)K(u,f)\right\|_{q_1,(t,s)}\right\|_{r_1,(t,\infty)}\right\|_{r,(0,\infty)},$$

and
$$I_4(\theta) := \left\|t^{-\theta-\frac{1}{r}}(\chi_0(t))^{1-\theta}(\chi_1(t))^\theta a(\rho(t))K(t,f)\right\|_{r,(0,\infty)} = \|f\|_{\theta,r;a^\#}.$$

Consider $I_1(\theta)$. Obviously,
$$I_1(0) = \left\|t^{-1/r}a(\rho(t))\left\|s^{-1/r_0}b_0(s)\left\|u^{-1/q_0}a_0(u)K(u,f)\right\|_{q_0,(0,s)}\right\|_{r_0,(0,t)}\right\|_{r,(0,\infty)}$$
$$= \|f\|_{\bar{A}^{\mathcal{L},\mathcal{L}}_{0,r,a^\circ\rho,r_0,b_0,q_0,a_0}}.$$

Let $\theta > 0$. By Lemma 22, we have
$$I_1(\theta) \approx \left\|t^{-\theta-1/r}\chi_0(t)^{-\theta}\chi_1(t)^\theta a(\rho(t))b_0(t)a_0(t)K(t,f)\right\|_{r,(0,\infty)}.$$

By Lemma 3,
$$b_0(t)a_0(t) \prec \left\|s^{-1/r_0}b_0(s)\left\|u^{-1/q_0}a_0(u)\right\|_{q_0,(t,s)}\right\|_{r_0,(t,\infty)} = \chi_0(t).$$

Hence,
$$\chi_0(t)^{-\theta}\chi_1(t)^\theta a(\rho(t))b_0(t)a_0(t) \prec (\chi_0(t))^{1-\theta}(\chi_1(t))^\theta a(\rho(t))$$
and
$$I_1(\theta) \prec I_4(\theta).$$

Consider $I_2(\theta)$. We observe,
$$I_2(0) = \left\|t^{-1/r}a(\rho(t))\left\|s^{-1/r_0}b_0(s)\right\|_{r_0,(t,\infty)}\left\|u^{-1/q_0}a_0(u)K(u,f)\right\|_{q_0,(0,t)}\right\|_{r,(0,\infty)}$$

$$= \|f\|_{\mathcal{L};0,r,B_0,q_0,a_0},$$
where $B_0(t) = a(\rho(t)) \|s^{-1/r_0} b_0(s)\|_{r_0,(t,\infty)}$. Let $\theta > 0$. Lemma 21, implies
$$I_2(\theta) \approx \left\| t^{-\theta-1/r} \left(\frac{\chi_1(t)}{\chi_0(t)}\right)^\theta a(\rho(t)) \|s^{-1/r_0} b_0(s)\|_{r_0,(t,\infty)} a_0(t) K(t,f) \right\|_{r,(0,\infty)}.$$
By Lemma 3,
$$\|s^{-1/r_0} b_0(s)\|_{r_0,(t,\infty)} a_0(t) \prec \left\| \|s^{-1/r_0} b_0(s)\|u^{-1/q_0} a_0(u)\|_{q_0,(t,s)} \right\|_{r_0,(t,\infty)} = \chi_0(t).$$
Hence,
$$\chi_0(t)^{-\theta} \chi_1(t)^\theta a(\rho(t)) \|s^{-1/r_0} b_0(s)\|_{r_0,(t,\infty)} a_0(t) \prec (\chi_0(t))^{1-\theta} (\chi_1(t))^\theta a(\rho(t))$$
and
$$I_2(\theta) \prec I_4(\theta).$$
Consider $I_3(\theta)$. If $0 \le \theta < 1$, by Lemma 23 (ii) and (17), we get
$$I_3(\theta) \approx \left\| t^{-\theta-1/r} \left(\frac{\chi_1(t)}{\chi_0(t)}\right)^{\theta-1} a(\rho(t)) \|s^{-1/r_1} b_1(s)\|_{r_1,(t,\infty)} a_1(t) K(t,f) \right\|_{r,(0,\infty)}$$
$$\prec \left\| t^{-\theta-1/r} \left(\frac{\chi_1(t)}{\chi_0(t)}\right)^{\theta-1} a(\rho(t)) \chi_1(t) K(t,f) \right\|_{r,(0,\infty)} = I_4(\theta).$$
Thus, the statements (i) and (ii) are proved. Furthermore,
$$I_3(1) = \left\| t^{-1/r} a(\rho(t)) \left\| \|s^{-1/r_1} b_1(s)\|u^{-1-1/q_1} a_1(u) K(u,f)\|_{q_1,(t,s)} \right\|_{r_1,(t,\infty)} \right\|_{r,(0,\infty)}$$
$$= \|f\|_{\mathcal{L},\mathcal{R};1,r,a\circ\rho,r_1,b_1,q_1,a_1}.$$
Thus, the statement (iii) is also proved. □

**Remark 32.** Note that Theorem 26 and Theorem 27 can be deduced from Theorem 30 and Theorem 31, respectively. Really, taking $q_0 = \infty$ and $a_0 = 1$ in Theorem 30 and using Lemma 13 and Lemma 15, we get Theorem 26. Similarly, Theorem 41 and Theorem 42 can be deduced from Theorem 45 and Theorem 46, correspondently. See below.

**Theorem 33.** Let $0 < \theta_0 < 1$, $0 < r_0, q_0 \le \infty$, $a_0, b_0 \in SV$, and $\|s^{-1/r_0} b_0(s)\|_{r_0,(1,\infty)} < \infty$. Put $\chi_0(t) = a_0(t) \|s^{-1/r_0} b_0(s)\|_{r_0,(t,\infty)}$ and $\rho(t) = t^{1-\theta_0} \frac{\chi_0(t)}{\chi_1(t)}$. Then, for all $f \in \bar{A}^{\mathcal{L}}_{\theta_0,r_0,b_0,q_0,a_0} + \bar{A}^{\mathcal{L}}_{1,r_1,b_1,q_1,a_1}$ and $t > 0$,
$$K\big(\rho(t), f; \bar{A}^{\mathcal{L}}_{\theta_0,r_0,b_0,q_0,a_0}, \bar{A}^{\mathcal{L}}_{1,r_1,b_1,q_1,a_1}\big)$$
$$\approx \left\| \|s^{-1/r_0} b_0(s)\|u^{-\theta_0-1/q_0} a_0(u) K(u,f)\|_{q_0,(0,s)} \right\|_{r_0,(0,t)}$$
$$+ \|s^{-1/r_0} b_0(s)\|_{r_0,(t,\infty)} \|u^{-\theta_0-1/q_0} a_0(u) K(u,f)\|_{q_0,(0,t)}$$
$$+ \rho(t) \left\| \|s^{-1/r_1} b_1(s)\|u^{-1-1/q_1} a_1(u) K(u,f)\|_{q_1,(t,s)} \right\|_{r_1,(t,\infty)}.$$

*Proof.* Let $\Phi_0$ be the function space corresponding to $\bar{A}^{\mathcal{L}}_{\theta_0,r_0,b_0,q_0,a_0}$:
$$\|F\|_{\Phi_0} = \left\| \|s^{-1/r_0} b_0(s)\|u^{-\theta_0-1/q_0} a_0(u) F(u)\|_{q_0,(0,s)} \right\|_{r_0,(0,\infty)}.$$
We know (see the proof of Theorem 22 [9]) that
$$I(t,f) \approx \left\| \|s^{-1/r_0} b_0(s)\|u^{-\theta_0-1/q_0} a_0(u) K(u,f)\|_{q_0,(0,s)} \right\|_{r_0,(0,t)}$$
$$+ \|s^{-1/r_0} b_0(s)\|_{r_0,(t,\infty)} \|u^{-\theta_0-1/q_0} a_0(u) K(u,f)\|_{q_0,(0,t)}$$
and
$$g_0(t) \prec h_0(t) \approx g_0(t) + h_0(t) \approx t^{1-\theta_0} \chi_0(t).$$

Now we estimate the functions $g(t)$ and $h(t)$.

$$\frac{1}{g(t)} := \left\|\chi_{(t,\infty)}(u)\frac{u}{g_0(u)+h_0(u)}\right\|_{\Phi_1}$$

$$\approx \left\|s^{-1/r_1}b_1(s)\left\|u^{\theta_0-1-1/q_1}\frac{a_1(u)}{\chi_0(u)}\right\|_{q_1,(t,s)}\right\|_{r_1,(t,\infty)}$$

$$< \left\|s^{-1/r_1}b_1(s)\left\|u^{\theta_0-1-1/q_1}\frac{a_1(u)}{\chi_0(u)}\right\|_{q_1,(t,\infty)}\right\|_{r_1,(t,\infty)}$$

$$\approx t^{\theta_0-1}\frac{a_1(t)}{\chi_0(t)}\left\|s^{-1/r_1}b_1(s)\right\|_{r_1,(t,\infty)} \prec t^{\theta_0-1}\frac{\chi_1(t)}{\chi_0(t)} = \frac{1}{\rho(t)}.$$

The last inequality follows from (17). Hence, $\rho(t) \prec g(t)$.

$$h(t) := \left\|\chi_{(0,t)}(u)\frac{u}{g_1(u)+h_1(u)}\right\|_{\Phi_0}$$

$$\approx \left\|s^{-1/r_0}b_0(s)\left\|u^{1-\theta_0-1/q_0}\frac{a_0(u)}{\chi_1(u)}\right\|_{q_0,(0,\min(t,s))}\right\|_{r_0,(0,\infty)}$$

$$\approx \left\|s^{-1/r_0}b_0(s)\left\|u^{1-\theta_0-1/q_0}\frac{a_0(u)}{\chi_1(u)}\right\|_{q_0,(0,s)}\right\|_{r_0,(0,t)}$$

$$+ \left\|s^{-1/r_0}b_0(s)\left\|u^{1-\theta_0-1/q_0}\frac{a_0(u)}{\chi_1(u)}\right\|_{q_0,(0,t)}\right\|_{r_0,(t,\infty)}$$

$$\approx \left\|s^{1-\theta_0-1/r_0}b_0(s)\frac{a_0(s)}{\chi_1(s)}\right\|_{r_0,(0,t)} + t^{1-\theta_0}\frac{a_0(t)}{\chi_1(t)}\left\|s^{-1/r_0}b_0(s)\right\|_{r_0,(t,\infty)}$$

$$\approx t^{1-\theta_0}b_0(t)\frac{a_0(t)}{\chi_1(t)} + t^{1-\theta_0}\frac{a_0(t)}{\chi_1(t)}\left\|s^{-1/r_0}b_0(s)\right\|_{r_0,(t,\infty)}$$

$$\approx t^{1-\theta_0}\frac{a_0(t)}{\chi_1(t)}\left\|s^{-1/r_0}b_0(s)\right\|_{r_0,(t,\infty)} = \rho(t).$$

Finally, note that

$$\frac{h_0(t)}{h_1(t)} \approx \frac{t^{1-\theta_0}\chi_0(t)}{\chi_1(t)} = \rho(t).$$

Theorem 4 (Case 4) of [2] completes the proof. □

**Theorem** 34. Let $0 < \theta_0 < 1$, $0 < r, r_0, q_0 \leq \infty$, $a, a_0, b_0 \in SV$, and $\left\|s^{-1/r_0}b_0(s)\right\|_{r_0,(1,\infty)} < \infty$. Put $\chi_0(t) = a_0(t)\left\|s^{-1/r_0}b_0(s)\right\|_{r_0,(t,\infty)}$ and $\rho(t) = t^{1-\theta_0}\frac{\chi_0(t)}{\chi_1(t)}$. For given $0 \leq \theta \leq 1$, denote $a^\# = \chi_0^{(1-\theta)}\chi_1^\theta a \circ \rho$ and $\eta = (1-\theta)\theta_0 + \theta$.

(i) If $0 < \theta < 1$, then
$$\left(\bar{A}^{\mathcal{L}}_{\theta_0,r_0,b_0,q_0,a_0}, \bar{A}^{\mathcal{L}}_{1,r_1,b_1,q_1,a_1}\right)_{\theta,r,a} = \bar{A}_{\eta,r,a^\#}.$$

(ii) If $\left\|s^{-1/r}a(s)\right\|_{r,(1,\infty)} < \infty$, then
$$\left(\bar{A}^{\mathcal{L}}_{\theta_0,r_0,b_0,q_0,a_0}, \bar{A}^{\mathcal{L}}_{1,r_1,b_1,q_1,a_1}\right)_{0,r,a} = \bar{A}^{\mathcal{L}}_{\theta_0,r,B_0,q_0,a_0} \cap \bar{A}^{\mathcal{L},\mathcal{L}}_{\theta_0,r,a\circ\rho,r_0,b_0,q_0,a_0},$$
where $B_0(t) = \left\|s^{-1/r_0}b_0(s)\right\|_{r_0,(t,\infty)}a(\rho(t))$.

(iii) If $\left\|s^{-1/r}a(s)\right\|_{r,(0,1)} < \infty$, then
$$\left(\bar{A}^{\mathcal{L}}_{\theta_0,r_0,b_0,q_0,a_0}, \bar{A}^{\mathcal{L}}_{1,r_1,b_1,q_1,a_1}\right)_{1,r,a} = \bar{A}_{1,r,a^\#} \cap \bar{A}^{\mathcal{L},\mathcal{R}}_{1,r,a\circ\rho,r_1,b_1,q_1,a_1}.$$

*Proof.* First, note that due to Lemma 5 (i) and Lemma 11 the space $\bar{A}^{\mathcal{L}}_{\theta_0,r,B_0,q_0,a_0}$ is an intermediate space for the couple $(A_0, A_1)$. Denote $\bar{K}(t,f) = K\big(t,f; \bar{A}^{\mathcal{L}}_{\theta_0,r_0,b_0,q_0,a_0}, \bar{A}^{\mathcal{L}}_{1,r_1,b_1,q_1,a_1}\big)$, and $Z = \big(\bar{A}^{\mathcal{L}}_{\theta_0,r_0,b_0,q_0,a_0}, \bar{A}^{\mathcal{L}}_{1,r_1,b_1,q_1,a_1}\big)_{\theta,r,a}$ $(0 \leq \theta \leq 1)$.
Using the change of variables and Theorem 33, we can write
$$\|f\|_Z \approx \left\|\rho(t)^{-\theta} t^{-1/r} a(\rho(t)) \bar{K}(\rho(t),f)\right\|_{r,(0,\infty)}$$
$$= \left\|t^{-\theta(1-\theta_0)-1/r} \left(\frac{\chi_1(t)}{\chi_0(t)}\right)^\theta a(\rho(t)) \bar{K}(\rho(t),f)\right\|_{r,(0,\infty)}$$
$$= \left\|t^{\theta_0-\eta-1/r} \left(\frac{\chi_1(t)}{\chi_0(t)}\right)^\theta a(\rho(t)) \bar{K}(\rho(t),f)\right\|_{r,(0,\infty)} \approx I_1(\theta) + I_2(\theta) + I_3(\theta),$$
where
$$I_1(\theta) :=$$
$$\left\|t^{\theta_0-\eta-\frac{1}{r}} \left(\frac{\chi_1(t)}{\chi_0(t)}\right)^\theta a(\rho(t)) \left\|s^{-\frac{1}{r_0}} b_0(s) \left\|u^{-\theta_0-\frac{1}{q_0}} a_0(u) K(u,f)\right\|_{q_0,(0,s)}\right\|_{r_0,(0,t)}\right\|_{r,(0,\infty)},$$
$$I_2(\theta) :=$$
$$\left\|t^{\theta_0-\eta-\frac{1}{r}} \left(\frac{\chi_1(t)}{\chi_0(t)}\right)^\theta a(\rho(t)) \left\|s^{-\frac{1}{r_0}} b_0(s)\right\|_{r_0,(t,\infty)} \left\|u^{-\theta_0-\frac{1}{q_0}} a_0(u) K(u,f)\right\|_{q_0,(0,t)}\right\|_{r,(0,\infty)},$$
and
$$I_3(\theta)$$
$$:= \left\|t^{\theta_0-\eta-\frac{1}{r}} \left(\frac{\chi_1(t)}{\chi_0(t)}\right)^\theta a(\rho(t))\rho(t) \left\|s^{-\frac{1}{r_1}} b_1(s) \left\|u^{-1-\frac{1}{q_1}} a_1(u) K(u,f)\right\|_{q_1,(t,s)}\right\|_{r_1,(t,\infty)}\right\|_{r,(0,\infty)}$$
$$= \left\|t^{1-\eta-\frac{1}{r}} \left(\frac{\chi_1(t)}{\chi_0(t)}\right)^{\theta-1} a(\rho(t)) \left\|s^{-\frac{1}{r_1}} b_1(s) \left\|u^{-1-\frac{1}{q_1}} a_1(u) K(u,f)\right\|_{q_1,(t,s)}\right\|_{r_1,(t,\infty)}\right\|_{r,(0,\infty)}.$$
Consider $I_1(\theta)$. We observe,
$$I_1(0) = \left\|t^{-1/r} a(\rho(t)) \left\|s^{-1/r_0} b_0(s) \left\|u^{-\theta_0-1/q_0} a_0(u) K(u,f)\right\|_{q_0,(0,s)}\right\|_{r_0,(0,t)}\right\|_{r,(0,\infty)}$$
$$= \|f\|_{\bar{A}^{\mathcal{L},\mathcal{L}}_{\theta_0,r,a\circ\rho,r_0,b_0,q_0,a_0}}.$$
Let $\theta > 0$. Because $b_0(t)a_0(t) \prec \chi_0(t)$, Lemma 22 implies
$$I_1(\theta) \approx \left\|t^{-\eta-1/r} \left(\frac{\chi_1(t)}{\chi_0(t)}\right)^\theta a(\rho(t)) b_0(t) a_0(t) K(t,f)\right\|_{r,(0,\infty)}$$
$$\prec \left\|t^{-\eta-1/r} (\chi_0(t))^{1-\theta} (\chi_1(t))^\theta a(\rho(t)) K(t,f)\right\|_{r,(0,\infty)} = \|f\|_{\eta,r;a^\#}.$$
Consider $I_2(\theta)$. First, notice that
$$I_2(0) := \left\|t^{-1/r} a(\rho(t)) \left\|s^{-1/r_0} b_0(s)\right\|_{r_0,(t,\infty)} \left\|u^{-\theta_0-1/q_0} a_0(u) K(u,f)\right\|_{q_0,(0,t)}\right\|_{r,(0,\infty)}$$
$$= \|f\|_{\mathcal{L};\theta_0,r,B_0,q_0,a_0},$$
where $B_0(t) = a(\rho(t)) \|s^{-1/r_0} b_0(s)\|_{r_0,(t,\infty)}$. Let $\theta > 0$. Lemma 21 implies
$$I_2(\theta) \approx \left\|t^{-\eta-1/r} \left(\frac{\chi_1(t)}{\chi_0(t)}\right)^\theta a(\rho(t)) \|s^{-1/r_0} b_0(s)\|_{r_0,(t,\infty)} a_0(t) K(t,f)\right\|_{r,(0,\infty)}$$
$$= \left\|t^{-\eta-1/r} (\chi_0(t))^{1-\theta} (\chi_1(t))^\theta a(\rho(t)) K(t,f)\right\|_{r,(0,\infty)} = \|f\|_{\eta,r;a^\#}.$$

Consider $I_3(\theta)$. If $0 < \theta < 1$, by Lemma 23 (ii) and (17), we get

$$I_3(\theta) \approx \left\| t^{-\eta-1/r} \left(\frac{\chi_1(t)}{\chi_0(t)}\right)^{\theta-1} a(\rho(t)) \left\| s^{-1/r_1} b_1(s) \right\|_{r_1,(t,\infty)} a_1(t) K(t,f) \right\|_{r,(0,\infty)}$$

$$\prec \left\| t^{-\eta-1/r} (\chi_0(t))^{1-\theta} (\chi_1(t))^{\theta} a(\rho(t)) K(t,f) \right\|_{r,(0,\infty)} = \|f\|_{\eta,r;a^\#}.$$

Thus, the statement (i) is proved. From Lemma 23 (ii), (17), and Lemma 16, it follows

$$I_3(0) \approx \left\| t^{-\theta_0-1/r} \frac{\chi_0(t)}{\chi_1(t)} a(\rho(t)) \left\| s^{-1/r_1} b_1(s) \right\|_{r_1,(t,\infty)} a_1(t) K(t,f) \right\|_{r,(0,\infty)}$$

$$\prec \left\| t^{-\theta_0-1/r} \chi_0(t) a(\rho(t)) K(t,f) \right\|_{r,(0,\infty)}$$

$$= \left\| t^{-\theta_0-1/r} a_0(t) \left\| s^{-1/r_0} b_0(s) \right\|_{r_0,(t,\infty)} a(\rho(t)) K(t,f) \right\|_{r,(0,\infty)}$$

$$\prec \left\| t^{-1/r} \left\| s^{-1/r_0} b_0(s) \right\|_{r_0,(t,\infty)} a(\rho(t)) \left\| u^{-\theta_0-1/q_0} a_0(u) K(u,f) \right\|_{q_0,(0,t)} \right\|_{r,(0,\infty)} = I_2(0).$$

So,
$$I_1(0) + I_2(0) + I_3(0) \approx I_1(0) + I_2(0).$$

This means that the statement (ii) is proved. Finally,

$$I_3(1) = \left\| t^{-1/r} a(\rho(t)) \left\| s^{-1/r_1} b_1(s) \left\| u^{-1-1/q_1} a_1(u) K(u,f) \right\|_{q_1,(t,s)} \right\|_{r_1,(t,\infty)} \right\|_{r,(0,\infty)}$$

$$= \|f\|_{\mathcal{L},\mathcal{R};1,r,a\circ\rho,r_1,b_1,q_1,a_1}$$

Thus, the statement (iii) is also proved. □

**Theorem 35.** Let $0 < r_0, q_0 \leq \infty$, $a_0, b_0 \in SV$, $\left\| s^{-1/r_0} b_0(s) \right\|_{r_0,(0,1)} < \infty$, and $\left\| s^{-1/r_0} b_0(s) \left\| u^{-1/q_0} a_0(u) \right\|_{q_0,(s,\infty)} \right\|_{r_0,(1,\infty)} < \infty$. Put

$$\chi_0(t) = \left\| s^{-1/r_0} b_0(s) \left\| u^{-1/q_0} a_0(u) \right\|_{q_0,(s,\infty)} \right\|_{r_0,(t,\infty)}$$
$$+ \left\| s^{-1/r_0} b_0(s) \right\|_{r_0,(0,t)} \left\| u^{-1/q_0} a_0(u) \right\|_{q_0,(t,\infty)}$$

and $\rho(t) = t \frac{\chi_0(t)}{\chi_1(t)}$. Then, for all $f \in \bar{A}^{\mathcal{R}}_{0,r_0,b_0,q_0,a_0} + \bar{A}^{\mathcal{L}}_{1,r_1,b_1,q_1,a_1}$ and $t > 0$,

$$K(\rho(t), f; \bar{A}^{\mathcal{R}}_{0,r_0,b_0,q_0,a_0}, \bar{A}^{\mathcal{L}}_{1,r_1,b_1,q_1,a_1}) \approx \left\| s^{-1/r_0} b_0(s) \left\| u^{-1/q_0} a_0(u) K(u,f) \right\|_{q_0,(s,t)} \right\|_{r_0,(0,t)}$$

$$+ \rho(t) \left\| s^{-1/r_1} b_1(s) \left\| u^{-1-1/q_1} a_1(u) K(u,f) \right\|_{q_1,(t,s)} \right\|_{r_1,(t,\infty)} + \chi_0(t) K(t,f).$$

*Proof.* Let $\Phi_0$ be the function space corresponding to $\bar{A}^{\mathcal{R}}_{0,r_0,b_0,q_0,a_0}$:

$$\|F\|_{\Phi_0} = \left\| s^{-1/r_0} b_0(s) \left\| u^{-1/q_0} a_0(u) F(u) \right\|_{q_0,(s,\infty)} \right\|_{r_0,(0,\infty)}.$$

Let us estimate the functions $I(t,f)$, $g_0(t)$, and $h_0(t)$.

$$I(t,f) := \left\| \chi_{(0,t)}(*) K(*,f) \right\|_{\Phi_0}$$
$$= \left\| s^{-1/r_0} b_0(s) \left\| u^{-1/q_0} a_0(u)(u) K(u,f) \right\|_{q_0,(s,t)} \right\|_{r_0,(0,t)}.$$

$$g_0(t) := t \left\| \chi_{(t,\infty)}(*) \right\|_{\Phi_0} = t \left\| s^{-1/r_0} b_0(s) \left\| u^{-1/q_0} a_0(u) \right\|_{q_0,(max(t,s),\infty)} \right\|_{r_0,(t,\infty)}$$

$$\approx t \left( \left\| s^{-1/r_0} b_0(s) \right\|_{r_0,(0,t)} \left\| u^{-1/q_0} a_0(u) \right\|_{q_0,(t,\infty)} + \right.$$

$$\left. \left\| s^{-1/r_0} b_0(s) \left\| u^{-1/q_0} a_0(u) \right\|_{q_0,(s,\infty)} \right\|_{r_0,(t,\infty)} \right) = t \chi_0(t).$$

$$h_0(t) := \|* \chi_{(0,t)}(*)\|_{\Phi_0} = \left\| \|s^{-1/r_0} b_0(s)\| u^{1-1/q_0} a_0(u)\|_{q_0,(s,t)} \right\|_{r_0,(0,t)}$$

$$\leq \left\| \|s^{-1/r_0} b_0(s)\| u^{1-1/q_0} a_0(u)\|_{q_0,(0,t)} \right\|_{r_0,(0,t)}$$

$$\approx t a_0(t) \|s^{-1/r_0} b_0(s)\|_{r_0,(0,t)} \prec t \|u^{-1/q_0} a_0(u)\|_{q_0,(t,\infty)} \|s^{-1/r_0} b_0(s)\|_{r_0,(0,t)} \prec t\chi_0(t).$$

Thus,
$$h_0(t) \prec g_0(t) \approx h_0(t) + g_0(t) \approx t\chi_0(t).$$

Now we estimate the functions $g(t)$ and $h(t)$.

$$\frac{1}{g(t)} := \left\| \chi_{(t,\infty)}(u) \frac{u}{g_0(u) + h_0(u)} \right\|_{\Phi_1} \approx$$

$$\approx \left\| \|s^{-1/r_1} b_1(s)\| u^{-1-1/q_1} \frac{a_1(u)}{\chi_0(u)} \|_{q_1,(t,s)} \right\|_{r_1,(t,\infty)}$$

$$\prec \left\| \|s^{-1/r_1} b_1(s)\| u^{-1-1/q_1} \frac{a_1(u)}{\chi_0(u)} \|_{q_1,(t,\infty)} \right\|_{r_1,(t,\infty)} \approx t^{-1} \frac{a_1(t)}{\chi_0(t)} \|s^{-1/r_1} b_1(s)\|_{r_1,(t,\infty)}$$

Using (17), we get $\frac{1}{g(t)} \prec t^{-1} \frac{\chi_1(t)}{\chi_0(t)} = \frac{1}{\rho(t)}$ and $\rho(t) \prec g(t)$.

$$h(t) := \left\| \chi_{(0,t)}(u) \frac{u}{g_1(u) + h_1(u)} \right\|_{\Phi_0} \approx \left\| \|s^{-1/r_0} b_0(s)\| u^{1-1/q_0} \frac{a_0(u)}{\chi_1(u)} \|_{q_0,(s,t)} \right\|_{r_0,(0,t)}$$

$$\prec \left\| u^{1-1/q_0} \frac{a_0(u)}{\chi_1(u)} \right\|_{q_0,(0,t)} \|s^{-1/r_0} b_0(s)\|_{r_0,(0,t)}$$

$$\approx t \frac{a_0(t)}{\chi_1(t)} \|s^{-1/r_0} b_0(s)\|_{r_0,(0,t)} \prec t \frac{\|u^{-1/q_0} a_0(u)\|_{q_0,(t,\infty)}}{\chi_1(t)} \|s^{-1/r_0} b_0(s)\|_{r_0,(0,t)}$$

$$\prec t \frac{\chi_0(t)}{\chi_1(t)} = \rho(t).$$

Thus, $h(t) \prec \rho(t) \prec g(t)$. Theorem 4 (Case 2) of [2] completes the proof. $\square$

**Theorem** 36. Let $0 < r, r_0, q_0 \leq \infty$, $a, a_0, b_0 \in SV$, $\|s^{-1/r_0} b_0(s)\|_{r_0,(0,1)} < \infty$, and $\left\| \|s^{-1/r_0} b_0(s)\| u^{-1/q_0} a_0(u)\|_{q_0,(s,\infty)} \right\|_{r_0,(1,\infty)} < \infty$. Put

$$\chi_0(t) = \left\| \|s^{-1/r_0} b_0(s)\| u^{-1/q_0} a_0(u)\|_{q_0,(s,\infty)} \right\|_{r_0,(t,\infty)}$$

$$+ \|s^{-1/r_0} b_0(s)\|_{r_0,(0,t)} \|u^{-1/q_0} a_0(u)\|_{q_0,(t,\infty)}$$

and $\rho(t) = t \frac{\chi_0(t)}{\chi_1(t)}$. For given $0 \leq \theta \leq 1$, denote $a^{\#} = \chi_0^{(1-\theta)} \chi_1^\theta a \circ \rho$.

(i) If $0 < \theta < 1$, then
$$\left( \bar{A}^{\mathcal{R}}_{0,r_0,b_0,q_0,a_0}, \bar{A}^{\mathcal{L}}_{1,r_1,b_1,q_1,a_1} \right)_{\theta,r,a} = \bar{A}_{\theta,r,a^{\#}}.$$

(ii) If $\|s^{-1/r} a(s)\|_{r,(1,\infty)} < \infty$, then
$$\left( \bar{A}^{\mathcal{R}}_{0,r_0,b_0,q_0,a_0}, \bar{A}^{\mathcal{L}}_{1,r_1,b_1,q_1,a_1} \right)_{0,r,a} = \bar{A}_{0,r,a^{\#}} \cap \bar{A}^{\mathcal{R},\mathcal{L}}_{0,r,a\circ\rho,r_0,b_0,q_0,a_0}.$$

(iii) If $\|s^{-1/r} a(s)\|_{r,(0,1)} < \infty$, then
$$\left( \bar{A}^{\mathcal{R}}_{0,r_0,b_0,q_0,a_0}, \bar{A}^{\mathcal{L}}_{1,r_1,b_1,q_1,a_1} \right)_{1,r,a} = \bar{A}_{1,r,a^{\#}} \cap \bar{A}^{\mathcal{L},\mathcal{R}}_{1,r,a\circ\rho,r_1,b_1,q_1,a_1}.$$

*Proof.* Denote $\bar{K}(t,f) = K\left(t,f; \bar{A}^{\mathcal{R}}_{0,r_0,b_0,q_0,a_0}, \bar{A}^{\mathcal{L}}_{1,r_1,b_1,q_1,a_1}\right)$, and $Z = \bar{X}_{\theta,r,a}$ ($0 \leq \theta \leq 1$). Using the change of variables $u = \rho(t)$ and Theorem 35, we can write
$$\|f\|_Z \approx \left\| \rho(t)^{-\theta} t^{-1/r} a(\rho(t)) \bar{K}(\rho(t), f) \right\|_{r,(0,\infty)}$$

$$= \left\|t^{-\theta-1/r}\left(\frac{\chi_1(t)}{\chi_0(t)}\right)^\theta a(\rho(t))\overline{K}(\rho(t),f)\right\|_{r,(0,\infty)} \approx I_1(\theta) + I_2(\theta) + I_3(\theta),$$

where

$$I_1(\theta) := \left\|t^{-\theta-\frac{1}{r}}\left(\frac{\chi_1(t)}{\chi_0(t)}\right)^\theta a(\rho(t))\left\|s^{-\frac{1}{r_0}}b_0(s)\left\|u^{-\frac{1}{q_0}}a_0(u)K(u,f)\right\|_{q_0,(s,t)}\right\|_{r_0,(0,t)}\right\|_{r,(0,\infty)},$$

$$I_2(\theta)$$
$$:= \left\|t^{-\theta-\frac{1}{r}}\left(\frac{\chi_1(t)}{\chi_0(t)}\right)^\theta a(\rho(t))\,\rho(t)\left\|s^{-\frac{1}{r_1}}b_1(s)\left\|u^{-1-\frac{1}{q_1}}a_1(u)K(u,f)\right\|_{q_1,(t,s)}\right\|_{r_1,(t,\infty)}\right\|_{r,(0,\infty)}$$
$$= \left\|t^{1-\theta-\frac{1}{r}}\left(\frac{\chi_1(t)}{\chi_0(t)}\right)^{\theta-1} a(\rho(t))\left\|s^{-\frac{1}{r_1}}b_1(s)\left\|u^{-1-\frac{1}{q_1}}a_1(u)K(u,f)\right\|_{q_1,(t,s)}\right\|_{r_1,(t,\infty)}\right\|_{r,(0,\infty)},$$

and

$$I_3(\theta) := \left\|t^{-\theta-1/r}\left(\frac{\chi_1(t)}{\chi_0(t)}\right)^\theta a(\rho(t))\chi_0(t)K(t,f)\right\|_{r,(0,\infty)} = \|f\|_{\theta,r;a^\#}.$$

Consider $I_1(\theta)$. We observe,

$$I_1(0) = \left\|t^{-1/r}a(\rho(t))\left\|s^{-1/r_0}b_0(s)\|u^{-1/q_0}a_0(u)K(u,f)\|_{q_0,(s,t)}\right\|_{r_0,(0,t)}\right\|_{r,(0,\infty)}$$
$$= \|f\|_{\mathcal{R},\mathcal{L};0,r,a\circ\rho,r_0,b_0,q_0,a_0}.$$

Let $\theta > 0$. Note that by Lemma 2 (iv), we have

$$a_0(t)\|s^{-1/r_0}b_0(s)\|_{r_0,(0,t)} \prec \|s^{-1/r_0}b_0(s)\|_{r_0,(0,t)}\|u^{-1/q_0}a_0(u)\|_{q_0,(t,\infty)} \prec \chi_0(t).$$

Hence, by Lemma 23 (i), we get

$$I_1(\theta) \approx \left\|t^{-\theta-1/r}\chi_0(t)^{-\theta}\chi_1(t)^\theta a(\rho(t))\|s^{-1/r_0}b_0(s)\|_{r_0,(0,t)}a_0(t)K(t,f)\right\|_{r,(0,\infty)}$$
$$\prec \left\|t^{-\theta-1/r}\chi_0(t)^{1-\theta}\chi_1(t)^\theta a(\rho(t))K(t,f)\right\|_{r,(0,\infty)} = \|f\|_{\theta,r,a^\#}.$$

Consider $I_2(\theta)$. If $\theta < 1$, by Lemma 23 (ii) and (17), we obtain

$$I_2(\theta) \approx \left\|t^{-\theta-1/r}\chi_0(t)^{1-\theta}\chi_1(t)^{\theta-1}a(\rho(t))\|s^{-1/r_1}b_1(s)\|_{r_1,(t,\infty)}a_1(t)K(t,f)\right\|_{r,(0,\infty)}$$
$$\prec \left\|t^{-\theta-1/r}\chi_0(t)^{1-\theta}\chi_1(t)^\theta a(\rho(t))K(t,f)\right\|_{r,(0,\infty)} = \|f\|_{\theta,r,a^\#}.$$

Thus, the statements (i) and (ii) are proved. Finally,

$$I_2(1) = \left\|t^{-1/r}a(\rho(t))\left\|s^{-1/r_1}b_1(s)\|u^{-1-1/q_1}a_1(u)K(u,f)\|_{q_1,(t,s)}\right\|_{r_1,(t,\infty)}\right\|_{r,(0,\infty)}$$
$$= \|f\|_{\mathcal{L},\mathcal{R};1,r,a\circ\rho,r_1,b_1,q_1,a_1}.$$

Hence, the statement (iii) is also proved. □

**Theorem** 37. Let $0 < \theta_0 < 1$, $0 < r_0, q_0 \leq \infty$, $a_0, b_0 \in SV$, and $\|s^{-1/r_0}b_0(s)\|_{r_0,(0,1)} < \infty$. Put

$$\chi_0(t) = a_0(t)\|s^{-1/r_0}b_0(s)\|_{r_0,(0,t)}$$

and $\rho(t) = t^{1-\theta_0}\frac{\chi_0(t)}{\chi_1(t)}$. Then, for all $f \in \bar{A}^{\mathcal{R}}_{\theta_0,r_0,b_0,q_0,a_0} + \bar{A}^{\mathcal{L}}_{1,r_1,b_1,q_1,a_1}$ and $t > 0$,

$$K\left(\rho(t),f;\bar{A}^{\mathcal{R}}_{\theta_0,r_0,b_0,q_0,a_0},\bar{A}^{\mathcal{L}}_{1,r_1,b_1,q_1,a_1}\right) \approx$$
$$= \left\|s^{-1/r_0}b_0(s)\|u^{-\theta_0-1/q_0}a_0(u)K(u,f)\|_{q_0,(s,t)}\right\|_{r_0,(0,t)} +$$
$$+ \rho(t)\left\|s^{-1/r_1}b_1(s)\|u^{-1-1/q_1}a_1(u)K(u,f)\|_{q_1,(t,s)}\right\|_{r_1,(t,\infty)}. \tag{23}$$

*Proof.* Let $\Phi_0$ be the function space corresponding to $\bar{A}^{\mathcal{R}}_{\theta_0,r_0,b_0,q_0,a_0}$:
$$\|F\|_{\Phi_0} = \left\| s^{-1/r_0} b_0(s) \left\| u^{-\theta_0-1/q_0} a_0(u) F(u) \right\|_{q_0,(s,\infty)} \right\|_{r_0,(0,\infty)}.$$
We know (see the proof of Theorem 37 [9]) that
$$I(t,f) = \left\| s^{-1/r_0} b_0(s) \left\| u^{-\theta_0-1/q_0} a_0(u) K(u,f) \right\|_{q_0,(s,t)} \right\|_{r_0,(0,t)}$$
and
$$h_0(t) \prec g_0(t) \approx g_0(t) + h_0(t) \approx t^{1-\theta_0} \chi_0(t).$$
Now we estimate the functions $g(t)$ and $h(t)$.
$$\frac{1}{g(t)} := \left\| \chi_{(t,\infty)}(u) \frac{u}{g_0(u) + h_0(u)} \right\|_{\Phi_1}$$
$$\approx \left\| s^{-1/r_1} b_1(s) \left\| u^{\theta_0-1-1/q_1} \frac{a_1(u)}{\chi_0(u)} \right\|_{q_1,(t,s)} \right\|_{r_1,(t,\infty)}$$
$$\prec \left\| s^{-1/r_1} b_1(s) \left\| u^{\theta_0-1-1/q_1} \frac{a_1(u)}{\chi_0(u)} \right\|_{q_1,(t,\infty)} \right\|_{r_1,(t,\infty)} \approx t^{\theta_0-1} \frac{a_1(t)}{\chi_0(t)} \left\| s^{-1/r_1} b_1(s) \right\|_{r_1,(t,\infty)}.$$
By (17), we get $\frac{1}{g(t)} \prec t^{\theta_0-1} \frac{\chi_1(t)}{\chi_0(t)} = \frac{1}{\rho(t)}$ and hence, $\rho(t) \prec g(t)$.
$$h(t) := \left\| \chi_{(0,t)}(u) \frac{u}{g_1(u) + h_1(u)} \right\|_{\Phi_0} \approx \left\| s^{-1/r_0} b_0(s) \left\| u^{1-\theta_0-1/q_0} \frac{a_0(u)}{\chi_1(u)} \right\|_{q_0,(s,t)} \right\|_{r_0,(0,t)}$$
$$\prec \left\| u^{1-\theta_0-1/q_0} \frac{a_0(u)}{\chi_1(u)} \right\|_{q_0,(0,t)} \left\| s^{-1/r_0} b_0(s) \right\|_{r_0,(0,t)}$$
$$\approx t^{1-\theta_0} \frac{a_0(t)}{\chi_1(t)} \left\| s^{-1/r_0} b_0(s) \right\|_{r_0,(0,t)} = t^{1-\theta_0} \frac{\chi_0(t)}{\chi_1(t)} = \rho(t).$$
Thus, $h(t) \prec \rho(t) \prec g(t)$. By [2, Theorem 4, Case 2] we get
$$K\left(\rho(t), f; \bar{A}^{\mathcal{R}}_{\theta_0,r_0,b_0,q_0,a_0}, \bar{A}^{\mathcal{L}}_{1,r_1,b_1,q_1,a_1}\right)$$
$$\approx \left\| s^{-1/r_0} b_0(s) \left\| u^{-\theta_0-1/q_0} a_0(u) K(u,f) \right\|_{q_0,(s,t)} \right\|_{r_0,(0,t)}$$
$$+ \rho(t) \left\| s^{-1/r_1} b_1(s) \left\| u^{-1-1/q_1} a_1(u) K(u,f) \right\|_{q_1,(t,s)} \right\|_{r_1,(t,\infty)} + t^{-\theta_0} \chi_0(t) K(t,f).$$
By (9), we have
$$t^{-\theta_0} \chi_0(t) K(t,f) = t^{-\theta_0} a_0(t) \left\| s^{-1/r_0} b_0(s) \right\|_{r_0,(0,t)} K(t,f)$$
$$\prec \left\| s^{-1/r_0} b_0(s) \left\| u^{-\theta_0-1/q_0} a_0(u) K(u,f) \right\|_{q_0,(s,t)} \right\|_{r_0,(0,t)}. \tag{24}$$
So, we arrive at (23). $\square$

**Theorem 38.** Let $0 < \theta_0 < 1$, $0 < r, r_0, q_0 \le \infty$, $a, a_0, b_0 \in SV$, and $\left\| s^{-1/r_0} b_0(s) \right\|_{r_0,(0,1)} < \infty$. Put
$$\chi_0(t) = a_0(t) \left\| s^{-1/r_0} b_0(s) \right\|_{r_0,(0,t)}$$
and $\rho(t) = t^{1-\theta_0} \frac{\chi_0(t)}{\chi_1(t)}$. For given $0 \le \theta \le 1$, denote $a^\# = \chi_0^{(1-\theta)} \chi_1^\theta a \circ \rho$ and $\eta = (1-\theta)\theta_0 + \theta$.

(i) If $0 < \theta < 1$, then
$$\left( \bar{A}^{\mathcal{R}}_{\theta_0,r_0,b_0,q_0,a_0}, \bar{A}^{\mathcal{L}}_{1,r_1,b_1,q_1,a_1} \right)_{\theta,r,a} = \bar{A}_{\eta,r,a^\#}.$$

(ii) If $\left\| s^{-1/r} a(s) \right\|_{r,(1,\infty)} < \infty$, then
$$\left( \bar{A}^{\mathcal{R}}_{\theta_0,r_0,b_0,q_0,a_0}, \bar{A}^{\mathcal{L}}_{1,r_1,b_1,q_1,a_1} \right)_{0,r,a} = \bar{A}^{\mathcal{R},\mathcal{L}}_{\theta_0,r,a \circ \rho,r_0,b_0,q_0,a_0}.$$

(iii) If $\left\|s^{-1/r}a(s)\right\|_{r,(0,1)} < \infty$, then
$$\left(\bar{A}^{\mathcal{R}}_{\theta_0,r_0,b_0,q_0,a_0}, \bar{A}^{\mathcal{L}}_{1,r_1,b_1,q_1,a_1}\right)_{1,r,\mathrm{a}} = \bar{A}_{1,r,a^\#} \cap \bar{A}^{\mathcal{L},\mathcal{R}}_{1,r,a^\circ\rho,r_1,b_1,q_1,a_1}.$$

*Proof.* Denote $\bar{K}(t,f) = K\left(t,f; \bar{A}^{\mathcal{R}}_{\theta_0,r_0,b_0,q_0,a_0}, \bar{A}^{\mathcal{L}}_{1,r_1,b_1,q_1,a_1}\right)$ and $Z = \bar{X}_{\theta,r,\mathrm{a}}$ ($0 \leq \theta \leq 1$). Using the change of variables $u = \rho(t)$ and Theorem 37, we can write
$$\|f\|_Z \approx \left\|\rho(t)^{-\theta} t^{-1/r} a(\rho(t)) \bar{K}(\rho(t), f)\right\|_{r,(0,\infty)}$$
$$= \left\|t^{-\theta(1-\theta_0)-1/r} \left(\frac{\chi_1(t)}{\chi_0(t)}\right)^\theta a(\rho(t)) \bar{K}(\rho(t), f)\right\|_{r,(0,\infty)} \approx I_1(\theta) + I_2(\theta),$$
where
$$I_1(\theta)$$
$$:= \left\|t^{-\theta(1-\theta_0)-\frac{1}{r}} \left(\frac{\chi_1(t)}{\chi_0(t)}\right)^\theta a(\rho(t)) \left\|s^{-\frac{1}{r_0}} b_0(s) \left\|u^{-\theta_0-\frac{1}{q_0}} a_0(u) K(u,f)\right\|_{q_0,(s,t)}\right\|_{r_0,(0,t)}\right\|_{r,(0,\infty)}$$
and
$$I_2(\theta)$$
$$:= \left\|t^{-\theta(1-\theta_0)-\frac{1}{r}} \left(\frac{\chi_1(t)}{\chi_0(t)}\right)^\theta a(\rho(t)) \rho(t) \left\|s^{-\frac{1}{r_1}} b_1(s) \left\|u^{-1-\frac{1}{q_1}} a_1(u) K(u,f)\right\|_{q_1,(t,s)}\right\|_{r_1,(t,\infty)}\right\|_{r,(0,\infty)}$$
$$= \left\|t^{(1-\theta)(1-\theta_0)-\frac{1}{r}} \left(\frac{\chi_1(t)}{\chi_0(t)}\right)^{\theta-1} a(\rho(t)) \left\|s^{-\frac{1}{r_1}} b_1(s) \left\|u^{-1-\frac{1}{q_1}} a_1(u) K(u,f)\right\|_{q_1,(t,s)}\right\|_{r_1,(t,\infty)}\right\|_{r,(0,\infty)}$$

Consider the case $\theta = 0$. We have
$$I_1(0) = \left\|t^{-1/r} a(\rho(t)) \left\|s^{-1/r_0} b_0(s) \left\|u^{-\theta_0-1/q_0} a_0(u) K(u,f)\right\|_{q_0,(s,t)}\right\|_{r_0,(0,t)}\right\|_{r,(0,\infty)} =$$
$$\|f\|_{\mathcal{R},\mathcal{L};\theta_0,r,a^\circ\rho,r_0,b_0,q_0,a_0}.$$

By Lemma 23 (ii) and (17), we get
$$I_2(0) = \left\|t^{1-\theta_0-\frac{1}{r}} \frac{\chi_0(t)}{\chi_1(t)} a(\rho(t)) \left\|s^{-1/r_1} b_1(s) \left\|u^{-1-1/q_1} a_1(u) K(u,f)\right\|_{q_1,(t,s)}\right\|_{r_1,(t,\infty)}\right\|_{r,(0,\infty)}$$
$$\approx \left\|t^{-\theta_0-1/r} \frac{\chi_0(t)}{\chi_1(t)} a(\rho(t)) \|s^{-1/r_1} b_1(s)\|_{r_1,(t,\infty)} a_1(t) K(t,f)\right\|_{r,(0,\infty)}$$
$$\prec \left\|t^{-\theta_0-1/r} \chi_0(t) a(\rho(t)) K(t,f)\right\|_{r,(0,\infty)}.$$

Using now (24), we arrive at
$$I_2(0) \prec \left\|t^{-1/r} a(\rho(t)) \left\|s^{-1/r} b_0(s) \left\|u^{-\theta_0-\frac{1}{q}} a_0(u) K(u,f)\right\|_{q,(s,t)}\right\|_{r,(0,t)}\right\|_{r,(0,\infty)} = I_1(0).$$

Thus, the statement (ii) is proved. Consider $I_1(\theta)$ for $\theta > 0$. By Lemma 23 (i), we get
$$I_1(\theta) \approx \left\|t^{-\eta-1/r} \left(\frac{\chi_1(t)}{\chi_0(t)}\right)^\theta a(\rho(t)) \|s^{-1/r_0} b_0(s)\|_{r_0,(0,t)} a_0(t) K(t,f)\right\|_{r,(0,\infty)}$$
$$= \left\|t^{-\eta-1/r} (\chi_0(t))^{1-\theta} (\chi_1(t))^\theta a(\rho(t)) K(t,f)\right\|_{r,(0,\infty)} = \|f\|_{\eta,r,a^\#}.$$

Consider $I_2(\theta)$ for $0 < \theta < 1$. By Lemma 23 (ii) and (17), we get
$$I_2(\theta) \approx \left\|t^{-\eta-\frac{1}{r}} (\chi_0(t))^{1-\theta} (\chi_1(t))^{\theta-1} a(\rho(t)) \|s^{-1/r_1} b_1(s)\|_{r_1,(t,\infty)} a_1(t) K(t,f)\right\|_{r,(0,\infty)}$$
$$\prec \left\|t^{-\eta-\frac{1}{r}} (\chi_0(t))^{1-\theta} (\chi_1(t))^\theta a(\rho(t)) K(t,f)\right\|_{r,(0,\infty)} = \|f\|_{\eta,r,a^\#}.$$

Hence, the statement (i) is proved. Finally,

$$I_2(1) = \left\| t^{-1/r} a(\rho(t)) \left\| s^{-1/r_1} b_1(s) \| u^{-1-1/q_1} a_1(u) K(u,f) \|_{q_1,(t,s)} \right\|_{r_1,(t,\infty)} \right\|_{r,(0,\infty)}$$
$$= \|f\|_{\mathcal{L},\mathcal{R};1,r,a\circ\rho,r_1,b_1,q_1,a_1}.$$

Hence, the statement (iii) is also proved. □

### 4.2. Couples of the form $\left(*, \overline{A}_{1,*}^{\mathcal{R}}\right)$

In this subsection we consider the space $\overline{A}_{1,r_1,b_1,q_1,a_1}^{\mathcal{R}}$ as the second operand. Hence, in view of Lemma 12, in all theorems of this subsection, we assume that $0 < q_1, r_1 \leq \infty$, $a_1, b_1 \in SV$, and $\left\| s^{-1/r_1} b_1(s) \| u^{-1/q_1} a_1(u) \|_{q_1,(s,1)} \right\|_{r_1,(0,1)} < \infty$. Additionally, we put

$$\chi_1(t) = \left\| s^{-1/r_1} b_1(s) \| u^{-1/q_1} a_1(u) \|_{q_1,(s,t)} \right\|_{r_1,(0,t)}.$$

Note that $\chi_1 \in SV$ and due to Lemma 3,
$$a_1(t) b_1(t) \prec \left\| s^{-1/r_1} b_1(s) \right\|_{r_1,(0,t)} a_1(t)$$
$$\prec \left\| s^{-1/r_1} b_1(s) \| u^{-1/q_1} a_1(u) \|_{q_1,(s,t)} \right\|_{r_1,(0,t)} = \chi_1(t). \tag{25}$$

In proving all Holmstedt-type formulae in this subsection, by $\Phi_1$ we denote the function space corresponding to $\overline{A}_{1,r_1,b_1,q_1,a_1}^{\mathcal{R}}$:

$$\|F\|_{\Phi_1} = \left\| s^{-\frac{1}{r_1}} b_1(s) \| u^{-1-1/q_1} a_1(u) F(u) \|_{q_1,(s,\infty)} \right\|_{r_1,(0,\infty)}.$$

Following [2], we consider the next functions $J(t,f)$, $g_1(t)$, and $h_1(t)$.

$$J(t,f) := \left\| \chi_{(t,\infty)}(*) K(*,f) \right\|_{\Phi_1} = \left\| s^{-\frac{1}{r_1}} b_1(s) \left\| u^{-1-\frac{1}{q_1}} a_1(u) K(u,f) \right\|_{q_1,(\max(s,t),\infty)} \right\|_{r_1,(0,\infty)}$$
$$\approx \left\| s^{-\frac{1}{r_1}} b_1(s) \right\|_{r_1,(0,t)} \left\| u^{-1-\frac{1}{q_1}} a_1(u) K(u,f) \right\|_{q_1,(t,\infty)}$$
$$+ \left\| s^{-\frac{1}{r_1}} b_1(s) \left\| u^{-1-\frac{1}{q_1}} a_1(u) K(u,f) \right\|_{q_1,(s,\infty)} \right\|_{r_1,(t,\infty)},$$

$$g_1(t) := t \left\| \chi_{(t,\infty)}(*) \right\|_{\Phi_1} = t \left\| s^{-\frac{1}{r_1}} b_1(s) \left\| u^{-1-\frac{1}{q_1}} a_1(u) \right\|_{q_1,(\max(s,t),\infty)} \right\|_{r_1,(0,\infty)}$$
$$\approx t \left( \left\| u^{-1-\frac{1}{q_1}} a_1(u) \right\|_{q_1,(t,\infty)} \left\| s^{-\frac{1}{r_1}} b_1(s) \right\|_{r_1,(0,t)} \right.$$
$$\left. + \left\| s^{-\frac{1}{r_1}} b_1(s) \left\| u^{-1-\frac{1}{q_1}} a_1(u) \right\|_{q_1,(s,\infty)} \right\|_{r_1,(t,\infty)} \right)$$
$$\approx t \left( t^{-1} a_1(t) \left\| s^{-\frac{1}{r_1}} b_1(s) \right\|_{r_1,(0,t)} + \left\| s^{-1-\frac{1}{r_1}} b_1(s) a_1(s) \right\|_{r_1,(t,\infty)} \right)$$
$$\approx a_1(t) \left( \left\| s^{-\frac{1}{r_1}} b_1(s) \right\|_{r_1,(0,t)} + b_1(t) \right) \approx a_1(t) \left\| s^{-\frac{1}{r_1}} b_1(s) \right\|_{r_1,(0,t)},$$

and
$$h_1(t) := \| * \chi_{(0,t)}(*) \|_{\Phi_1} = \left\| s^{-\frac{1}{r_1}} b_1(s) \left\| u^{-\frac{1}{q_1}} a_1(u) \right\|_{q_1,(s,t)} \right\|_{r_1,(0,t)} = \chi_1(t).$$

Because of (25), we conclude that
$$g_1(t) \prec h_1(t) \approx g_1(t) + h_1(t) \approx \chi_1(t). \tag{26}$$

**Theorem 39.** Put $\rho(t) = \frac{t}{\chi_1(t)}$. For all $f \in A_0 + \bar{A}^{\mathcal{R}}_{1,r_1,b_1,q_1,a_1}$ and $t > 0$,
$$K(\rho(t), f; A_0, \bar{A}^{\mathcal{R}}_{1,r_1,b_1,q_1,a_1}) \approx K(t,f)$$
$$+ \rho(t) \left( \left\| s^{-\frac{1}{r_1}} b_1(s) \right\|_{r_1,(0,t)} \left\| u^{-1-\frac{1}{q_1}} a_1(u) K(u,f) \right\|_{q_1,(t,\infty)} + \right.$$
$$\left. \left\| s^{-\frac{1}{r_1}} b_1(s) \left\| u^{-1-\frac{1}{q_1}} a_1(u) K(u,f) \right\|_{q_1,(s,\infty)} \right\|_{r_1,(t,\infty)} \right).$$

*Proof.* By (26), we arrive at $\rho(t) \approx \frac{t}{g_1(t)+h_1(t)}$. Using of [2, Theorem 3], completes the proof. □

**Theorem 40.** Let $0 < r \leq \infty$ and $a \in SV$. Put $\rho(t) = \frac{t}{\chi_1(t)}$. For given $0 \leq \theta \leq 1$, denote $a^\# = \chi_1^\theta a \circ \rho$.

(i) If $0 < \theta < 1$, then
$$\left(A_0, \bar{A}^{\mathcal{R}}_{1,r_1,b_1,q_1,a_1}\right)_{\theta,r,a} = \bar{A}_{\theta,r,a^\#}.$$

(ii) If $\left\| s^{-1/r} a(s) \right\|_{r,(1,\infty)} < \infty$, then
$$\left(A_0, \bar{A}^{\mathcal{R}}_{1,r_1,b_1,q_1,a_1}\right)_{0,r,a} = \bar{A}_{0,r,a\circ\rho},$$

(iii) If $\left\| s^{-1/r} a(s) \right\|_{r,(0,1)} < \infty$, then
$$\left(A_0, \bar{A}^{\mathcal{R}}_{1,r_1,b_1,q_1,a_1}\right)_{1,r,a} = \bar{A}_{1,r,a^\#} \cap \bar{A}^{\mathcal{R}}_{1,r,B_1,q_1,a_1} \cap \bar{A}^{\mathcal{R},\mathcal{R}}_{1,r,a\circ\rho,r_1,b_1,q_1,a_1},$$
where $B_1(t) = \left\| s^{-1/r_1} b_1(s) \right\|_{r_1,(0,t)} a(\rho(t))$.

*Proof.* First, note that due to Lemma 5 (ii) and Lemma 12 the space $\bar{A}^{\mathcal{R}}_{1,r,B_1,q_1,a_1}$ is an intermediate space for the couple $(A_0, A_1)$. Denote $\bar{K}(t,f) = K(t,f; A_0, \bar{A}^{\mathcal{R}}_{1,r_1,b_1,q_1,a_1})$ and $Z = \left(A_0, \bar{A}^{\mathcal{R}}_{1,r_1,b_1,q_1,a_1}\right)_{\theta,r,a}$ $(0 \leq \theta \leq 1)$. By a change of variables and Theorem 39 we can write
$$\|f\|_Z \approx \left\| \rho(t)^{-\theta} t^{-1/r} a(\rho(t)) \bar{K}(\rho(t),f) \right\|_{r,(0,\infty)} \approx I_1(\theta) + I_2(\theta) + I_3(\theta),$$
where
$$I_1(\theta) := \left\| t^{-\theta-1/r} \chi_1(t)^\theta a(\rho(t)) K(t,f) \right\|_{r,(0,\infty)} = \|f\|_{\theta,r;a^\#},$$
$$I_2(\theta) := \left\| \rho(t)^{1-\theta} t^{-1/r} a(\rho(t)) \left\| s^{-1/r_1} b_1(s) \right\|_{r_1,(0,t)} \left\| u^{-1-1/q_1} a_1(u) K(u,f) \right\|_{q_1,(t,\infty)} \right\|_{r,(0,\infty)}$$
$$= \left\| t^{1-\theta-1/r} a(\rho(t)) \chi_1(t)^{\theta-1} \left\| s^{-1/r_1} b_1(s) \right\|_{r_1,(0,t)} \left\| u^{-1-1/q_1} a_1(u) K(u,f) \right\|_{q_1,(t,\infty)} \right\|_{r,(0,\infty)},$$
and
$$I_3(\theta) := \left\| \rho(t)^{1-\theta} t^{-\frac{1}{r}} a(\rho(t)) \left\| s^{-1/r_1} b_1(s) \left\| u^{-1-1/q_1} a_1(u) K(u,f) \right\|_{q_1,(s,\infty)} \right\|_{r_1,(t,\infty)} \right\|_{r,(0,\infty)}$$
$$= \left\| t^{1-\theta-1/r} a(\rho(t)) \chi_1(t)^{\theta-1} \left\| s^{-1/r_1} b_1(s) \left\| u^{-1-1/q_1} a_1(u) K(u,f) \right\|_{q_1,(s,\infty)} \right\|_{r_1,(t,\infty)} \right\|_{r,(0,\infty)}.$$

Consider $I_2$. Let $0 \leq \theta < 1$. By Lemma 21 and (25), we have
$$I_2(\theta) \approx \left\| t^{-\theta-1/r} a(\rho(t)) \chi_1(t)^{\theta-1} \left\| s^{-1/r_1} b_1(s) \right\|_{r_1,(0,t)} a_1(t) K(t,f) \right\|_{r,(0,\infty)}$$
$$\prec \left\| t^{-\theta-1/r} a(\rho(t)) \chi_1(t)^\theta K(t,f) \right\|_{r,(0,\infty)} = \|f\|_{\theta,r;a^\#}.$$

Additionally,
$$I_2(1) = \|f\|_{\mathcal{R};1,r,B_1,q_1,a_1},$$
where $B_1(t) = \left\| s^{-1/r_1} b_1(s) \right\|_{r_1,(0,t)} a(\rho(t))$. Consider $I_3$. If $0 \leq \theta < 1$, using Lemma 22 and (25), we get
$$I_3(\theta) \approx \left\| t^{-\theta-1/r} a(\rho(t)) \chi_1(t)^{\theta-1} a_1(t) b_1(t) K(t,f) \right\|_{r,(0,\infty)}$$

$$\prec \left\|t^{-\theta-1/r}a(\rho(t))\chi_1(t)^\theta K(t,f)\right\|_{r,(0,\infty)} = \|f\|_{\theta,r;a^{\#}}.$$

If $\theta = 1$, we have
$$I_3(1) = \|f\|_{\mathcal{R},\mathcal{R};1,r,a\circ\rho,r_1,b_1,q_1,a_1}.$$

This completes the proof. $\square$

**Theorem** 41. Let $0 < r_0 \leq \infty$, $b_0 \in SV$, and $\left\|s^{-1/r_0}b_0(s)\right\|_{r_0,(1,\infty)} < \infty$. Put $\chi_0(t) = \left\|s^{-1/r_0}b_0(s)\right\|_{r_0,(t,\infty)}$ and $\rho(t) = t\frac{\chi_0(t)}{\chi_1(t)}$. Then, for all $f \in A_{0,r_0,b_0} + \bar{A}^{\mathcal{R}}_{1,r_1,b_1,q_1,a_1}$ and $t > 0$,
$$K\left(\rho(t),f;A_{0,r_0,b_0},\bar{A}^{\mathcal{R}}_{1,r_1,b_1,q_1,a_1}\right) \approx$$
$$\left\|u^{-1/r_0}b_0(u)K(u,f)\right\|_{r_0,(0,t)} + \rho(t)\left(\left\|s^{-1/r_1}b_1(s)\right\|_{r_1,(0,t)}\left\|u^{-1-1/q_1}a_1(u)K(u,f)\right\|_{q_1,(t,\infty)} + \right.$$
$$\left. \left\|s^{-1/r_1}b_1(s)\left\|u^{-1-1/q_1}a_1(u)K(u,f)\right\|_{q_1,(s,\infty)}\right\|_{r_1,(t,\infty)}\right) + \chi_0(t)K(t,f).$$

*Proof.* Let $\Phi_0$ be the function space corresponding to $\bar{A}_{0,r_0,b_0}$:
$$\|F\|_{\Phi_0} = \left\|u^{-1/r_0}b_0(u)F(u)\right\|_{r_0,(0,\infty)}.$$

From the proof of Theorem 26, we know that
$$I(t,f) := \left\|\chi_{(0,t)}(*)K(*,f)\right\|_{\Phi_0} = \left\|u^{-1/r_0}b_0(u)K(u,f)\right\|_{r_0,(0,t)}$$

and
$$h_0(t) \prec g_0(t) \approx g_0(t) + h_0(t) \approx t\chi_0(t).$$

Now we estimate the functions $g(t)$ and $h(t)$.
$$\frac{1}{g(t)} := \left\|\chi_{(t,\infty)}(u)\frac{u}{g_0(u)+h_0(u)}\right\|_{\Phi_1}$$
$$\approx \left\|s^{-1/r_1}b_1(s)\left\|u^{-1-1/q_1}\frac{a_1(u)}{\chi_0(u)}\right\|_{q_1,(\max(t,s),\infty)}\right\|_{r_1,(0,\infty)}$$
$$\approx \left\|u^{-1-1/q_1}\frac{a_1(u)}{\chi_0(u)}\right\|_{q_1,(t,\infty)}\left\|s^{-1/r_1}b_1(s)\right\|_{r_1,(0,t)}$$
$$+ \left\|s^{-1/r_1}b_1(s)\left\|u^{-1-1/q_1}\frac{a_1(u)}{\chi_0(u)}\right\|_{q_1,(s,\infty)}\right\|_{r_1,(t,\infty)}$$
$$\approx t^{-1}\frac{a_1(t)}{\chi_0(t)}\left\|s^{-1/r_1}b_1(s)\right\|_{r_1,(0,t)} + \left\|s^{-1-1/r_1}b_1(s)\frac{a_1(s)}{\chi_0(s)}\right\|_{r_1,(t,\infty)}$$
$$\approx t^{-1}\frac{a_1(t)}{\chi_0(t)}\left(\left\|s^{-1/r_1}b_1(s)\right\|_{r_1,(0,t)} + b_1(t)\right) \approx t^{-1}\frac{a_1(t)\left\|s^{-1/r_1}b_1(s)\right\|_{r_1,(0,t)}}{\chi_0(t)} \prec t^{-1}\frac{\chi_1(t)}{\chi_0(t)} = \frac{1}{\rho(t)}.$$

Here we have used estimate (25). Repeating the arguments of the proof of Theorem 26, we get
$$h(t) := \left\|\chi_{(0,t)}(u)\frac{u}{g_1(u)+h_1(u)}\right\|_{\Phi_0} \prec \rho(t).$$

Hence,
$$h(t) \prec \rho(t) \prec g(t).$$

Note that
$$\frac{g_0(t)}{g_1(t)+h_1(t)} \approx \frac{g_0(t)+h_0(t)}{h_1(t)} \approx \frac{t\chi_0(t)}{\chi_1(t)} = \rho(t)$$

and
$$g_0(t) + \rho(t)h_1(t) \approx t\chi_0(t) + t\frac{\chi_0(t)}{\chi_1(t)}\chi_1(t) \approx t\chi_0(t).$$

Theorem 4 (Case 1) of [2] completes the proof. □

**Theorem** 42. Let $0 < r, r_0 \leq \infty$, $a, b_0 \in SV$, and $\left\|s^{-1/r_0}b_0(s)\right\|_{r_0,(1,\infty)} < \infty$. Put $\chi_0(t) = \left\|s^{-1/r_0}b_0(s)\right\|_{r_0,(t,\infty)}$ and $\rho(t) = t\frac{\chi_0(t)}{\chi_1(t)}$. For a given $0 \leq \theta \leq 1$, denote $a^\# = \chi_0^{1-\theta}\chi_1^\theta a \circ \rho$.

(i) If $0 < \theta < 1$, then
$$\left(A_{0,r_0,b_0}, \bar{A}^{\mathcal{R}}_{1,r_1,b_1,q_1,a_1}\right)_{\theta,r,a} = \bar{A}_{\theta,r,a^\#}.$$

(ii) If $\left\|s^{-1/r}a(s)\right\|_{r,(1,\infty)} < \infty$, then
$$\left(A_{0,r_0,b_0}, \bar{A}^{\mathcal{R}}_{1,r_1,b_1,q_1,a_1}\right)_{0,r,a} = \bar{A}_{0,r,a^\#} \cap \bar{A}^{\mathcal{L}}_{0,r,a \circ \rho,r_0,b_0}.$$

(iii) If $\left\|s^{-1/r}a(s)\right\|_{r,(0,1)} < \infty$, then
$$\left(A_{0,r_0,b_0}, \bar{A}^{\mathcal{R}}_{1,r_1,b_1,q_1,a_1}\right)_{1,r,a} = \bar{A}_{1,r,a^\#} \cap \bar{A}^{\mathcal{R}}_{1,r,B_1,q_1,a_1} \cap \bar{A}^{\mathcal{R},\mathcal{R}}_{1,r,a \circ \rho,r_1,b_1,q_1,a_1},$$
where $B_1(t) = \left\|s^{-1/r_1}b_1(s)\right\|_{r_1,(0,t)}a(\rho(t))$.

*Proof.* First, note that due to Lemma 5 (ii) and Lemma 12 the space $\bar{A}^{\mathcal{R}}_{1,r,B_1,q_1,a_1}$ is an intermediate space for the couple $(A_0, A_1)$. Denote $\bar{K}(t,f) = K\left(t,f; A_{0,r_0,b_0}, \bar{A}^{\mathcal{R}}_{1,r_1,b_1,q_1,a_1}\right)$ and $Z = \left(A_{0,r_0,b_0}, \bar{A}^{\mathcal{R}}_{1,r_1,b_1,q_1,a_1}\right)_{\theta,r,a}$ ($0 \leq \theta \leq 1$). By a change of variables and Theorem 41 we can write
$$\|f\|_Z \approx \left\|\rho(t)^{-\theta}t^{-1/r}a(\rho(t))\bar{K}(\rho(t),f)\right\|_{r,(0,\infty)} \approx I_1(\theta) + I_2(\theta) + I_3(\theta) + I_4(\theta),$$
where
$$I_1(\theta) := \left\|\rho(t)^{-\theta}t^{-1/r}a(\rho(t))\left\|u^{-1/r_0}b_0(u)K(u,f)\right\|_{r_0,(0,t)}\right\|_{r,(0,\infty)}$$
$$= \left\|t^{-\theta-1/r}a(\rho(t))\left(\frac{\chi_1(t)}{\chi_0(t)}\right)^\theta\left\|u^{-1/r_0}b_0(u)K(u,f)\right\|_{r_0,(0,t)}\right\|_{r,(0,\infty)},$$
$$I_2(\theta) := \left\|\rho(t)^{1-\theta}t^{-1/r}a(\rho(t))\left\|s^{-1/r_1}b_1(s)\right\|_{r_1,(0,t)}\left\|u^{-1-1/q_1}a_1(u)K(u,f)\right\|_{q_1,(t,\infty)}\right\|_{r,(0,\infty)}$$
$$= \left\|t^{1-\theta-1/r}a(\rho(t))\left(\frac{\chi_1(t)}{\chi_0(t)}\right)^{\theta-1}\left\|s^{-1/r_1}b_1(s)\right\|_{r_1,(0,t)}\left\|u^{-1-1/q_1}a_1(u)K(u,f)\right\|_{q_1,(t,\infty)}\right\|_{r,(0,\infty)},$$
$$I_3(\theta) := \left\|\rho(t)^{1-\theta}t^{-\frac{1}{r}}a(\rho(t))\left\|s^{-1/r_1}b_1(s)\left\|u^{-1-1/q_1}a_1(u)K(u,f)\right\|_{q_1,(s,\infty)}\right\|_{r_1,(t,\infty)}\right\|_{r,(0,\infty)}$$
$$= \left\|t^{1-\theta-\frac{1}{r}}a(\rho(t))\left(\frac{\chi_1(t)}{\chi_0(t)}\right)^{\theta-1}\left\|s^{-1/r_1}b_1(s)\left\|u^{-1-1/q_1}a_1(u)K(u,f)\right\|_{q_1,(s,\infty)}\right\|_{r_1,(t,\infty)}\right\|_{r,(0,\infty)}.$$

and
$$I_4(\theta) := \left\|\rho(t)^{-\theta}t^{-1/r}a(\rho(t))\chi_0(t)K(t,f)\right\|_{r,(0,\infty)}$$
$$= \left\|t^{-\theta-1/r}\chi_0(t)^{1-\theta}\chi_1(t)^\theta a(\rho(t))K(t,f)\right\|_{r,(0,\infty)} = \|f\|_{\theta,r;a^\#}.$$

Consider $I_1$. We have
$$I_1(0) = \left\|t^{-1/r}a(\rho(t))\left\|u^{-1/r_0}b_0(u)K(u,f)\right\|_{r_0,(0,t)}\right\|_{r,(0,\infty)} = \|f\|_{\mathcal{L};0,r,a \circ \rho,r_0,b_0}.$$

Let $0 < \theta \leq 1$. Because $b_0(t) \prec \chi_0(t)$, using the first estimate from Lemma 21, we conclude that
$$I_1(\theta) \approx \left\|t^{-\theta-1/r}a(\rho(t))\chi_0(t)^{-\theta}\chi_1(t)^\theta b_0(t)K(t,f)\right\|_{r,(0,\infty)}$$
$$\prec \left\|t^{-\theta-1/r}a(\rho(t))\chi_0(t)^{1-\theta}\chi_1(t)^\theta K(t,f)\right\|_{r,(0,\infty)} = \|f\|_{\theta,r;a^\#}.$$

Consider $I_2$. Let $0 \leq \theta < 1$. Lemma 21 and (25) imply
$$I_2(\theta) \approx \left\|t^{-\theta-1/r}a(\rho(t))\chi_0(t)^{1-\theta}\chi_1(t)^{\theta-1}\left\|s^{-1/r_1}b_1(s)\right\|_{r_1,(0,t)}a_1(t)K(t,f)\right\|_{r,(0,\infty)}$$

$$\prec \left\|t^{-\theta-1/r}a(\rho(t))\chi_0(t)^{1-\theta}\chi_1(t)^\theta K(t,f)\right\|_{r,(0,\infty)} = \|f\|_{\theta,r;a^\#}.$$
Additionally,
$$I_2(1) = \|f\|_{\mathcal{R};1,r,B_1,q_1,a_1},$$
where
$$B_1(t) = \left\|s^{-1/r_1}b_1(s)\right\|_{r_1,(0,t)}a(\rho(t)).$$
Consider $I_3$. Let $0 \leq \theta < 1$. Using Lemma 22 and (25), we get
$$I_3(\theta) \approx \left\|t^{-\theta-1/r}a(\rho(t))\left(\frac{\chi_1(t)}{\chi_0(t)}\right)^{\theta-1}a_1(t)b_1(t)K(t,f)\right\|_{r,(0,\infty)}$$
$$\prec \left\|t^{-\theta-1/r}a(\rho(t))\chi_0(t)^{1-\theta}\chi_1(t)^\theta K(t,f)\right\|_{r,(0,\infty)} = \|f\|_{\theta,r;a^\#}.$$
Finally,
$$I_3(1) = \|f\|_{\mathcal{R},\mathcal{R};1,r,a^\circ\rho,r_1,b_1,q_1,a_1}.$$
This completes the proof. $\square$

**Theorem 43.** Let $0 < \theta_0 < 1$, $0 < r_0 \leq \infty$, and $b_0 \in SV$. Put $\rho(t) = t^{1-\theta_0}\frac{b_0(t)}{\chi_1(t)}$. Then, for all $f \in \bar{A}_{\theta_0,r_0,b_0} + \bar{A}^{\mathcal{R}}_{1,r,b,q,a}$ and $t > 0$,
$$K\left(\rho(t),f;\bar{A}_{\theta_0,r_0,b_0},\bar{A}^{\mathcal{R}}_{1,r_1,b_1,q_1,a_1}\right) \approx \left\|u^{-\theta_0-1/r_0}b_0(u)K(u,f)\right\|_{r_0,(0,t)}$$
$$+\rho(t)\left(\left\|s^{-1/r_1}b_1(s)\right\|_{r_1,(0,t)}\left\|u^{-1-1/q_1}a_1(u)K(u,f)\right\|_{q_1,(t,\infty)} + \right.$$
$$\left.\left\|s^{-1/r_1}b_1(s)\left\|u^{-1-1/q_1}a_1(u)K(u,f)\right\|_{q_1,(s,\infty)}\right\|_{r_1,(t,\infty)}\right).$$

*Proof.* Let $\Phi_0$ be the function space corresponding to $\bar{A}_{\theta_0,r_0,b_0}$:
$$\|F\|_{\Phi_0} = \left\|u^{-\theta_0-1/r_0}b_0(u)F(u)\right\|_{r_0,(0,\infty)}.$$
From the proof of [8, Theorem 18], we know that
$$I(t,f) := \left\|\chi_{(0,t)}(*)K(*,f)\right\|_{\Phi_0} = \left\|u^{-\theta_0-1/r_0}b_0(u)K(u,f)\right\|_{r_0,(0,t)}$$
and
$$g_0(t) + h_0(t) \approx h_0(t) \approx g_0(t) \approx t^{1-\theta_0}b_0(t).$$
Now we estimate the functions $g(t)$ and $h(t)$.
$$\frac{1}{g(t)} := \left\|\chi_{(t,\infty)}(u)\frac{u}{g_0(u)+h_0(u)}\right\|_{\Phi_1}$$
$$\approx \left\|s^{-1/r_1}b_1(s)\left\|u^{\theta_0-1-1/q_1}\frac{a_1(u)}{b_0(u)}\right\|_{q_1,(\max(t,s),\infty)}\right\|_{r_1,(0,\infty)}$$
$$\approx \left\|u^{\theta_0-1-1/q_1}\frac{a_1(u)}{b_0(u)}\right\|_{q_1,(t,\infty)}\left\|s^{-1/r_1}b_1(s)\right\|_{r_1,(0,t)}$$
$$+ \left\|s^{-1/r_1}b_1(s)\left\|u^{\theta_0-1-1/q_1}\frac{a_1(u)}{b_0(u)}\right\|_{q_1,(s,\infty)}\right\|_{r_1,(t,\infty)}$$
$$\approx t^{\theta_0-1}\frac{a_1(t)}{b_0(t)}\left\|s^{-1/r_1}b_1(s)\right\|_{r_1,(0,t)} + \left\|s^{\theta_0-1-1/r_1}\frac{b_1(s)a_1(s)}{b_0(s)}\right\|_{r_1,(t,\infty)}$$
$$\approx t^{\theta_0-1}\frac{a_1(t)}{b_0(t)}\left\|s^{-1/r_1}b_1(s)\right\|_{r_1,(0,t)} + t^{\theta_0-1}\frac{b_1(t)a_1(t)}{b_0(t)}$$
$$= t^{\theta_0-1}\frac{a_1(t)}{b_0(t)}\left(\left\|s^{-1/r_1}b_1(s)\right\|_{r_1,(0,t)} + b_1(t)\right)$$

$$\approx t^{\theta_0-1}\frac{a_1(t)\|s^{-1/r_1}b_1(s)\|_{r_1,(0,t)}}{b_0(t)} \prec t^{\theta_0-1}\frac{\chi_1(t)}{b_0(t)} = \frac{1}{\rho(t)}.$$

Here we have used estimate (25). Repeating the arguments of the proof of Theorem 28, we get

$$h(t) := \left\|\chi_{(0,t)}(u)\frac{u}{g_1(u)+h_1(u)}\right\|_{\Phi_0} \approx \rho(t).$$

Thus,
$$h(t) \approx \rho(t) \prec g(t).$$

Theorem 4 (Case 1) of [2] completes the proof. □

**Theorem** 44. Let $0 < \theta_0 < 1$, $0 < r, r_0 \le \infty$, and $a, b_0 \in SV$. Put $\rho(t) = t^{1-\theta_0}\frac{b_0(t)}{\chi_1(t)}$. For a given $0 \le \theta \le 1$, denote $a^\# = b_0^{1-\theta}\chi_1^\theta a\circ\rho$ and $\eta = (1-\theta)\theta_0 + \theta$.

(i) If $0 < \theta < 1$, then
$$(\bar{A}_{\theta_0,r_0,b_0}, \bar{A}^{\mathcal{R}}_{1,r_1,b_1,q_1,a_1})_{\theta,r,a} = \bar{A}_{\eta,r,a^\#}.$$

(ii) If $\|s^{-1/r}a(s)\|_{r,(1,\infty)} < \infty$, then
$$(\bar{A}_{\theta_0,r_0,b_0}, \bar{A}^{\mathcal{R}}_{1,r_1,b_1,q_1,a_1})_{0,r,a} = \bar{A}^{\mathcal{L}}_{\theta_0,r,a\circ\rho,r_0,b_0}.$$

(iii) If $\|s^{-1/r}a(s)\|_{r,(0,1)} < \infty$, then
$$(\bar{A}_{\theta_0,r_0,b_0}, \bar{A}^{\mathcal{R}}_{1,r_1,b_1,q_1,a_1})_{1,r,a} = \bar{A}_{1,r,a^\#} \cap \bar{A}^{\mathcal{R}}_{1,r,B_1,q_1,a_1} \cap \bar{A}^{\mathcal{R},\mathcal{R}}_{1,r,a\circ\rho,r_1,b_1,q_1,a_1},$$

where $B_1(t) = \|s^{-1/r_1}b_1(s)\|_{r_1,(0,t)}a(\rho(t))$.

*Proof.* First, note that due to Lemma 5 (ii) and Lemma 12 the space $\bar{A}^{\mathcal{R}}_{1,r,B_1,q_1,a_1}$ is an intermediate space for the couple $(A_0, A_1)$. Denote $\bar{K}(t, f) = K(t, f; \bar{A}_{\theta_0,r_0,b_0}, \bar{A}^{\mathcal{R}}_{1,r_1,b_1,q_1,a_1})$ and $Z = (\bar{A}_{\theta_0,r_0,b_0}, \bar{A}^{\mathcal{R}}_{1,r_1,b_1,q_1,a_1})_{\theta,r,a}$ ($0 \le \theta \le 1$). By a change of variables and **Theorem** 43, we can write

$$\|f\|_Z \approx \left\|\rho(t)^{-\theta}t^{-1/r}a(\rho(t))\bar{K}(\rho(t),f)\right\|_{r,(0,\infty)} \approx I_1(\theta) + I_2(\theta) + I_3(\theta),$$

where

$$I_1(\theta) := \left\|\rho(t)^{-\theta}t^{-1/r}a(\rho(t))\left\|u^{-\theta_0-1/r_0}b_0(u)K(u,f)\right\|_{r_0,(0,t)}\right\|_{r,(0,\infty)}$$
$$= \left\|t^{-\theta(1-\theta_0)-1/r}a(\rho(t))b_0(t)^{-\theta}\chi_1(t)^\theta\left\|u^{-\theta_0-1/r_0}b_0(u)K(u,f)\right\|_{r_0,(0,t)}\right\|_{r,(0,\infty)},$$

$$I_2(\theta) := \left\|\rho(t)^{1-\theta}t^{-1/r}a(\rho(t))\left\|s^{-1/r_1}b_1(s)\right\|_{r_1,(0,t)}\left\|u^{-1-1/q_1}a_1(u)K(u,f)\right\|_{q_1,(t,\infty)}\right\|_{r,(0,\infty)}$$
$$=$$
$$\left\|t^{(1-\theta)(1-\theta_0)-\frac{1}{r}}a(\rho(t))\left(\frac{\chi_1(t)}{b_0(t)}\right)^{\theta-1}\left\|s^{-\frac{1}{r_1}}b_1(s)\right\|_{r_1,(0,t)}\left\|u^{-1-\frac{1}{q_1}}a_1(u)K(u,f)\right\|_{q_1,(t,\infty)}\right\|_{r,(0,\infty)},$$

and

$$I_3(\theta) := \left\|\rho(t)^{1-\theta}t^{-\frac{1}{r}}a(\rho(t))\left\|s^{-\frac{1}{r_1}}b_1(s)\left\|u^{-1-\frac{1}{q_1}}a_1(u)K(u,f)\right\|_{q_1,(s,\infty)}\right\|_{r_1,(t,\infty)}\right\|_{r,(0,\infty)}$$
$$= \left\|t^{(1-\theta)(1-\theta_0)-\frac{1}{r}}a(\rho(t))\left(\frac{\chi_1(t)}{b_0(t)}\right)^{\theta-1}\left\|s^{-\frac{1}{r_1}}b_1(s)\left\|u^{-1-\frac{1}{q_1}}a_1(u)K(u,f)\right\|_{q_1,(s,\infty)}\right\|_{r_1,(t,\infty)}\right\|_{r,(0,\infty)}$$

Consider $I_1$. We have
$$I_1(0) = \left\|t^{-1/r}a(\rho(t))\left\|u^{-\theta_0-1/r_0}b_0(u)K(u,f)\right\|_{r_0,(0,t)}\right\|_{r,(0,\infty)} = \|f\|_{\mathcal{L};\theta_0,r,a\circ\rho,r_0,b_0}.$$

Let $0 < \theta \le 1$. Using the first estimate from Lemma 21, we conclude that

$$I_1(\theta) \approx \left\| t^{-\eta-1/r} a(\rho(t)) b_0(t)^{-\theta} \chi_1(t)^\theta b_0(t) K(t,f) \right\|_{r,(0,\infty)}$$
$$= \left\| t^{-\eta-1/r} a(\rho(t)) b_0(t)^{1-\theta} \chi_1(t)^\theta K(t,f) \right\|_{r,(0,\infty)} = \|f\|_{\eta,r;a^\#}.$$

Consider $I_2$. Let $0 \le \theta < 1$. Lemma 21 (ii) and (25) imply
$$I_2(\theta) \approx \left\| t^{-\eta-1/r} a(\rho(t)) b_0(t)^{1-\theta} \chi_1(t)^{\theta-1} \left\| s^{-1/r_1} b_1(s) \right\|_{r_1,(0,t)} a_1(t) K(t,f) \right\|_{r,(0,\infty)}$$
$$\prec \left\| t^{-\eta-1/r} a(\rho(t)) b_0(t)^{1-\theta} \chi_1(t)^\theta K(t,f) \right\|_{r,(0,\infty)} = \|f\|_{\eta,r;a^\#}.$$

If $\theta = 1$, we have
$$I_2(1) = \|f\|_{\mathcal{R};1,r,B_1,q_1,a_1},$$
where
$$B_1(t) = \left\| s^{-1/r_1} b_1(s) \right\|_{r_1,(0,t)} a(\rho(t)).$$

Consider $I_3$. Let $0 \le \theta < 1$. Using Lemma 22 and (25), we get
$$I_3(\theta) \approx \left\| t^{-\eta-1/r} a(\rho(t)) b_0(t)^{1-\theta} \chi_1(t)^{\theta-1} b_1(t) a_1(t) K(t,f) \right\|_{r,(0,\infty)}$$
$$\prec \left\| t^{-\eta-1/r} a(\rho(t)) b_0(t)^{1-\theta} \chi_1(t)^\theta K(t,f) \right\|_{r,(0,\infty)} = \|f\|_{\eta,r;a^\#}.$$

Finally,
$$I_3(1) = \|f\|_{\mathcal{R},\mathcal{R};1,r,a^\circ\rho,r_1,b_1,q_1,a_1}.$$
Thus, (i) and (iii) are proved. To prove (ii), it is enough to show that
$$\|f\|_{\theta_0,r,b_0 a^\circ \rho} \prec \|f\|_{\mathcal{L};\theta_0,r,a^\circ\rho,r_0,b_0}.$$
But this follows from Corollary 17 (i). □

**Theorem** 45. Let $0 < r_0, q_0 \le \infty$, $a_0, b_0 \in SV$, and
$$\left\| s^{-1/r_0} b_0(s) \left\| u^{-1/q_0} a_0(u) \right\|_{q_0,(1,s)} \right\|_{r_0,(1,\infty)} < \infty.$$
Put
$$\chi_0(t) = \left\| s^{-1/r_0} b_0(s) \left\| u^{-1/q_0} a_0(u) \right\|_{q_0,(t,s)} \right\|_{r_0,(t,\infty)}$$
and $\rho(t) = t \frac{\chi_0(t)}{\chi_1(t)}$. Then, for all $f \in \bar{A}^{\mathcal{L}}_{0,r_0,b_0,q_0,a_0} + \bar{A}^{\mathcal{R}}_{1,r_1,b_1,q_1,a_1}$ and $t > 0$,
$$K\left(\rho(t), f; \bar{A}^{\mathcal{L}}_{0,r_0,b_0,q_0,a_0}, \bar{A}^{\mathcal{R}}_{1,r_1,b_1,q_1,a_1}\right)$$
$$\approx \left\| s^{-1/r_0} b_0(s) \left\| u^{-1/q_0} a_0(u) K(u,f) \right\|_{q_0,(0,s)} \right\|_{r_0,(0,t)}$$
$$+ \left\| s^{-1/r_0} b_0(s) \right\|_{r_0,(t,\infty)} \left\| u^{-1/q_0} a_0(u) K(u,f) \right\|_{q_0,(0,t)}$$
$$+ \rho(t) \left( \left\| s^{-1/r_1} b_1(s) \right\|_{r_1,(0,t)} \left\| u^{-1-1/q_1} a_1(u) K(u,f) \right\|_{q_1,(t,\infty)} + \right.$$
$$\left. \left\| s^{-1/r_1} b_1(s) \left\| u^{-1-1/q_1} a_1(u) K(u,f) \right\|_{q_1,(s,\infty)} \right\|_{r_1,(t,\infty)} \right) + \chi_0(t) K(t,f).$$

*Proof.* Let $\Phi_0$ be the function space corresponding to $\bar{A}^{\mathcal{L}}_{0,r_0,b_0,q_0,a_0}$:
$$\|F\|_{\Phi_0} = \left\| s^{-1/r_0} b_0(s) \left\| u^{-1/q_0} a_0(u) F(u) \right\|_{q_0,(0,s)} \right\|_{r_0,(0,\infty)}.$$
From the proof of Theorem 30, we know that
$$I(t,f) := \left\| \chi_{(0,t)}(*) K(*,f) \right\|_{\Phi_0} \approx \left\| s^{-1/r_0} b_0(s) \left\| u^{-1/q_0} a_0(u) K(u,f) \right\|_{q_0,(0,s)} \right\|_{r_0,(0,t)}$$
$$+ \left\| s^{-1/r_0} b_0(s) \right\|_{r_0,(t,\infty)} \left\| u^{-1/q_0} a_0(u) K(u,f) \right\|_{q_0,(0,t)}$$
and
$$h_0(t) \prec g_0(t) \approx g_0(t) + h_0(t) \approx t \chi_0(t).$$
Now we estimate the functions $g(t)$ and $h(t)$.

$$\frac{1}{g(t)} := \left\| \chi_{(t,\infty)}(u) \frac{u}{g_0(u) + h_0(u)} \right\|_{\Phi_1}$$

$$\approx \left\| s^{-1/r_1} b_1(s) \left\| u^{-1-1/q_1} \frac{a_1(u)}{\chi_0(u)} \right\|_{q_1,(\max(t,s),\infty)} \right\|_{r_1,(0,\infty)}$$

$$\approx \left\| u^{-1-1/q_1} \frac{a_1(u)}{\chi_0(u)} \right\|_{q_1,(t,\infty)} \left\| s^{-1/r_1} b_1(s) \right\|_{r_1,(0,t)}$$

$$+ \left\| s^{-1/r_1} b_1(s) \left\| u^{-1-1/q_1} \frac{a_1(u)}{\chi_0(u)} \right\|_{q_1,(s,\infty)} \right\|_{r_1,(t,\infty)}$$

$$\approx t^{-1} \frac{a_1(t)}{\chi_0(t)} \left\| s^{-1/r_1} b_1(s) \right\|_{r_1,(0,t)} + \left\| s^{-1-1/r_1} \frac{b_1(s) a_1(s)}{\chi_0(s)} \right\|_{r_1,(t,\infty)}$$

$$\approx t^{-1} \frac{a_1(t)}{\chi_0(t)} \left( \left\| s^{-1/r_1} b_1(s) \right\|_{r_1,(0,t)} + b_1(t) \right) \approx t^{-1} \frac{a_1(t) \left\| s^{-1/r_1} b_1(s) \right\|_{r_1,(0,t)}}{\chi_0(t)} \prec t^{-1} \frac{\chi_1(t)}{\chi_0(t)} = \frac{1}{\rho(t)}.$$

Here we have used estimate (25). Repeating the arguments of the proof of Theorem 30, we get

$$h(t) := \left\| \chi_{(0,t)}(u) \frac{u}{g_1(u) + h_1(u)} \right\|_{\Phi_0} \prec t \frac{\chi_0(t)}{\chi_1(t)} = \rho(t).$$

Hence,

$$h(t) \prec \rho(t) \prec g(t).$$

In addition, we have

$$\frac{g_0(t)}{g_1(t) + h_1(t)} \approx \frac{g_0(t) + h_0(t)}{h_1(t)} \approx \frac{t \chi_0(t)}{\chi_1(t)} = \rho(t).$$

Note that

$$g_0(t) + \rho(t) h_1(t) \approx t \chi_0(t) + t \frac{\chi_0(t)}{\chi_1(t)} \chi_1(t) \approx t \chi_0(t).$$

Theorem 4 (Case 1) of [2] completes the proof. □

**Theorem** 46. Let $0 < r, r_0, q_0 \leq \infty$, $a, a_0, b_0 \in SV$, and

$$\left\| s^{-1/r_0} b_0(s) \left\| u^{-1/q_0} a_0(u) \right\|_{q_0,(1,s)} \right\|_{r_0,(1,\infty)} < \infty.$$

Put $\chi_0(t) = \left\| s^{-1/r_0} b_0(s) \left\| u^{-1/q_0} a_0(u) \right\|_{q_0,(t,s)} \right\|_{r_0,(t,\infty)}$ and $\rho(t) = t \frac{\chi_0(t)}{\chi_1(t)}$. For given $0 \leq \theta \leq 1$, denote $a^\# = \chi_0^{1-\theta} \chi_1^\theta a \circ \rho$.

(i) If $0 < \theta < 1$, then

$$\left( \bar{A}^{\mathcal{L}}_{0,r_0,b_0,q_0,a_0}, \bar{A}^{\mathcal{R}}_{1,r_1,b_1,q_1,a_1} \right)_{\theta,r,a} = \bar{A}_{\theta,r,a^\#}.$$

(ii) If $\left\| s^{-1/r} a(s) \right\|_{r,(1,\infty)} < \infty$, then

$$\left( \bar{A}^{\mathcal{L}}_{0,r_0,b_0,q_0,a_0}, \bar{A}^{\mathcal{R}}_{1,r_1,b_1,q_1,a_1} \right)_{0,r,a} = \bar{A}_{0,r,a^\#} \cap \bar{A}^{\mathcal{L}}_{0,r,B_0,q_0,a_0} \cap \bar{A}^{\mathcal{L},\mathcal{L}}_{0,r,a \circ \rho,r_0,b_0,q_0,a_0},$$

where $B_0(t) = \left\| s^{-1/r_0} b_0(s) \right\|_{r_0,(t,\infty)} a(\rho(t))$.

(iii) If $\left\| s^{-1/r} a(s) \right\|_{r,(0,1)} < \infty$, then

$$\left( \bar{A}^{\mathcal{L}}_{0,r_0,b_0,q_0,a_0}, \bar{A}^{\mathcal{R}}_{1,r_1,b_1,q_1,a_1} \right)_{1,r,a} = \bar{A}_{1,r,a^\#} \cap \bar{A}^{\mathcal{R}}_{1,r,B_1,q_1,a_1} \cap \bar{A}^{\mathcal{R},\mathcal{R}}_{1,r,a \circ \rho,r_1,b_1,q_1,a_1},$$

where $B_1(t) = \left\| s^{-1/r_1} b_1(s) \right\|_{r_1,(0,t)} a(\rho(t))$.

*Proof.* First, note that due to Lemma 5, Lemma 11, and Lemma 12 the spaces $\bar{A}^{\mathcal{L}}_{0,r,B_0,r_0,b_0}$ and $\bar{A}^{\mathcal{R}}_{1,r,B_1,q_1,a_1}$ are an intermediate space for the couple $(A_0, A_1)$. Denote $\bar{K}(t,f) =$

$K(t, f; \bar{A}^{\mathcal{L}}_{0,r_0,b_0,q_0,a_0}, \bar{A}^{\mathcal{R}}_{1,r_1,b_1,q_1,a_1})$ and $Z = (\bar{A}^{\mathcal{L}}_{0,r_0,b_0,q_0,a_0}, \bar{A}^{\mathcal{R}}_{1,r_1,b_1,q_1,a_1})_{\theta,r,a}$ $(0 \le \theta \le 1)$. By a change of variables and Theorem 45 we can write
$$\|f\|_Z \approx \left\|\rho(t)^{-\theta} t^{-1/r} a(\rho(t)) \bar{K}(\rho(t), f)\right\|_{r,(0,\infty)} \approx I_1(\theta) + I_2(\theta) + I_3(\theta) + I_4(\theta) + I_5(\theta),$$
where
$$I_1(\theta) := \left\| \rho(t)^{-\theta} t^{-1/r} a(\rho(t)) \left\| s^{-1/r_0} b_0(s) \| u^{-1/q_0} a_0(u) K(u,f) \|_{q_0,(0,s)} \right\|_{r_0,(0,t)} \right\|_{r,(0,\infty)}$$
$$= \left\| t^{-\theta-1/r} a(\rho(t)) \left(\frac{\chi_1(t)}{\chi_0(t)}\right)^{\theta} \left\| s^{-1/r_0} b_0(s) \| u^{-1/q_0} a_0(u) K(u,f) \|_{q_0,(0,s)} \right\|_{r_0,(0,t)} \right\|_{r,(0,\infty)},$$
$$I_2(\theta) := \left\| \rho(t)^{-\theta} t^{-1/r} a(\rho(t)) \| s^{-1/r_0} b_0(s) \|_{r_0,(t,\infty)} \| u^{-1/q_0} a_0(u) K(u,f) \|_{q_0,(0,t)} \right\|_{r,(0,\infty)}$$
$$= \left\| t^{-\theta-1/r} a(\rho(t)) \left(\frac{\chi_1(t)}{\chi_0(t)}\right)^{\theta} \| s^{-1/r_0} b_0(s) \|_{r_0,(t,\infty)} \| u^{-1/q_0} a_0(u) K(u,f) \|_{q_0,(0,t)} \right\|_{r,(0,\infty)},$$
$$I_3(\theta) := \left\| \rho(t)^{1-\theta} t^{-\frac{1}{r}} a(\rho(t)) \left\| s^{-\frac{1}{r_1}} b_1(s) \right\|_{r_1,(0,t)} \left\| u^{-1-\frac{1}{q_1}} a_1(u) K(u,f) \right\|_{q_1,(t,\infty)} \right\|_{r,(0,\infty)}$$
$$= \left\| t^{(1-\theta)-\frac{1}{r}} a(\rho(t)) \left(\frac{\chi_1(t)}{\chi_0(t)}\right)^{\theta-1} \left\| s^{-\frac{1}{r_1}} b_1(s) \right\|_{r_1,(0,t)} \left\| u^{-1-\frac{1}{q_1}} a_1(u) K(u,f) \right\|_{q_1,(t,\infty)} \right\|_{r,(0,\infty)},$$
$$I_4(\theta) := \left\| \rho(t)^{1-\theta} t^{-\frac{1}{r}} a(\rho(t)) \left\| s^{-\frac{1}{r_1}} b_1(s) \left\| u^{-1-\frac{1}{q_1}} a_1(u) K(u,f) \right\|_{q_1,(s,\infty)} \right\|_{r_1,(t,\infty)} \right\|_{r,(0,\infty)}$$
$$= \left\| t^{(1-\theta)-\frac{1}{r}} a(\rho(t)) \left(\frac{\chi_1(t)}{\chi_0(t)}\right)^{\theta-1} \left\| s^{-\frac{1}{r_1}} b_1(s) \left\| u^{-1-\frac{1}{q_1}} a_1(u) K(u,f) \right\|_{q_1,(s,\infty)} \right\|_{r_1,(t,\infty)} \right\|_{r,(0,\infty)},$$
and
$$I_5(\theta) := \left\| \rho(t)^{-\theta} t^{-1/r} a(\rho(t)) \chi_0(t) K(t,f) \right\|_{r,(0,\infty)}$$
$$= \left\| t^{-\theta-1/r} \left(\frac{\chi_1(t)}{\chi_0(t)}\right)^{\theta} a(\rho(t)) \chi_0(t) K(t,f) \right\|_{r,(0,\infty)}$$
$$= \left\| t^{-\theta-1/r} \chi_0(t)^{1-\theta} \chi_1(t)^{\theta} a(\rho(t)) K(t,f) \right\|_{r,(0,\infty)} = \|f\|_{\theta,r,a^{\#}}.$$
Repeating the proof of Theorem 31, we conclude that
$$I_1(0) = \|f\|_{\bar{A}^{\mathcal{L},\mathcal{L}}_{0,r,a^{\circ}\rho,r_0,b_0,q_0,a_0}},$$
$$I_1(\theta) \prec \left\| t^{-\theta-1/r} a(\rho(t)) \chi_0(t)^{1-\theta} \chi_1(t)^{\theta} K(t,f) \right\|_{r,(0,\infty)} = \|f\|_{\theta,r;a^{\#}}, \text{ provided } 0 < \theta \le 1,$$
$$I_2(0) = \|f\|_{\mathcal{L};0,r,B_0,q_0,a_0}, \text{ where } B_0(t) = \left\| s^{-1/r_0} b_0(s) \right\|_{r_0,(t,\infty)} a(\rho(t)),$$
and
$$I_2(\theta) \prec \|f\|_{\theta,r;a^{\#}}, \text{ provided } 0 < \theta \le 1.$$
Repeating the proof of Theorem 42, we get
$$I_3(\theta) \prec \|f\|_{\theta,r;a^{\#}}, \text{ provided } 0 \le \theta < 1,$$
$$I_3(1) = \|f\|_{\mathcal{R};1,r,B_1,q_1,a_1}, \text{ where } B_1(t) = \left\| s^{-1/r_1} b_1(s) \right\|_{r_1,(0,t)} a(\rho(t)),$$
$$I_4(\theta) \prec \|f\|_{\theta,r;a^{\#}}, \text{ provided } 0 \le \theta < 1,$$
and
$$I_4(1) = \|f\|_{\mathcal{R},\mathcal{R};1,r,a^{\circ}\rho,r_1,b_1,q_1,a_1}.$$
This completes the proof. □

**Theorem** 47. Let $0 < \theta_0 < 1$, $0 < r_0, q_0 \le \infty$, $a_0, b_0 \in SV$, and $\left\| s^{-1/r_0} b_0(s) \right\|_{r_0,(1,\infty)} < \infty$. Put

$$\chi_0(t) = a_0(t)\|s^{-1/r_0}b_0(s)\|_{r_0,(t,\infty)} \text{ and } \rho(t) = t^{1-\theta_0}\frac{\chi_0(t)}{\chi_1(t)}.$$

Then, for all $f \in \bar{A}^{\mathcal{L}}_{\theta_0,r_0,b_0,q_0,a_0} + \bar{A}^{\mathcal{R}}_{1,r_1,b_1,q_1,a_1}$ and $t > 0$,

$$K\big(\rho(t), f; \bar{A}^{\mathcal{L}}_{\theta_0,r_0,b_0,q_0,a_0}, \bar{A}^{\mathcal{R}}_{1,r_1,b_1,q_1,a_1}\big)$$
$$\approx \Big\|\|s^{-1/r_0}b_0(s)\|\|u^{-\theta_0-1/q_0}a_0(u)K(u,f)\|_{q_0,(0,s)}\Big\|_{r_0,(0,t)}$$
$$+ \|s^{-1/r_0}b_0(s)\|_{r_0,(t,\infty)}\|u^{-\theta_0-1/q_0}a_0(u)K(u,f)\|_{q_0,(0,t)}$$
$$+\rho(t)\bigg(\|s^{-1/r_1}b_1(s)\|_{r_1,(0,t)}\|u^{-1-1/q_1}a_1(u)K(u,f)\|_{q_1,(t,\infty)} +$$
$$\Big\|\|s^{-1/r_1}b_1(s)\|\|u^{-1-1/q_1}a_1(u)K(u,f)\|_{q_1,(s,\infty)}\Big\|_{r_1,(t,\infty)}\bigg).$$

*Proof.* Let $\Phi_0$ be the function space corresponding to $\bar{A}^{\mathcal{L}}_{\theta_0,r_0,b_0,q_0,a_0}$:
$$\|F\|_{\Phi_0} = \Big\|\|s^{-1/r_0}b_0(s)\|\|u^{-\theta_0-1/q_0}a_0(u)F(u)\|_{q_0,(0,s)}\Big\|_{r_0,(0,\infty)}.$$

From the proof of Theorem 33, we know that
$$I(t,f) := \|\chi_{(0,t)}(*)K(*,f)\|_{\Phi_0} \approx \Big\|\|s^{-1/r_0}b_0(s)\|\|u^{-\theta_0-1/q_0}a_0(u)K(u,f)\|_{q_0,(0,s)}\Big\|_{r_0,(0,t)} +$$
$$\|s^{-1/r_0}b_0(s)\|_{r_0,(t,\infty)}\|u^{-\theta_0-1/q_0}a_0(u)K(u,f)\|_{q_0,(0,t)}$$

and
$$g_0(t) \prec h_0(t) \approx g_0(t) + h_0(t) \approx t^{1-\theta_0}\chi_0(t).$$

Now we estimate the functions $g(t)$ and $h(t)$. Repeating the arguments of the proof of Theorem 43 and Theorem 33, we get
$$h(t) \prec \rho(t) \prec g(t).$$

Finally, note that
$$\frac{h_0(t)}{h_1(t)} \approx \frac{t^{1-\theta_0}b_0(t)}{\chi_1(t)} = \rho(t).$$

Theorem 4 (Case 4) of [2] completes the proof. □

**Theorem** 48. Let $0 < \theta_0 < 1$, $0 < r, r_0, q_0 \leq \infty$, $a, a_0, b_0 \in SV$, and $\|s^{-1/r_0}b_0(s)\|_{r_0,(1,\infty)} < \infty$. Put $\chi_0(t) = a_0(t)\|s^{-1/r_0}b_0(s)\|_{r_0,(t,\infty)}$ and $\rho(t) = t^{1-\theta_0}\frac{\chi_0(t)}{\chi_1(t)}$. For a given $0 \leq \theta \leq 1$, denote $a^{\#} = \chi_0^{1-\theta}\chi_1^{\theta}a \circ \rho$ and $\eta = (1-\theta)\theta_0 + \theta$.

(i) If $0 < \theta < 1$, then
$$\big(\bar{A}^{\mathcal{L}}_{\theta_0,r_0,b_0,q_0,a_0}, \bar{A}^{\mathcal{R}}_{1,r_1,b_1,q_1,a_1}\big)_{\theta,r,a} = \bar{A}_{\eta,r,a^{\#}}.$$

(ii) If $\|s^{-1/r}a(s)\|_{r,(1,\infty)} < \infty$, then
$$\big(\bar{A}^{\mathcal{L}}_{\theta_0,r_0,b_0,q_0,a_0}, \bar{A}^{\mathcal{R}}_{1,r_1,b_1,q_1,a_1}\big)_{0,r,a} = \bar{A}^{\mathcal{L}}_{\theta_0,r,B_0,q_0,a_0} \cap \bar{A}^{\mathcal{L},\mathcal{L}}_{\theta_0,r,a\circ\rho,r_0,b_0,q_0,a_0},$$
where $B_0(t) = \|s^{-1/r_0}b_0(s)\|_{r_0,(t,\infty)}a(\rho(t))$.

(iii) If $\|s^{-1/r}a(s)\|_{r,(0,1)} < \infty$, then
$$\big(\bar{A}^{\mathcal{L}}_{\theta_0,r_0,b_0,q_0,a_0}, \bar{A}^{\mathcal{R}}_{1,r_1,b_1,q_1,a_1}\big)_{1,r,a} = \bar{A}_{1,r,a^{\#}} \cap \bar{A}^{\mathcal{R}}_{1,r,B_1,q_1,a_1} \cap \bar{A}^{\mathcal{R},\mathcal{R}}_{1,r,a\circ\rho,r_1,b_1,q_1,a_1},$$
where $B_1(t) = \|s^{-1/r_1}b_1(s)\|_{r_1,(0,t)}a(\rho(t))$.

*Proof.* First, note that due to Lemma 5, Lemma 11, and Lemma 12 the spaces $\bar{A}^{\mathcal{L}}_{\theta_0,r,B_0,r_0,a_0}$ and $\bar{A}^{\mathcal{R}}_{1,r,B_1,q_1,a_1}$ are an intermediate space for the couple $(A_0, A_1)$. Denote $\bar{K}(t,f) = K\big(t, f; \bar{A}^{\mathcal{L}}_{\theta_0,r_0,b_0,q_0,a_0}, \bar{A}^{\mathcal{R}}_{1,r_1,b_1,q_1,a_1}\big)$ and $Z = \big(\bar{A}^{\mathcal{L}}_{\theta_0,r_0,b_0,q_0,a_0}, \bar{A}^{\mathcal{R}}_{1,r_1,b_1,q_1,a_1}\big)_{\theta,r,a}$ $(0 \leq \theta \leq 1)$. By a change of variables and Theorem 47 we can write

$$\|f\|_Z \approx \left\|\rho(t)^{-\theta} t^{-1/r} a(\rho(t)) \overline{K}(\rho(t), f)\right\|_{r,(0,\infty)} \approx I_1(\theta) + I_2(\theta) + I_3(\theta) + I_4(\theta),$$

where

$$I_1(\theta) := \left\|\rho(t)^{-\theta} t^{-1/r} a(\rho(t)) \left\|s^{-1/r_0} b_0(s) \left\|u^{-\theta_0 - 1/q_0} a_0(u) K(u,f)\right\|_{q_0,(0,s)}\right\|_{r_0,(0,t)}\right\|_{r,(0,\infty)}$$

$$= \left\|t^{\theta_0 - \eta - 1/r} a(\rho(t)) \left(\frac{\chi_1(t)}{\chi_0(t)}\right)^{\theta} \left\|s^{-1/r_0} b_0(s) \left\|u^{-\theta_0 - 1/q_0} a_0(u) K(u,f)\right\|_{q_0,(0,s)}\right\|_{r_0,(0,t)}\right\|_{r,(0,\infty)},$$

$$I_2(\theta) := \left\|\rho(t)^{-\theta} t^{-1/r} a(\rho(t)) \left\|s^{-1/r_0} b_0(s)\right\|_{r_0,(t,\infty)} \left\|u^{-\theta_0 - 1/q_0} a_0(u) K(u,f)\right\|_{q_0,(0,t)}\right\|_{r,(0,\infty)}$$

$$= \left\|t^{\theta_0 - \eta - 1/r} a(\rho(t)) \left(\frac{\chi_1(t)}{\chi_0(t)}\right)^{\theta} \left\|s^{-1/r_0} b_0(s)\right\|_{r_0,(t,\infty)} \left\|u^{-\theta_0 - 1/q_0} a_0(u) K(u,f)\right\|_{q_0,(0,t)}\right\|_{r,(0,\infty)},$$

$$I_3(\theta) := \left\|\rho(t)^{1-\theta} t^{-1/r} a(\rho(t)) \left\|s^{-1/r_1} b_1(s)\right\|_{r_1,(0,t)} \left\|u^{-1-1/q_1} a_1(u) K(u,f)\right\|_{q_1,(t,\infty)}\right\|_{r,(0,\infty)}$$

$$= \left\|t^{(1-\theta)(1-\theta_0)-1/r} a(\rho(t)) \left(\frac{\chi_1(t)}{\chi_0(t)}\right)^{\theta-1} \left\|s^{-1/r_1} b_1(s)\right\|_{r_1,(0,t)} \left\|u^{-1-1/q_1} a_1(u) K(u,f)\right\|_{q_1,(t,\infty)}\right\|_{r,(0,\infty)}$$

$$= \left\|t^{1-\eta-1/r} a(\rho(t)) \left(\frac{\chi_1(t)}{\chi_0(t)}\right)^{\theta-1} \left\|s^{-1/r_1} b_1(s)\right\|_{r_1,(0,t)} \left\|u^{-1-1/q_1} a_1(u) K(u,f)\right\|_{q_1,(t,\infty)}\right\|_{r,(0,\infty)},$$

and

$$I_4(\theta) := \left\|\rho(t)^{1-\theta} t^{-\frac{1}{r}} a(\rho(t)) \left\|s^{-1/r_1} b_1(s) \left\|u^{-1-1/q_1} a_1(u) K(u,f)\right\|_{q_1,(s,\infty)}\right\|_{r_1,(t,\infty)}\right\|_{r,(0,\infty)}$$

$$= \left\|t^{(1-\theta)(1-\theta_0)-\frac{1}{r}} a(\rho(t)) \left(\frac{\chi_1(t)}{\chi_0(t)}\right)^{\theta-1} \left\|s^{-1/r_1} b_1(s) \left\|u^{-1-1/q_1} a_1(u) K(u,f)\right\|_{q_1,(s,\infty)}\right\|_{r_1,(t,\infty)}\right\|_{r,(0,\infty)}$$

$$= \left\|t^{1-\eta-\frac{1}{r}} a(\rho(t)) \left(\frac{\chi_1(t)}{\chi_0(t)}\right)^{\theta-1} \left\|s^{-1/r_1} b_1(s) \left\|u^{-1-1/q_1} a_1(u) K(u,f)\right\|_{q_1,(s,\infty)}\right\|_{r_1,(t,\infty)}\right\|_{r,(0,\infty)}.$$

Repeating the proof of Theorem 34, we conclude that
$$I_1(0) = \|f\|_{\bar{A}^{\mathcal{L},\mathcal{L}}_{\theta_0,r,a°\rho,r_0,b_0,q_0,a_0}},$$
$$I_1(\theta) \prec \|f\|_{\eta,r;a^\#}, \text{ provided } 0 < \theta \leq 1,$$
$$I_2(0) = \|f\|_{\mathcal{L};\theta_0,r,B_0,q_0,a_0}, \text{ where } B_0(t) = \left\|s^{-1/r_0} b_0(s)\right\|_{r_0,(t,\infty)} a(\rho(t)),$$

and
$$I_2(\theta) \approx \|f\|_{\eta,r;a^\#}, \text{ provided } 0 < \theta \leq 1.$$

Repeating the proof of Theorem 44, we conclude that
$$I_3(\theta) \prec \|f\|_{\eta,r;a^\#}, \text{ provided } 0 \leq \theta < 1,$$
$$I_3(1) = \|f\|_{\mathcal{R};1,r,B_1,q_1,a_1}, \text{ where } B_1(t) = \left\|s^{-1/r_1} b_1(s)\right\|_{r_1,(0,t)} a(\rho(t)),$$
$$I_4(\theta) \prec \|f\|_{\eta,r;a^\#}, \text{ provided } 0 \leq \theta < 1,$$

and
$$I_4(1) = \|f\|_{\mathcal{R},\mathcal{R};1,r,a°\rho,r_1,b_1,q_1,a_1}.$$

Thus (i) and (iii) are proved. Consider the case $\theta = 0$. Because
$$B_0(t) a_0(t) = a_0(t) \left\|s^{-1/r_0} b_0(s)\right\|_{r_0,(t,\infty)} a(\rho(t)) = a^\#(t),$$

by Corollary 17, we get $\bar{A}^{\mathcal{L}}_{\theta_0,r,B_0,q_0,a_0} \subset \bar{A}_{\theta_0,r,a^\#}$. This means, that
$$\max(I_3(0), I_4(0)) \prec \|f\|_{\theta_0,r;a^\#} \prec I_2(0).$$

This completes the proof. □

**Theorem** 49. Let $0 < \theta_0 < 1$, $0 < r_0, q_0 \le \infty$, $a_0, b_0 \in SV$, and $\|s^{-1/r_0} b_0(s)\|_{r_0,(0,1)} < \infty$. Put $\chi_0(t) = a_0(t)\|s^{-1/r_0} b_0(s)\|_{r_0,(0,t)}$ and $\rho(t) = t^{1-\theta_0} \frac{\chi_0(t)}{\chi_1(t)}$. Then, for all $f \in \bar{A}^{\mathcal{R}}_{\theta_0,r_0,b_0,q_0,a_0} + \bar{A}^{\mathcal{R}}_{1,r_1,b_1,q_1,a_1}$ and $t > 0$,

$$K(\rho(t), f; \bar{A}^{\mathcal{R}}_{\theta_0,r_0,b_0,q_0,a_0}, \bar{A}^{\mathcal{R}}_{1,r_1,b_1,q_1,a_1})$$
$$\approx \left\| s^{-1/r_0} b_0(s) \left\| u^{-\theta_0 - 1/q_0} a_0(u) K(u, f) \right\|_{q_0,(s,t)} \right\|_{r_0,(0,t)}$$
$$+ \rho(t) \left( \left\| s^{-1/r_1} b_1(s) \right\|_{r_1,(0,t)} \left\| u^{-1-1/q_1} a_1(u) K(u, f) \right\|_{q_1,(t,\infty)} + \right.$$
$$\left. \left\| s^{-1/r_1} b_1(s) \left\| u^{-1-1/q_1} a_1(u) K(u, f) \right\|_{q_1,(s,\infty)} \right\|_{r_1,(t,\infty)} \right). \tag{27}$$

*Proof.* Let $\Phi_0$ be the function space corresponding to $\bar{A}^{\mathcal{R}}_{\theta_0,r_0,b_0,q_0,a_0}$:
$$\|F\|_{\Phi_0} = \left\| s^{-1/r_0} b_0(s) \left\| u^{-\theta_0 - 1/q_0} a_0(u) F(u) \right\|_{q_0,(s,\infty)} \right\|_{r_0,(0,\infty)}.$$
We know (see the proof of Theorem 37 [9]) that
$$I(t, f) = \left\| s^{-1/r_0} b_0(s) \left\| u^{-\theta_0 - 1/q_0} a_0(u) K(u, f) \right\|_{q_0,(s,t)} \right\|_{r_0,(0,t)}$$
and
$$h_0(t) \prec g_0(t) \approx g_0(t) + h_0(t) \approx t^{1-\theta_0} \chi_0(t).$$
Now we estimate the functions $g(t)$ and $h(t)$. Repeating the arguments of the proofs of Theorem 47 and Theorem 37, we get
$$\frac{1}{g(t)} := \left\| \chi_{(t,\infty)}(u) \frac{u}{g_0(u) + h_0(u)} \right\|_{\Phi_1} \prec \frac{1}{\rho(t)}$$
and
$$h(t) := \left\| \chi_{(0,t)}(u) \frac{u}{g_1(u) + h_1(u)} \right\|_{\Phi_0} \approx \rho(t).$$
Hence, we have
$$h(t) \prec \rho(t) \prec g(t).$$
Note that
$$g_0(t) + h_1(t)\rho(t) \approx t^{1-\theta_0} \chi_0(t) + \chi_1(t) t^{1-\theta_0} \frac{\chi_0(t)}{\chi_1(t)} \approx t^{1-\theta_0} \chi_0(t).$$
Thus, by [2, Theorem 4, Case 1] we arrive at
$$K(\rho(t), f; \bar{A}^{\mathcal{R}}_{\theta_0,r_0,b_0,q_0,a_0}, \bar{A}^{\mathcal{R}}_{1,r_1,b_1,q_1,a_1})$$
$$\approx \left\| s^{-1/r_0} b_0(s) \left\| u^{-\theta_0 - 1/q_0} a_0(u) K(u, f) \right\|_{q_0,(s,t)} \right\|_{r_0,(0,t)}$$
$$+ \rho(t) \left( \left\| s^{-1/r_1} b_1(s) \right\|_{r_1,(0,t)} \left\| u^{-1-1/q_1} a_1(u) K(u, f) \right\|_{q_1,(t,\infty)} + \right.$$
$$\left. \left\| s^{-1/r_1} b_1(s) \left\| u^{-1-1/q_1} a_1(u) K(u, f) \right\|_{q_1,(s,\infty)} \right\|_{r_1,(t,\infty)} \right) + t^{-\theta_0} \chi_0(t) K(t, f).$$
Repeating the arguments of the proof of Theorem 37, we can show that
$$t^{-\theta_0} \chi_0(t) K(t, f) \prec \left\| s^{-1/r_0} b_0(s) \left\| u^{-\theta_0 - \frac{1}{q_0}} a_0(u) K(u, f) \right\|_{q_0,(s,t)} \right\|_{r_0,(0,t)}.$$
Finally, we obtain formula (27). □

**Theorem** 50. Let $0 < \theta_0 < 1$, $0 < r, r_0, q_0 \le \infty$, $a, a_0, b_0 \in SV$, and $\|s^{-1/r_0} b_0(s)\|_{r_0,(0,1)} < \infty$. Put $\chi_0(t) = a_0(t)\|s^{-1/r_0} b_0(s)\|_{r_0,(0,t)}$ and $\rho(t) = t^{1-\theta_0} \frac{\chi_0(t)}{\chi_1(t)}$. For a given $0 \le \theta \le 1$, denote $a^{\#} = \chi_0^{1-\theta} \chi_1^{\theta} a \circ \rho$ and $\eta = (1-\theta)\theta_0 + \theta$.

(i) If $0 < \theta < 1$, then
$$\left(\bar{A}^{\mathcal{R}}_{\theta_0,r_0,b_0,q_0,a_0}, \bar{A}^{\mathcal{R}}_{1,r_1,b_1,q_1,a_1}\right)_{\theta,r,a} = \bar{A}_{\eta,r,a^\#}.$$

(ii) If $\left\|s^{-1/r}a(s)\right\|_{r,(1,\infty)} < \infty$, then
$$\left(\bar{A}^{\mathcal{R}}_{\theta_0,r_0,b_0,q_0,a_0}, \bar{A}^{\mathcal{R}}_{1,r_1,b_1,q_1,a_1}\right)_{0,r,a} = \bar{A}^{\mathcal{R},\mathcal{L}}_{\theta_0,r,a^\circ\rho,r_0,b_0,q_0,a_0}.$$

(iii) If $\left\|s^{-1/r}a(s)\right\|_{r,(0,1)} < \infty$, then
$$\left(\bar{A}^{\mathcal{R}}_{\theta_0,r_0,b_0,q_0,a_0}, \bar{A}^{\mathcal{R}}_{1,r_1,b_1,q_1,a_1}\right)_{1,r,a} = \bar{A}_{1,r,a^\#} \cap \bar{A}^{\mathcal{R}}_{1,r,B_1,q_1,a_1} \cap \bar{A}^{\mathcal{R},\mathcal{R}}_{1,r,a^\circ\rho,r_1,b_1,q_1,a_1},$$

where $B_1(t) = \left\|s^{-1/r_1}b_1(s)\right\|_{r_1,(0,t)} a(\rho(t))$.

*Proof.* First, note that due to Lemma 5 (ii) and Lemma 12 the space $\bar{A}^{\mathcal{R}}_{1,r,B_1,q_1,a_1}$ is an intermediate space for the couple $(A_0, A_1)$. Denote $K(t,f) = K(t,f; \bar{A}^{\mathcal{R}}_{\theta_0,r_0,b_0,q_0,a_0}, \bar{A}^{\mathcal{R}}_{1,r_1,b_1,q_1,a_1})$ and $Z = \left(\bar{A}^{\mathcal{R}}_{\theta_0,r_0,b_0,q_0,a_0}, \bar{A}^{\mathcal{R}}_{1,r_1,b_1,q_1,a_1}\right)_{\theta,r,a}$ $(0 \le \theta \le 1)$. By a change of variables and Theorem 49 we can write
$$\|f\|_Z \approx \left\|\rho(t)^{-\theta} t^{-1/r} a(\rho(t)) \bar{K}(\rho(t), f)\right\|_{r,(0,\infty)} \approx I_1(\theta) + I_2(\theta) + I_3(\theta),$$

where
$$I_1(\theta) := \left\|\rho(t)^{-\theta} t^{-1/r} a(\rho(t)) \left\|s^{-1/r_0}b_0(s)\left\|u^{-\theta_0-1/q_0}a_0(u)K(u,f)\right\|_{q_0,(s,t)}\right\|_{r_0,(0,t)}\right\|_{r,(0,\infty)}$$
$$= \left\|t^{\theta_0-\eta-1/r} a(\rho(t)) \left(\frac{\chi_1(t)}{\chi_0(t)}\right)^\theta \left\|s^{-1/r_0}b_0(s)\left\|u^{-\theta_0-1/q_0}a_0(u)K(u,f)\right\|_{q_0,(s,t)}\right\|_{r_0,(0,t)}\right\|_{r,(0,\infty)},$$

$$I_2(\theta) := \left\|\rho(t)^{1-\theta} t^{-1/r} a(\rho(t)) \left\|s^{-1/r_1}b_1(s)\right\|_{r_1,(0,t)} \left\|u^{-1-1/q_1}a_1(u)K(u,f)\right\|_{q_1,(t,\infty)}\right\|_{r,(0,\infty)}$$
$$= \left\|t^{1-\eta-1/r} a(\rho(t)) \left(\frac{\chi_1(t)}{\chi_0(t)}\right)^{\theta-1} \left\|s^{-1/r_1}b_1(s)\right\|_{r_1,(0,t)} \left\|u^{-1-1/q_1}a_1(u)K(u,f)\right\|_{q_1,(t,\infty)}\right\|_{r,(0,\infty)},$$

and
$$I_3(\theta) := \left\|\rho(t)^{1-\theta} t^{-\frac{1}{r}} a(\rho(t)) \left\|s^{-1/r_1}b_1(s) \left\|u^{-1-1/q_1}a_1(u)K(u,f)\right\|_{q_1,(s,\infty)}\right\|_{r_1,(t,\infty)}\right\|_{r,(0,\infty)}$$
$$= \left\|t^{1-\eta-\frac{1}{r}} a(\rho(t)) \left(\frac{\chi_1(t)}{\chi_0(t)}\right)^{\theta-1} \left\|s^{-1/r_1}b_1(s)\left\|u^{-1-1/q_1}a_1(u)K(u,f)\right\|_{q_1,(s,\infty)}\right\|_{r_1,(t,\infty)}\right\|_{r,(0,\infty)}.$$

Repeating the proofs of Theorem 38 and Theorem 44, we conclude that
$$I_1(0) = \|f\|_{\bar{A}^{\mathcal{R},\mathcal{L}}_{\theta_0,r,a^\circ\rho,r_0,b_0,q_0,a_0}},$$
$$I_1(\theta) \approx \|f\|_{\eta,r,a^\#}, \text{ provided } 0 < \theta \le 1,$$
$$I_2(\theta) \prec \|f\|_{\eta,r;a^\#}, \text{ provided } 0 \le \theta < 1,$$
$$I_2(1) = \|f\|_{\mathcal{R};1,r,B_1,q_1,a_1}, \quad \text{where } B_1(t) = \left\|s^{-1/r_1}b_1(s)\right\|_{r_1,(0,t)} a(\rho(t)),$$
$$I_3(\theta) \prec \|f\|_{\eta,r;a^\#}, \text{ provided } 0 \le \theta < 1,$$

and
$$I_3(1) = \|f\|_{\mathcal{R},\mathcal{R};1,r,a^\circ\rho,r_1,b_1,q_1,a_1}.$$

Thus (i) and (iii) are proved. Consider the case $\theta = 0$. Because
$$(a^\circ \rho)(t) a_0(t) \left\|s^{-1/r_0}b_0(s)\right\|_{r_0,(0,t)} = (\chi_0 a^\circ \rho)(t) = a^\#(t),$$

by Corollary 19 (v), we get $\bar{A}^{\mathcal{R},\mathcal{L}}_{\theta_0,r,a^\circ\rho,r_0,b_0,q_0,a_0} \subset \bar{A}_{\theta_0,r,a^\#}$. This means, that
$$\max(I_2(0), I_3(0)) \prec \|f\|_{\theta_0,r;a^\#} \prec I_1(0).$$

This completes the proof. □

## 5. The limiting case $\theta_0 = 0$

Using symmetry arguments based on formulae (6), (7), and (8), the counterparts to all theorems from previous section can be formulated and proved. See, e.g., the proof of Theorems 3.4 and 4.7 in [16] or the proof of Theorem 34 in [9]. Thereby, it is not hard to prove Holmstedt-type and reiteration formulae for the combinations $(\bar{A}^{\mathcal{L}}_{0,*},*)$ and $(\bar{A}^{R}_{0,*},*)$ corresponding to the limiting case $\theta_0 = 0$. We leave this to the reader and formulate only the symmetry counterpart of Theorem 31 which completes Subsection 5.2.

**Theorem** 51. Let $0 < r, r_0, q_0, q_1, r_1 \le \infty$, $a, a_0, b_0, a_1, b_1 \in SV$, $\left\|s^{-1/r_0}b_0(s)\right\|_{r_0,(0,1)} < \infty$,

$\left\|s^{-1/r_0}b_0(s)\left\|u^{-1/q_0}a_0(u)\right\|_{q_0,(s,\infty)}\right\|_{r_0,(1,\infty)} < \infty$, and

$\left\|s^{-1/r_1}b_1(s)\left\|u^{-1/q_1}a_1(u)\right\|_{q_1,(s,1)}\right\|_{r_1,(0,1)} < \infty$. Put

$$\chi_0(t) = \left\|s^{-1/r_0}b_0(s)\left\|u^{-1/q_0}a_0(u)\right\|_{q_0,(s,\infty)}\right\|_{r_0,(t,\infty)} +$$
$$\left\|s^{-1/r_0}b_0(s)\right\|_{r_0,(0,t)}\left\|u^{-1/q_0}a_0(u)\right\|_{q_0,(t,\infty)},$$
$$\chi_1(t) = \left\|s^{-1/r_1}b_1(s)\left\|u^{-1/q_1}a_1(u)\right\|_{q_1,(s,t)}\right\|_{r_1,(0,t)},$$

and $\rho(t) = t\frac{\chi_0(t)}{\chi_1(t)}$. For given $0 \le \theta \le 1$, denote $a^{\#} = \chi_0^{(1-\theta)}\chi_1^{\theta}a\circ\rho$.

(i) If $0 < \theta < 1$, then
$$\left(\bar{A}^{\mathcal{R}}_{0,r_0,b_0,q_0,a_0}, \bar{A}^{\mathcal{R}}_{1,r_1,b_1,q_1,a_1}\right)_{\theta,r,a} = \bar{A}_{\theta,r,a^{\#}}.$$

(ii) If $\left\|s^{-1/r}a(s)\right\|_{r,(1,\infty)} < \infty$, then
$$\left(\bar{A}^{\mathcal{R}}_{0,r_0,b_0,q_0,a_0}, \bar{A}^{\mathcal{R}}_{1,r_1,b_1,q_1,a_1}\right)_{0,r,a} = \bar{A}_{0,r,a^{\#}} \cap \bar{A}^{\mathcal{R},\mathcal{L}}_{0,r,a\circ\rho,r_0,b_0,q_0,a_0}.$$

(iii) If $\left\|s^{-1/r}a(s)\right\|_{r,(0,1)} < \infty$, then
$$\left(\bar{A}^{\mathcal{R}}_{0,r_0,b_0,q_0,a_0}, \bar{A}^{\mathcal{R}}_{1,r_1,b_1,q_1,a_1}\right)_{1,r,a} = \bar{A}_{1,r,a^{\#}} \cap \bar{A}^{\mathcal{R}}_{1,r,B_1,q_1,a_1} \cap \bar{A}^{\mathcal{R},\mathcal{R}}_{1,r,a\circ\rho,r_1,b_1,q_1,a_1},$$
where $B_1(t) = \left\|s^{-1/r_1}b_1(s)\right\|_{r_1,(0,t)}a(\rho(t))$.

## 6. Applications

Here we demonstrate how our general reiteration theorems can be used to establish limiting interpolation results for the grand and small Lorentz spaces as well as for Lorentz–Karamata spaces. For the sake of shortness, we present only some possible results.

Let $(\Omega, \mu)$ denote a $\sigma$-finite measure space with a non-atomic measure $\mu$. We consider functions $f$ from the set $\mathfrak{M}(\Omega, \mu)$ of all $\mu$-measurable functions on $\Omega$. As conventional (see e.g. [4]), $f^*(t)$ ($t$>0) denotes the non-increasing rearrangement of $f$ and the maximal function of $f^*$ is defined by
$$f^{**}(t) := \tfrac{1}{t}\int_0^t f^*(u)du.$$

For simplicity, we focus our attention below on interpolation spaces between $L^1$ and $L^{\infty}$. For functions from these spaces J. Peetre has shown that [4, Theorem V.1.6]
$$K(t,f;L_1,L_{\infty}) = \int_0^t f^*(s)ds = tf^{**}(t). \qquad (28)$$

The following assertion is a modification of [20, Lemma 5.2] and can be proved similarly. See also [9] and [10].

**Lemma** 52. Let $0 < \theta \le 1$, $p = \tfrac{1}{1-\theta}$, $0 < q \le \infty$, and $b \in SV$. If $\theta = 1$, we additionally suppose that $\left\|u^{-1/q}b(u)\right\|_{q,(0,1)} < \infty$. Then, for all $f \in L_1 + L_{\infty}$ and $t > 0$,

$$\left\|u^{-\theta-\frac{1}{q}}b(u)K(u,f;L_1,L_\infty)\right\|_{q,(0,t)} = \left\|u^{\frac{1}{p}-\frac{1}{q}}b(u)f^{**}(u)\right\|_{q,(0,t)} \approx \left\|u^{\frac{1}{p}-\frac{1}{q}}b(u)f^{*}(u)\right\|_{q,(0,t)}.$$

### 6.1. Function spaces

First, we define function spaces under consideration.

**Definition 53.** Let $0<p,q\leq\infty$ and $b \in SV$. The Lorentz–Karamata spaces $L_{p,q;b}$ and $L_{(p,q;b)}$ are the sets of all $f \in \mathfrak{M}(\Omega,\mu)$ such that

$$\|f\|_{p,q;b} = \left\|t^{\frac{1}{p}-\frac{1}{q}}b(t)f^{*}(t)\right\|_{q,(0,\infty)} < \infty$$

and

$$\|f\|_{(p,q;b)} = \left\|t^{\frac{1}{p}-\frac{1}{q}}b(t)f^{**}(t)\right\|_{q,(0,\infty)} < \infty$$

respectively.

The Lorentz–Karamata spaces form an important scale of spaces. It contains e.g. the Lebesgue spaces $L_p$, Lorentz space $L_{p,q}$, Lorentz-Zygmund, and the generalized Lorentz–Zygmund spaces. For further information about Lorentz–Karamata spaces we refer to e.g. [20, 22, 15]. Lemma 52 implies the following interpolation result.

**Lemma 54.** (Cf. [9, Lemma 53; 10; 15; 20, Corollary 5.3].) Let $0 < \theta \leq 1$, $p = \frac{1}{1-\theta}$, $0 < q \leq \infty$, and $b \in SV$. Then $(L_1,L_\infty)_{\theta,q,b} = L_{p,q,b}$. Moreover, $(L_1,L_\infty)_{0,q,b} = L_{(1,q,b)}$, provided $\left\|u^{-1/q}b(u)\right\|_{q,(1,\infty)} < \infty$.

**Definition 55.** (Cf. [20, (5.21), (5.33)].) Let $0 < p,q,r \leq \infty$ and $a,b \in SV$. The spaces $L^{\mathcal{L}}_{p,r,b,q,a}$ and $L^{\mathcal{R}}_{p,r,b,q,a}$ are the sets of all $f \in \mathfrak{M}(\Omega,\mu)$ such that

$$\|f\|_{L^{\mathcal{L}}_{p,r,b,q,a}} = \left\|t^{-1/r}b(t)\left\|u^{\frac{1}{p}-\frac{1}{q}}a(u)f^{*}(u)\right\|_{q,(0,t)}\right\|_{r,(0,\infty)} < \infty$$

or

$$\|f\|_{L^{\mathcal{R}}_{p,r,b,q,a}} = \left\|t^{-1/r}b(t)\left\|u^{\frac{1}{p}-\frac{1}{q}}a(u)f^{*}(u)\right\|_{q,(t,\infty)}\right\|_{r,(0,\infty)} < \infty,$$

respectively.

Similar definitions can be found in [5, 7, 9, 12, 15, 18]. Note that the spaces $L^{\mathcal{L}}_*$ are a special cases of Generalized Gamma space with double weights [19].

**Definition 56.** Let $0 < p,q,r \leq \infty$ and $b \in SV$. The small Lorentz space $L^{(p,q,r)}_{b}$ and the grand Lorentz space $L^{p),q,r}_{b}$ we define as the sets of all $f \in \mathfrak{M}(\Omega,\mu)$ such that

$$\|f\|_{L^{(p,q,r)}_{b}} = \left\|t^{-\frac{1}{r}}b(t)\left\|u^{\frac{1}{p}-\frac{1}{q}}f^{*}(u)\right\|_{q,(0,t)}\right\|_{r,(0,\infty)} < \infty$$

or

$$\|f\|_{L^{p),q,r}_{b}} = \left\|t^{-\frac{1}{r}}b(t)\left\|u^{\frac{1}{p}-\frac{1}{q}}f^{*}(u)\right\|_{q,(t,\infty)}\right\|_{r,(0,\infty)} < \infty,$$

respectively.

It is clear, that $L^{(p,q,r)}_{b} = L^{\mathcal{L}}_{p,r,b,q,1}$ and $L^{p),q,r}_{b} = L^{\mathcal{R}}_{p,r,b,q,1}$.

**Definition 57.** [9, 10]. Let $0 < p,q,r,s \leq \infty$ and $a,b,c \in SV$. The spaces $L^{\mathcal{L},\mathcal{L}}_{p,(s,c,r,b,q,a)}$, $L^{\mathcal{R},\mathcal{L}}_{p,(s,c,r,b,q,a)}$, and $L^{\mathcal{R},\mathcal{R}}_{p,(s,c,r,b,q,a)}$, are the sets of all $f \in \mathfrak{M}(\Omega,\mu)$ such that

$$\|f\|_{L^{\mathcal{L},\mathcal{L}}_{p,(s,c,r,b,q,a)}} = \left\| u^{-1/s}c(u) \left\| t^{-1/r}b(t) \left\| v^{\frac{1}{p}-\frac{1}{q}}a(v)f^*(v) \right\|_{q,(0,t)} \right\|_{r,(0,u)} \right\|_{s,(0,\infty)} < \infty,$$

$$\|f\|_{L^{\mathcal{R},\mathcal{L}}_{p,(s,c,r,b,q,a)}} = \left\| u^{-1/s}c(u) \left\| t^{-1/r}b(t) \left\| v^{\frac{1}{p}-\frac{1}{q}}a(v)f^*(v) \right\|_{q,(t,u)} \right\|_{r,(0,u)} \right\|_{s,(0,\infty)} < \infty,$$

or

$$\|f\|_{L^{\mathcal{R},\mathcal{R}}_{p,(s,c,r,b,q,a)}} = \left\| u^{-1/s}c(u) \left\| t^{-1/r}b(t) \left\| v^{\frac{1}{p}-\frac{1}{q}}a(v)f^*(v) \right\|_{q,(t,\infty)} \right\|_{r,(u,\infty)} \right\|_{s,(0,\infty)} < \infty,$$

respectively.

The following lemmas characterise the $L^{\mathcal{L}}$, $L^{\mathcal{R}}$, $L^{\mathcal{L},\mathcal{L}}$, $L^{\mathcal{R},\mathcal{L}}$, and $L^{\mathcal{R},\mathcal{R}}$ spaces as appropriate interpolation spaces.

**Lemma 58.** [10]. (Cf. [20, Lemmas 5.4, 5.9].) Let $0 < \theta \leq 1$, $p = \frac{1}{1-\theta}$, $0 < q,r \leq \infty$, and $a,b \in SV$. Then
$$L^{\mathcal{L}}_{p,r,b,q,a} = (L_1, L_\infty)^{\mathcal{L}}_{\theta,r,b,q,a} \quad \text{and} \quad L^{\mathcal{R}}_{p,r,b,q,a} = (L_1, L_\infty)^{\mathcal{R}}_{\theta,r,b,q,a}.$$
In particular, $L^{(p,q,r}_b = (L_1, L_\infty)^{\mathcal{L}}_{\theta,r,b,q,1}$ and $L^{p),q,r}_b = (L_1, L_\infty)^{\mathcal{R}}_{\theta,r,b,q,1}$.

**Lemma 59.** [10]. Let $0 < \theta \leq 1$, $p = \frac{1}{1-\theta}$, $0 < q,r,s \leq \infty$, and $a,b,c \in SV$. Then
$$L^{\mathcal{L},\mathcal{L}}_{p,(s,c,r,b,q,a)} = (L_1, L_\infty)^{\mathcal{L},\mathcal{L}}_{\theta,s,c,r,b,q,a},$$
$$L^{\mathcal{R},\mathcal{L}}_{p,(s,c,r,b,q,a)} = (L_1, L_\infty)^{\mathcal{R},\mathcal{L}}_{\theta,s,c,r,b,q,a},$$
and
$$L^{\mathcal{R},\mathcal{R}}_{p,(s,c,r,b,q,a)} = (L_1, L_\infty)^{\mathcal{R},\mathcal{R}}_{\theta,s,c,r,b,q,a}.$$

For further information about defined spaces we refer to e.g. [10]. We are going to use the theorems from Section 5 to the couple $(L_1, L_\infty)$ to get some interpolation results for the grand and small Lorentz spaces. Due to Lemma 11 and Lemma 13 the space $L^{(\infty,q,r}_b = (L_1, L_\infty)^{\mathcal{L}}_{1,r,b,q,1}$ is trivial if $q < \infty$ or coincides with $L_\infty + \infty L_1 = L_\infty$ if $q = \infty$ [4, Ch. 5. (1.30)]. Therefore, we will use only theorems from Subsection 5.2 for the couples where the second operand is $(L_1, L_\infty)^{\mathcal{R}}_{1,r_1,b_1,q_1,1} = L^{\infty),q_1,r_1}_{b_1}$. Note that due to Lemma 58 and Lemma 12, $L^{\infty),q,r}_b$ is not trivial only if $\left\|s^{-1/r}b(s)|\ln s|^{1/q}\right\|_{r,(0,1)} < \infty$.

### 6.2. Some interpolation results for couples, where the first operand is $L_1$ or a Lorentz-Karamata space

**Corollary 60.** Let $0 < q_1, r_1, r \leq \infty$, $a, b_1 \in SV$, and $\left\|s^{-1/r_1}b_1(s)|\ln s|^{1/q_1}\right\|_{r_1,(0,1)} < \infty$. Put

$$\chi_1(t) = \left\|s^{-1/r_1}b_1(s)\left\{\ln\frac{t}{s}\right\}^{1/q_1}\right\|_{r_1,(0,t)}$$

and $\rho(t) = \frac{t}{\chi_1(t)}$. For given $0 \leq \theta \leq 1$, denote $a^\# = \chi_1^\theta a \circ \rho$ and $p = \frac{1}{1-\theta}$.

(i) If $0 < \theta < 1$, then

$$\left(L_1, L_{b_1}^{\infty),q_1,r_1}\right)_{\theta,r,\mathrm{a}} = L_{p,r,a^\#}.$$

(ii) If $\left\|s^{-1/r}a(s)\right\|_{r,(1,\infty)} < \infty$, then
$$\left(L_1, L_{b_1}^{\infty),q_1,r_1}\right)_{0,r,\mathrm{a}} = L_{(1,r,a^\#)}.$$

(iii) If $\left\|s^{-1/r}a(s)\right\|_{r,(0,1)} < \infty$, then
$$\left(L_1, L_{b_1}^{\infty),q_1,r_1}\right)_{1,r,\mathrm{a}} = L_{\infty,r,a^\#} \cap L_{B_1}^{\infty),q_1,r} \cap L_{\infty,(r,a^\circ\rho,r_1,b_1,q_1,1)}^{\mathcal{R},\mathcal{R}},$$

where $B_1(t) = \left\|s^{-1/r_1}b_1(s)\right\|_{r_1,(0,t)} a(\rho(t))$.

*Proof.* From Lemma 54 and Lemma 59, we know that $(L_1, L_\infty)_{\theta,r,a^\#} = L_{p,r,a^\#}$ ($0 < \theta \le 1$), $(L_1, L_\infty)_{0,r,a^\circ\rho} = L_{(1,r,a^\circ\rho)}$, and $(L_1, L_\infty)_{1,r,a^\circ\rho,r_1,b_1,q_1,1}^{\mathcal{R},\mathcal{R}} = L_{\infty,(r,a^\circ\rho,r_1,b_1,q_1,1)}^{\mathcal{R},\mathcal{R}}$. Thus, the proof follows from Theorem 40. □

**Corollary 61.** (Cf. [15, Corollary 7.4].) Let $0 < r_1, r \le \infty$, $a, b_1 \in SV$, and $\left\|s^{-1/r_1}b_1(s)\right\|_{r_1,(0,1)} < \infty$. Put $\chi_1(t) = \left\|s^{-1/r_1}b_1(s)\right\|_{r_1,(0,t)}$ and $\rho(t) = \frac{t}{\chi_1(t)}$. For a given $0 \le \theta \le 1$, denote $a^\# = \chi_1^\theta a^\circ \rho$ and $p = \frac{1}{1-\theta}$.

(i) If $0 < \theta < 1$, then
$$\left(L_1, L_{\infty,r_1,b_1}\right)_{\theta,r,\mathrm{a}} \cong L_{p,r,a^\#}.$$

(ii) If $\left\|s^{-1/r}a(s)\right\|_{r,(1,\infty)} < \infty$, then
$$\left(L_1, L_{\infty,r_1,b_1}\right)_{0,r,\mathrm{a}} \cong L_{(1,r,a^\circ\rho)}.$$

(iii) If $\left\|s^{-1/r}a(s)\right\|_{r,(0,1)} < \infty$, then
$$\left(L_1, L_{\infty,r_1,b_1}\right)_{1,r,\mathrm{a}} \cong L_{\infty,r,a^\#} \cap L_{\infty,r,a^\circ\rho,r_1,b_1}^{\mathcal{R}}.$$

*Proof.* It is not difficult to show that $L_{b_1}^{\infty),\infty,r_1} = L_{\infty,r_1,b_1}$. Thus, the proof follows from Corollary 60, if we take $q_1 = \infty$, because by Lemma 59, Lemma 15, and Lemma 58, we have
$$L_{\infty,(r,a^\circ\rho,r_1,b_1,\infty,1)}^{\mathcal{R},\mathcal{R}} = (L_1, L_\infty)_{1,r,a^\circ\rho,r_1,b_1,\infty,1}^{\mathcal{R},\mathcal{R}} = (L_1, L_\infty)_{1,r,a^\circ\rho,r_1,b_1}^{\mathcal{R}} = L_{\infty,r,a^\circ\rho,r_1,b_1}^{\mathcal{R}}. \quad \square$$

**Corollary 62.** (Cf. [8, Corollaries 44 and 45].) Let $1 < p_0 < \infty$, $0 < r, r_0, q_1, r_1 \le \infty$, $a, b_0, b_1 \in SV$, and $\left\|s^{-1/r_1}b_1(s)|\ln s|^{1/q_1}\right\|_{r_1,(0,1)} < \infty$. Put
$$\chi_1(t) = \left\|s^{-1/r_1}b_1(s)\left\{\ln\frac{t}{s}\right\}^{1/q_1}\right\|_{r_1,(0,t)}$$

and $\rho(t) = t^{\frac{1}{p_0}} \frac{b_0(t)}{\chi_1(t)}$. For given $0 \le \theta \le 1$, denote $a^\# = b_0^{1-\theta} \chi_1^\theta a^\circ \rho$ and $p = \frac{p_0}{1-\theta}$.

(i) If $0 < \theta < 1$, then
$$\left(L_{p_0,r_0,b_0}, L_{b_1}^{\infty),q_1,r_1}\right)_{\theta,r,\mathrm{a}} = L_{p,r,a^\#}.$$

(ii) If $\left\|s^{-1/r}a(s)\right\|_{r,(1,\infty)} < \infty$, then
$$\left(L_{p_0,r_0,b_0}, L_{b_1}^{\infty),q_1,r_1}\right)_{0,r,\mathrm{a}} = L_{p_0,r,a^\circ\rho,r_0,b_0}^{\mathcal{L}}.$$

(iii) If $\left\|s^{-1/r}a(s)\right\|_{r,(0,1)} < \infty$, then
$$\left(L_{p_0,r_0,b_0}, L_{b_1}^{\infty),q_1,r_1}\right)_{1,r,\mathrm{a}} = L_{\infty,r,a^\#} \cap L_{B_1}^{\infty),q_1,r} \cap L_{\infty,(r,a^\circ\rho,r_1,b_1,q_1,1)}^{\mathcal{R},\mathcal{R}},$$

where $B_1(t) = \left\|s^{-1/r_1}b_1(s)\right\|_{r_1,(0,t)} a(\rho(t))$.

*Proof.* Let $\theta_0 = 1 - \frac{1}{p_0}$. The proof follows from Theorem 44 and the formulae $(L_1, L_\infty)_{\theta_0, r_0, b_0} = L_{p_0, r_0, b_0}$, $(L_1, L_\infty)_{\theta_0, r, a^\circ \rho, r_0, b_0}^{\mathcal{L}} = L_{p_0, r, a^\circ \rho, r_0, b_0}^{\mathcal{L}}$, $(L_1, L_\infty)_{1, r, B_1, q_1, 1}^{\mathcal{R}} = L_{B_1}^{\infty), q_1, r}$, and $(L_1, L_\infty)_{1, r, a^\circ \rho, r_1, b_1, q_1, 1}^{\mathcal{R}, \mathcal{R}} = L_{\infty, (r, a^\circ \rho, r_1, b_1, q_1, 1)}^{\mathcal{R}, \mathcal{R}}$. □

Taking $q_1 = \infty$ in the previous corollary, we get the next one. See the proof of Corollary 61.

**Corollary** 63. (Cf. [15, Corollary 6.3].) Let $1 < p_0 < \infty$, $0 < r, r_0, r_1 \leq \infty$, $a, b_0, b_1 \in SV$, and $\|s^{-1/r_1} b_1(s)\|_{r_1, (0,1)} < \infty$. Put
$$\chi_1(t) = \|s^{-1/r_1} b_1(s)\|_{r_1, (0, t)}$$
and $\rho(t) = t^{\frac{1}{p_0}} \frac{b_0(t)}{\chi_1(t)}$. For a given $0 \leq \theta \leq 1$, denote $a^\# = b_0^{1-\theta} \chi_1^\theta a^\circ \rho$ and $p = \frac{p_0}{1-\theta}$.

(i) If $0 < \theta < 1$, then
$$(L_{p_0, r_0, b_0}, L_{\infty, r_1, b_1})_{\theta, r, a} = L_{p, r, a^\#}.$$

(ii) If $\|s^{-1/r} a(s)\|_{r, (1, \infty)} < \infty$, then
$$(L_{p_0, r_0, b_0}, L_{\infty, r_1, b_1})_{0, r, a} = L_{p_0, r, a^\circ \rho, r_0, b_0}^{\mathcal{L}}.$$

(iii) If $\|s^{-1/r} a(s)\|_{r, (0, 1)} < \infty$, then
$$(L_{p_0, r_0, b_0}, L_{\infty, r_1, b_1})_{1, r, a} = L_{\infty, r, a^\#} \cap L_{\infty, r, a^\circ \rho, r_1, b_1}^{\mathcal{R}}.$$

### 6.3. Some interpolation results for couples, where the first operand is a grand or a small Lorentz space

**Corollary** 64. (Cf. [8, Corollary 52; 9, Corollary 65; 16, Corollary 5.7; 17, Corollary 5.5].) Let $1 < p_0 < \infty$, $0 < r, r_0, q_1, r_1 \leq \infty$, $a, b_0, b_1 \in SV$, $\|s^{-1/r_0} b_0(s)\|_{r_0, (0, 1)} < \infty$, and $\|s^{-1/r_1} b_1(s) |\ln s|^{1/q_1}\|_{r_1, (0, 1)} < \infty$. Put
$$\chi_0(t) = \|s^{-1/r_0} b_0(s)\|_{r_0, (0, t)}, \qquad \chi_1(t) = \left\|s^{-1/r_1} b_1(s) \left\{\ln \frac{t}{s}\right\}^{1/q_1}\right\|_{r_1, (0, t)},$$
and $\rho(t) = t^{\frac{1}{p_0}} \frac{\chi_0(t)}{\chi_1(t)}$. For a given $0 \leq \theta \leq 1$, denote $a^\# = \chi_0^{1-\theta} \chi_1^\theta a^\circ \rho$ and $p = \frac{p_0}{1-\theta}$.

(i) If $0 < \theta < 1$, then
$$\left(L_{b_0}^{p_0), q_0, r_0}, L_{b_1}^{\infty), q_1, r_1}\right)_{\theta, r, a} = L_{p, r, a^\#}.$$

(ii) If $\|s^{-1/r} a(s)\|_{r, (1, \infty)} < \infty$, then
$$\left(L_{b_0}^{p_0), q_0, r_0}, L_{b_1}^{\infty), q_1, r_1}\right)_{0, r, a} = L_{p_0, (r, a^\circ \rho, r_0, b_0, q_0, 1)}^{\mathcal{R}, \mathcal{L}}.$$

(iii) If $\|s^{-1/r} a(s)\|_{r, (0, 1)} < \infty$, then
$$\left(L_{b_0}^{p_0), q_0, r_0}, L_{b_1}^{\infty), q_1, r_1}\right)_{1, r, a} = L_{\infty, r, a^\#} \cap L_{B_1}^{\infty), q_1, r_1} \cap L_{\infty, (r, a^\circ \rho, r_1, b_1, q_1, 1)}^{\mathcal{R}, \mathcal{R}},$$
where $B_1(t) = \|s^{-1/r_1} b_1(s)\|_{r_1, (0, t)} a(\rho(t))$.

*Proof.* Put $\theta_0 = 1 - \frac{1}{p_0}$. The proof follows from Theorem 50 if we take $a_0 = 1$ and use the formulae $(L_1, L_\infty)_{\theta_0, r_0, b_0, q_0, 1}^{\mathcal{R}} = L_{b_0}^{p_0), q_0, r_0}$, $(L_1, L_\infty)_{1, r, B_1, q_1, 1}^{\mathcal{R}} = L_{B_1}^{\infty), q_1, r}$, $(L_1, L_\infty)_{\theta_0, r, a^\circ \rho, r_0, b_0, q_0, 1}^{\mathcal{R}, \mathcal{L}} = L_{p_0, (r, a^\circ \rho, r_0, b_0, q_0, 1)}^{\mathcal{R}, \mathcal{L}}$, and $(L_1, L_\infty)_{1, r, a^\circ \rho, r_1, b_1, q_1, 1}^{\mathcal{R}, \mathcal{R}} = L_{\infty, (r, a^\circ \rho, r_1, b_1, q_1, 1)}^{\mathcal{R}, \mathcal{R}}$. □

Similarly, using Theorem 48, the next corollary can be proved.

**Corollary** 65. (Cf. [8, Corollary 51; 9, Corollary 66; 17, Theorem 5.7].) Let $1 < p_0 < \infty$, $0 < r, r_0, q_1, r_1 \leq \infty$, $a, b_0, b_1 \in SV$, $\left\|s^{-1/r_0} b_0(s)\right\|_{r_0,(1,\infty)} < \infty$, and $\left\|s^{-1/r_1} b_1(s) |\ln s|^{1/q_1}\right\|_{r_1,(0,1)} < \infty$. Put

$$\chi_0(t) = \left\|s^{-1/r_0} b_0(s)\right\|_{r_0,(t,\infty)}, \qquad \chi_1(t) = \left\|s^{-1/r_1} b_1(s) \left\{\ln \frac{t}{s}\right\}^{1/q_1}\right\|_{r_1,(0,t)},$$

and $\rho(t) = t^{\frac{1}{p_0}} \frac{\chi_0(t)}{\chi_1(t)}$. For a given $0 \leq \theta \leq 1$, denote $a^\# = \chi_0^{1-\theta} \chi_1^\theta a \circ \rho$ and $p = \frac{p_0}{1-\theta}$.

(i) If $0 < \theta < 1$, then

$$\left(L_{b_0}^{(p_0,q_0,r_0)}, L_{b_1}^{\infty),q_1,r_1}\right)_{\theta,r,a} = L_{p,r,a^\#}.$$

(ii) If $\left\|s^{-1/r} a(s)\right\|_{r,(1,\infty)} < \infty$, then

$$\left(L_{b_0}^{(p_0,q_0,r_0)}, L_{b_1}^{\infty),q_1,r_1}\right)_{0,r,a} = L_{a^\#}^{(p_0,q_0,r_0} \cap L_{p_0,(r,a\circ\rho,r_0,b_0,q_0,1)}^{\mathcal{L},\mathcal{L}}.$$

(iii) If $\left\|s^{-1/r} a(s)\right\|_{r,(0,1)} < \infty$, then

$$\left(L_{b_0}^{(p_0,q_0,r_0)}, L_{b_1}^{\infty),q_1,r_1}\right)_{1,r,a} = L_{\infty,r,a^\#} \cap L_{B_1}^{\infty),q_1,r_1} \cap L_{\infty,(r,a\circ\rho,r_1,b_1,q_1,1)}^{\mathcal{R},\mathcal{R}},$$

where $B_1(t) = \left\|s^{-1/r_1} b_1(s)\right\|_{r_1,(0,t)} a(\rho(t))$.

(Leo R. Ya. Doktorski) Department Object Recognition, Fraunhofer Institute of Optronics, System Technologies and Image Exploitation IOSB, Gutleuthausstr. 1, 76275 Ettlingen, Germany, doktorskileo@gmail.com.

(Pedro Fernández-Martínez and Teresa M. Signes) Departamento de Matemáticas, Facultad de Matemáticas, Universidad de Murcia, Campus de Espinardo, 30071 Espinardo (Murcia), Spain.